
\input amstex

\expandafter\ifx\csname beta.def\endcsname\relax \else\endinput\fi
\expandafter\edef\csname beta.def\endcsname{%
 \catcode`\noexpand\@=\the\catcode`\@\space}

\let\atbefore @

\catcode`\@=11

\overfullrule\z@
\hsize 6.25truein
\vsize 9.63truein

\let\@ft@\expandafter \let\@tb@f@\atbefore

\newif\ifMag
\ifnum\mag>1000 \Magtrue\fi
=\ifMag cmr8\else cmr9\fi

\newdimen\p@@ \p@@\p@
\def\m@ths@r{\ifnum\mathsurround=\z@\z@\else\maths@r\fi}
\def\maths@r{1.6\p@@} \def\mathsurzero{\def\maths@r{\z@}}

\mathsurround\maths@r
\font\Brm=cmr12 \font\Bbf=cmbx12 \font\Bit=cmti12 \font\ssf=cmss10
\font\Bsl=cmsl10 scaled 1200 \font\Bmmi=cmmi10 scaled 1200
\font\BBf=cmbx12 scaled 1200 \font\BMmi=cmmi10 scaled 1440

\def\atletter{\edef\atrestore{\catcode`\noexpand\@=\the\catcode`\@\space}
 \catcode`\@=11}

\newread\@ux \newwrite\@@x \newwrite\@@cd
\let\@np@@\input
\def\@np@t#1{\openin\@ux#1\relax\ifeof\@ux\else\closein\@ux\relax\@np@@ #1\fi}
\def\input#1 {\openin\@ux#1\relax\ifeof\@ux\wrs@x{No file #1}\else
 \closein\@ux\relax\@np@@ #1\fi}
\def\Input#1 {\relax} 

\def\wr@@x#1{} \def\wrs@x{\immediate\write\sixt@@n}

\def\readldf{\@np@t{\jobname.ldf}}
\def\writeldf{\def\wr@@x{\immediate\write\@@x}
 \def\cl@selbl{\wr@@x{\string\Snodef{\the\Sno}}\wr@@x{\string\endinput}%
 \immediate\closeout\@@x} \immediate\openout\@@x\jobname.ldf}
\let\cl@selbl\relax

\def\tod@y{\ifcase\month\or
 January\or February\or March\or April\or May\or June\or July\or
 August\or September\or October\or November\or December\fi\space\,
\number\day,\space\,\number\year}

\newcount\c@time
\def\h@@r{hh}\def\m@n@te{mm}
\def\wh@tt@me{\c@time\time\divide\c@time 60\edef\h@@r{\number\c@time}%
 \multiply\c@time -60\advance\c@time\time\edef
 \m@n@te{\ifnum\c@time<10 0\fi\number\c@time}}
\def\t@me{\h@@r\/{\rm:}\m@n@te}  \let\whattime\wh@tt@me
\def\today{\tod@y\wr@@x{\string\todaydef{\tod@y}}}
\def\nowtime{\t@me{\let\/\ic@\wr@@x{\string\nowtimedef{\t@me}}}}
\def\todaydef#1{} \def\nowtimedef#1{}

\def\em#1{{\it #1\/}} \def\emph#1{{\sl #1\/}}

\def\itemleft#1{\par\setbox\z@\hbox{\rm #1\enspace}\hangindent\wd\z@
 \hglue-2\parindent\kern\wd\z@\textindent{\rm#1}}
\def\itemflat#1{\par\setbox\z@\hbox{\rm #1\enspace}\hang\ifnum\wd\z@>\parindent
 \noindent\unhbox\z@\ignore\else\textindent{\rm#1}\fi}

\newcount\itemlet
\def\newbi{\itemlet 96} \newbi
\def\bitem{\gad\itemlet \par\hangindent1.5\parindent
 \hglue-.5\parindent\textindent{\rm\rlap{\char\the\itemlet}\hp{b})}}

\newcount\itemrm

\def\iitem{\gad\itemrm \par\hangindent1.5\parindent
 \hglue-.5\parindent\textindent{\rm\hp{v}\llap{\romannumeral\the\itemrm})}}

\def\center{\par\begingroup\leftskip\z@ plus \hsize \rightskip\leftskip
 \parindent\z@\parfillskip\z@skip \def\\{\unskip\break}}
\def\endcenter{\endgraf\endgroup}

\let\b@gr@@\begingroup \let\B@gr@@\begingroup
\def\b@gr@{\b@gr@@\let\b@gr@@\undefined}
\def\B@gr@{\B@gr@@\let\B@gr@@\undefined}

\def\@fn@xt#1#2#3{\let\@ch@r=#1\def\n@xt{\ifx\t@st@\@ch@r
 \def\n@@xt{#2}\else\def\n@@xt{#3}\fi\n@@xt}\futurelet\t@st@\n@xt}

\def\@fwd@@#1#2#3{\setbox\z@\hbox{#1}\ifdim\wd\z@>\z@#2\else#3\fi}
\def\s@twd@#1#2{\setbox\z@\hbox{#2}#1\wd\z@}

\def\r@st@re#1{\let#1\s@v@} \def\s@v@d@f{\let\s@v@}

\def\p@sk@p#1#2{\par\skip@#2\relax\ifdim\lastskip<\skip@\relax\removelastskip
 \ifnum#1=\z@\else\penalty#1\relax\fi\vskip\skip@
 \else\ifnum#1=\z@\else\penalty#1\relax\fi\fi}
\def\sk@@p#1{\par\skip@#1\relax\ifdim\lastskip<\skip@\relax\removelastskip
 \vskip\skip@\fi}

\newbox\p@b@ld
\def\poorbold#1{\setbox\p@b@ld\hbox{#1}\kern-.01em\copy\p@b@ld\kern-\wd\p@b@ld
 \kern.02em\copy\p@b@ld\kern-\wd\p@b@ld\kern-.012em\raise.02em\box\p@b@ld}

\ifx\plainfootnote\undefined \let\plainfootnote\footnote \fi

\let\s@v@\proclaim \let\proclaim\relax
\def\r@R@fs#1{\let#1\s@R@fs} \let\s@R@fs\Refs \let\Refs\relax
\def\r@endd@#1{\let#1\s@endd@} \let\s@endd@\enddocument
\let\bye\relax

\def\myR@fs{\@fn@xt[\m@R@f@\m@R@fs} \def\m@R@fs{\@fn@xt*\m@r@f@@\m@R@f@@}
\def\m@R@f@@{\m@R@f@[References]} \def\m@r@f@@*{\m@R@f@[]}

\def\Twelvepoint{\twelvepoint \let\Bbf\BBf \let\Bmmi\BMmi
\font\Brm=cmr12 scaled 1200 \font\Bit=cmti12 scaled 1200
\font\ssf=cmss10 scaled 1200 \font\Bsl=cmsl10 scaled 1440
\font\BBf=cmbx12 scaled 1440 \font\BMmi=cmmi10 scaled 1728}

\newif\ifamsppt

\newdimen\b@gsize

\newdimen\r@f@nd \newbox\r@f@b@x \newbox\adjb@x
\newbox\p@nct@ \newbox\k@yb@x \newcount\rcount
\newbox\b@b@x \newbox\p@p@rb@x \newbox\j@@rb@x \newbox\y@@rb@x
\newbox\v@lb@x \newbox\is@b@x \newbox\p@g@b@x \newif\ifp@g@ \newif\ifp@g@s
\newbox\inb@@kb@x \newbox\b@@kb@x \newbox\p@blb@x \newbox\p@bl@db@x
\newbox\ed@b@x \newif\ifed@ \newif\ifed@s \newif\if@fl@b \newif\if@fn@m
\newbox\p@p@nf@b@x \newbox\inf@b@x \newbox\b@@nf@b@x

\newif\ifp@gen@

\@ft@\ifx\csname amsppt.sty\endcsname\relax

\headline={\hfil}
\footline={\ifp@gen@\ifnum\pageno=\z@\else\hfil\foliorm\folio\fi\else
 \ifnum\pageno=\z@\hfil\foliorm\folio\fi\fi\hfil\global\p@gen@true}
\parindent1pc

\font@\tensmc=cmcsc10
\font@\sevenex=cmex7
\font@\sevenit=cmti7
\font@\eightrm=cmr8
\font@\sixrm=cmr6
\font@\eighti=cmmi8 \skewchar\eighti='177
\font@\sixi=cmmi6 \skewchar\sixi='177
\font@\eightsy=cmsy8 \skewchar\eightsy='60
\font@\sixsy=cmsy6 \skewchar\sixsy='60
\font@\eightex=cmex8
\font@\eightbf=cmbx8
\font@\sixbf=cmbx6
\font@\eightit=cmti8
\font@\eightsl=cmsl8
\font@\eightsmc=cmcsc8
\font@\eighttt=cmtt8
\font@\ninerm=cmr9
\font@\ninei=cmmi9 \skewchar\ninei='177
\font@\ninesy=cmsy9 \skewchar\ninesy='60
\font@\nineex=cmex9
\font@\ninebf=cmbx9
\font@\nineit=cmti9
\font@\ninesl=cmsl9
\font@\ninesmc=cmcsc9
\font@\ninemsa=msam9
\font@\ninemsb=msbm9
\font@\nineeufm=eufm9
\font@\eightmsa=msam8
\font@\eightmsb=msbm8
\font@\eighteufm=eufm8
\font@\sixmsa=msam6
\font@\sixmsb=msbm6
\font@\sixeufm=eufm6

\loadmsam\loadmsbm\loadeufm
\input amssym.tex

\def\footnoterule{\kern-3\p@\hrule width5pc\kern 2.6\p@}
\def\m@k@foot#1{\insert\footins
 {\interlinepenalty\interfootnotelinepenalty
 \eightpoint\splittopskip\ht\strutbox\splitmaxdepth\dp\strutbox
 \floatingpenalty\@MM\leftskip\z@\rightskip\z@
 \spaceskip\z@\xspaceskip\z@
 \leavevmode\footstrut\ignore#1\unskip\lower\dp\strutbox
 \vbox to\dp\strutbox{}}}
\def\ftext#1{\m@k@foot{\vsk-.8>\nt #1}}
\def\pr@cl@@m#1{\p@sk@p{-100}\medskipamount\b@gr@\nt\ignore
 \bf #1\unskip.\enspace\sl\ignore}
\outer\def\proclaim{\pr@cl@@m} \s@v@d@f\proclaim \let\proclaim\relax
\def\endproclaim{\endgroup\p@sk@p{55}\medskipamount}
\def\demo#1{\sk@@p\medskipamount\nt{\ignore\it #1\unskip.}\enspace
 \ignore}
\def\enddemo{\sk@@p\medskipamount}

\def\cite#1{{\rm[#1]}} \let\nofrills\relax
\let\NoRunningHeads\relax \def\Refs#1#2{\relax}

\def\big@#1#2{{\hbox{$\left#2\vcenter to#1\b@gsize{}%
 \right.\nulldelimiterspace\z@\m@th$}}}
\def\big{\big@\@ne}
\def\Big{\big@{1.5}}
\def\bigg{\big@\tw@}
\def\Bigg{\big@{2.5}}
\normallineskiplimit\p@

\def\tenpoint{\p@@\p@ \normallineskiplimit\p@@
 \mathsurround\m@ths@r \normalbaselineskip12\p@@
 \abovedisplayskip12\p@@ plus3\p@@ minus9\p@@
 \belowdisplayskip\abovedisplayskip
 \abovedisplayshortskip\z@ plus3\p@@
 \belowdisplayshortskip7\p@@ plus3\p@@ minus4\p@@
 \textonlyfont@\rm\tenrm \textonlyfont@\it\tenit
 \textonlyfont@\sl\tensl \textonlyfont@\bf\tenbf
 \textonlyfont@\smc\tensmc \textonlyfont@\tt\tentt
 \ifsyntax@ \def\big##1{{\hbox{$\left##1\right.$}}}%
  \let\Big\big \let\bigg\big \let\Bigg\big
 \else
  \textfont\z@\tenrm \scriptfont\z@\sevenrm \scriptscriptfont\z@\fiverm
  \textfont\@ne\teni \scriptfont\@ne\seveni \scriptscriptfont\@ne\fivei
  \textfont\tw@\tensy \scriptfont\tw@\sevensy \scriptscriptfont\tw@\fivesy
  \textfont\thr@@\tenex \scriptfont\thr@@\sevenex
	\scriptscriptfont\thr@@\sevenex
  \textfont\itfam\tenit \scriptfont\itfam\sevenit
	\scriptscriptfont\itfam\sevenit
  \textfont\bffam\tenbf \scriptfont\bffam\sevenbf
	\scriptscriptfont\bffam\fivebf
  \textfont\msafam\tenmsa \scriptfont\msafam\sevenmsa
	\scriptscriptfont\msafam\fivemsa
  \textfont\msbfam\tenmsb \scriptfont\msbfam\sevenmsb
	\scriptscriptfont\msbfam\fivemsb
  \textfont\eufmfam\teneufm \scriptfont\eufmfam\seveneufm
	\scriptscriptfont\eufmfam\fiveeufm
  \setbox\strutbox\hbox{\vrule height8.5\p@@ depth3.5\p@@ width\z@}%
  \setbox\strutbox@\hbox{\lower.5\normallineskiplimit\vbox{%
	\kern-\normallineskiplimit\copy\strutbox}}%
   \setbox\z@\vbox{\hbox{$($}\kern\z@}\b@gsize1.2\ht\z@
  \fi
  \normalbaselines\rm\dotsspace@1.5mu\ex@.2326ex\jot3\ex@}

\def\eightpoint{\p@@.8\p@ \normallineskiplimit\p@@
 \mathsurround\m@ths@r \normalbaselineskip10\p@
 \abovedisplayskip10\p@ plus2.4\p@ minus7.2\p@
 \belowdisplayskip\abovedisplayskip
 \abovedisplayshortskip\z@ plus3\p@@
 \belowdisplayshortskip7\p@@ plus3\p@@ minus4\p@@
 \textonlyfont@\rm\eightrm \textonlyfont@\it\eightit
 \textonlyfont@\sl\eightsl \textonlyfont@\bf\eightbf
 \textonlyfont@\smc\eightsmc \textonlyfont@\tt\eighttt
 \ifsyntax@\def\big##1{{\hbox{$\left##1\right.$}}}%
  \let\Big\big \let\bigg\big \let\Bigg\big
 \else
  \textfont\z@\eightrm \scriptfont\z@\sixrm \scriptscriptfont\z@\fiverm
  \textfont\@ne\eighti \scriptfont\@ne\sixi \scriptscriptfont\@ne\fivei
  \textfont\tw@\eightsy \scriptfont\tw@\sixsy \scriptscriptfont\tw@\fivesy
  \textfont\thr@@\eightex \scriptfont\thr@@\sevenex
	\scriptscriptfont\thr@@\sevenex
  \textfont\itfam\eightit \scriptfont\itfam\sevenit
	\scriptscriptfont\itfam\sevenit
  \textfont\bffam\eightbf \scriptfont\bffam\sixbf
	\scriptscriptfont\bffam\fivebf
  \textfont\msafam\eightmsa \scriptfont\msafam\sixmsa
	\scriptscriptfont\msafam\fivemsa
  \textfont\msbfam\eightmsb \scriptfont\msbfam\sixmsb
	\scriptscriptfont\msbfam\fivemsb
  \textfont\eufmfam\eighteufm \scriptfont\eufmfam\sixeufm
	\scriptscriptfont\eufmfam\fiveeufm
 \setbox\strutbox\hbox{\vrule height7\p@ depth3\p@ width\z@}%
 \setbox\strutbox@\hbox{\raise.5\normallineskiplimit\vbox{%
   \kern-\normallineskiplimit\copy\strutbox}}%
 \setbox\z@\vbox{\hbox{$($}\kern\z@}\b@gsize1.2\ht\z@
 \fi
 \normalbaselines\eightrm\dotsspace@1.5mu\ex@.2326ex\jot3\ex@}

\def\ninepoint{\p@@.9\p@ \normallineskiplimit\p@@
 \mathsurround\m@ths@r \normalbaselineskip11\p@
 \abovedisplayskip11\p@ plus2.7\p@ minus8.1\p@
 \belowdisplayskip\abovedisplayskip
 \abovedisplayshortskip\z@ plus3\p@@
 \belowdisplayshortskip7\p@@ plus3\p@@ minus4\p@@
 \textonlyfont@\rm\ninerm \textonlyfont@\it\nineit
 \textonlyfont@\sl\ninesl \textonlyfont@\bf\ninebf
 \textonlyfont@\smc\ninesmc \textonlyfont@\tt\ninett
 \ifsyntax@ \def\big##1{{\hbox{$\left##1\right.$}}}%
  \let\Big\big \let\bigg\big \let\Bigg\big
 \else
  \textfont\z@\ninerm \scriptfont\z@\sevenrm \scriptscriptfont\z@\fiverm
  \textfont\@ne\ninei \scriptfont\@ne\seveni \scriptscriptfont\@ne\fivei
  \textfont\tw@\ninesy \scriptfont\tw@\sevensy \scriptscriptfont\tw@\fivesy
  \textfont\thr@@\nineex \scriptfont\thr@@\sevenex
	\scriptscriptfont\thr@@\sevenex
  \textfont\itfam\nineit \scriptfont\itfam\sevenit
	\scriptscriptfont\itfam\sevenit
  \textfont\bffam\ninebf \scriptfont\bffam\sevenbf
	\scriptscriptfont\bffam\fivebf
  \textfont\msafam\ninemsa \scriptfont\msafam\sevenmsa
	\scriptscriptfont\msafam\fivemsa
  \textfont\msbfam\ninemsb \scriptfont\msbfam\sevenmsb
	\scriptscriptfont\msbfam\fivemsb
  \textfont\eufmfam\nineeufm \scriptfont\eufmfam\seveneufm
	\scriptscriptfont\eufmfam\fiveeufm
  \setbox\strutbox\hbox{\vrule height8.5\p@@ depth3.5\p@@ width\z@}%
  \setbox\strutbox@\hbox{\lower.5\normallineskiplimit\vbox{%
	\kern-\normallineskiplimit\copy\strutbox}}%
   \setbox\z@\vbox{\hbox{$($}\kern\z@}\b@gsize1.2\ht\z@
  \fi
  \normalbaselines\rm\dotsspace@1.5mu\ex@.2326ex\jot3\ex@}

\font@\twelverm=cmr10 scaled 1200
\font@\twelveit=cmti10 scaled 1200
\font@\twelvesl=cmsl10 scaled 1200
\font@\twelvebf=cmbx10 scaled 1200
\font@\twelvesmc=cmcsc10 scaled 1200
\font@\twelvett=cmtt10 scaled 1200
\font@\twelvei=cmmi10 scaled 1200 \skewchar\twelvei='177
\font@\twelvesy=cmsy10 scaled 1200 \skewchar\twelvesy='60
\font@\twelveex=cmex10 scaled 1200
\font@\twelvemsa=msam10 scaled 1200
\font@\twelvemsb=msbm10 scaled 1200
\font@\twelveeufm=eufm10 scaled 1200

\def\twelvepoint{\p@@1.2\p@ \normallineskiplimit\p@@
 \mathsurround\m@ths@r \normalbaselineskip12\p@@
 \abovedisplayskip12\p@@ plus3\p@@ minus9\p@@
 \belowdisplayskip\abovedisplayskip
 \abovedisplayshortskip\z@ plus3\p@@
 \belowdisplayshortskip7\p@@ plus3\p@@ minus4\p@@
 \textonlyfont@\rm\twelverm \textonlyfont@\it\twelveit
 \textonlyfont@\sl\twelvesl \textonlyfont@\bf\twelvebf
 \textonlyfont@\smc\twelvesmc \textonlyfont@\tt\twelvett
 \ifsyntax@ \def\big##1{{\hbox{$\left##1\right.$}}}%
  \let\Big\big \let\bigg\big \let\Bigg\big
 \else
  \textfont\z@\twelverm \scriptfont\z@\eightrm \scriptscriptfont\z@\sixrm
  \textfont\@ne\twelvei \scriptfont\@ne\eighti \scriptscriptfont\@ne\sixi
  \textfont\tw@\twelvesy \scriptfont\tw@\eightsy \scriptscriptfont\tw@\sixsy
  \textfont\thr@@\twelveex \scriptfont\thr@@\eightex
	\scriptscriptfont\thr@@\sevenex
  \textfont\itfam\twelveit \scriptfont\itfam\eightit
	\scriptscriptfont\itfam\sevenit
  \textfont\bffam\twelvebf \scriptfont\bffam\eightbf
	\scriptscriptfont\bffam\sixbf
  \textfont\msafam\twelvemsa \scriptfont\msafam\eightmsa
	\scriptscriptfont\msafam\sixmsa
  \textfont\msbfam\twelvemsb \scriptfont\msbfam\eightmsb
	\scriptscriptfont\msbfam\sixmsb
  \textfont\eufmfam\twelveeufm \scriptfont\eufmfam\eighteufm
	\scriptscriptfont\eufmfam\sixeufm
  \setbox\strutbox\hbox{\vrule height8.5\p@@ depth3.5\p@@ width\z@}%
  \setbox\strutbox@\hbox{\lower.5\normallineskiplimit\vbox{%
	\kern-\normallineskiplimit\copy\strutbox}}%
  \setbox\z@\vbox{\hbox{$($}\kern\z@}\b@gsize1.2\ht\z@
  \fi
  \normalbaselines\rm\dotsspace@1.5mu\ex@.2326ex\jot3\ex@}

\font@\twelvetrm=cmr10 at 12truept
\font@\twelvetit=cmti10 at 12truept
\font@\twelvetsl=cmsl10 at 12truept
\font@\twelvetbf=cmbx10 at 12truept
\font@\twelvetsmc=cmcsc10 at 12truept
\font@\twelvettt=cmtt10 at 12truept
\font@\twelveti=cmmi10 at 12truept \skewchar\twelveti='177
\font@\twelvetsy=cmsy10 at 12truept \skewchar\twelvetsy='60
\font@\twelvetex=cmex10 at 12truept
\font@\twelvetmsa=msam10 at 12truept
\font@\twelvetmsb=msbm10 at 12truept
\font@\twelveteufm=eufm10 at 12truept

\def\twelvetruepoint{\p@@1.2truept \normallineskiplimit\p@@
 \mathsurround\m@ths@r \normalbaselineskip12\p@@
 \abovedisplayskip12\p@@ plus3\p@@ minus9\p@@
 \belowdisplayskip\abovedisplayskip
 \abovedisplayshortskip\z@ plus3\p@@
 \belowdisplayshortskip7\p@@ plus3\p@@ minus4\p@@
 \textonlyfont@\rm\twelvetrm \textonlyfont@\it\twelvetit
 \textonlyfont@\sl\twelvetsl \textonlyfont@\bf\twelvetbf
 \textonlyfont@\smc\twelvetsmc \textonlyfont@\tt\twelvettt
 \ifsyntax@ \def\big##1{{\hbox{$\left##1\right.$}}}%
  \let\Big\big \let\bigg\big \let\Bigg\big
 \else
  \textfont\z@\twelvetrm \scriptfont\z@\eightrm \scriptscriptfont\z@\sixrm
  \textfont\@ne\twelveti \scriptfont\@ne\eighti \scriptscriptfont\@ne\sixi
  \textfont\tw@\twelvetsy \scriptfont\tw@\eightsy \scriptscriptfont\tw@\sixsy
  \textfont\thr@@\twelvetex \scriptfont\thr@@\eightex
	\scriptscriptfont\thr@@\sevenex
  \textfont\itfam\twelvetit \scriptfont\itfam\eightit
	\scriptscriptfont\itfam\sevenit
  \textfont\bffam\twelvetbf \scriptfont\bffam\eightbf
	\scriptscriptfont\bffam\sixbf
  \textfont\msafam\twelvetmsa \scriptfont\msafam\eightmsa
	\scriptscriptfont\msafam\sixmsa
  \textfont\msbfam\twelvetmsb \scriptfont\msbfam\eightmsb
	\scriptscriptfont\msbfam\sixmsb
  \textfont\eufmfam\twelveteufm \scriptfont\eufmfam\eighteufm
	\scriptscriptfont\eufmfam\sixeufm
  \setbox\strutbox\hbox{\vrule height8.5\p@@ depth3.5\p@@ width\z@}%
  \setbox\strutbox@\hbox{\lower.5\normallineskiplimit\vbox{%
	\kern-\normallineskiplimit\copy\strutbox}}%
  \setbox\z@\vbox{\hbox{$($}\kern\z@}\b@gsize1.2\ht\z@
  \fi
  \normalbaselines\rm\dotsspace@1.5mu\ex@.2326ex\jot3\ex@}

\font@\elevenrm=cmr10 scaled 1095
\font@\elevenit=cmti10 scaled 1095
\font@\elevensl=cmsl10 scaled 1095
\font@\elevenbf=cmbx10 scaled 1095
\font@\elevensmc=cmcsc10 scaled 1095
\font@\eleventt=cmtt10 scaled 1095
\font@\eleveni=cmmi10 scaled 1095 \skewchar\eleveni='177
\font@\elevensy=cmsy10 scaled 1095 \skewchar\elevensy='60
\font@\elevenex=cmex10 scaled 1095
\font@\elevenmsa=msam10 scaled 1095
\font@\elevenmsb=msbm10 scaled 1095
\font@\eleveneufm=eufm10 scaled 1095

\def\elevenpoint{\p@@1.1\p@ \normallineskiplimit\p@@
 \mathsurround\m@ths@r \normalbaselineskip12\p@@
 \abovedisplayskip12\p@@ plus3\p@@ minus9\p@@
 \belowdisplayskip\abovedisplayskip
 \abovedisplayshortskip\z@ plus3\p@@
 \belowdisplayshortskip7\p@@ plus3\p@@ minus4\p@@
 \textonlyfont@\rm\elevenrm \textonlyfont@\it\elevenit
 \textonlyfont@\sl\elevensl \textonlyfont@\bf\elevenbf
 \textonlyfont@\smc\elevensmc \textonlyfont@\tt\eleventt
 \ifsyntax@ \def\big##1{{\hbox{$\left##1\right.$}}}%
  \let\Big\big \let\bigg\big \let\Bigg\big
 \else
  \textfont\z@\elevenrm \scriptfont\z@\eightrm \scriptscriptfont\z@\sixrm
  \textfont\@ne\eleveni \scriptfont\@ne\eighti \scriptscriptfont\@ne\sixi
  \textfont\tw@\elevensy \scriptfont\tw@\eightsy \scriptscriptfont\tw@\sixsy
  \textfont\thr@@\elevenex \scriptfont\thr@@\eightex
	\scriptscriptfont\thr@@\sevenex
  \textfont\itfam\elevenit \scriptfont\itfam\eightit
	\scriptscriptfont\itfam\sevenit
  \textfont\bffam\elevenbf \scriptfont\bffam\eightbf
	\scriptscriptfont\bffam\sixbf
  \textfont\msafam\elevenmsa \scriptfont\msafam\eightmsa
	\scriptscriptfont\msafam\sixmsa
  \textfont\msbfam\elevenmsb \scriptfont\msbfam\eightmsb
	\scriptscriptfont\msbfam\sixmsb
  \textfont\eufmfam\eleveneufm \scriptfont\eufmfam\eighteufm
	\scriptscriptfont\eufmfam\sixeufm
  \setbox\strutbox\hbox{\vrule height8.5\p@@ depth3.5\p@@ width\z@}%
  \setbox\strutbox@\hbox{\lower.5\normallineskiplimit\vbox{%
	\kern-\normallineskiplimit\copy\strutbox}}%
  \setbox\z@\vbox{\hbox{$($}\kern\z@}\b@gsize1.2\ht\z@
  \fi
  \normalbaselines\rm\dotsspace@1.5mu\ex@.2326ex\jot3\ex@}

\def\m@R@f@[#1]{\mathsurzero{
 \s@ct{}{#1}}\wr@@c{\string\Refcd{#1}{\the\pageno}}\B@gr@
 \frenchspacing\rcount\z@\refkey{[##1]}\refno{[##1]}\widest{AZ}\keyright
 \let\Key\key\let\refin\relax}
\def\widest#1{\s@twd@\r@f@nd{\r@fk@y{#1}\enspace}}
\def\widestno#1{\s@twd@\r@f@nd{\r@fn@{#1}\enspace}}
\def\widestlabel#1{\s@twd@\r@f@nd{#1\enspace}}
\def\refkey{\def\r@fk@y##1} \def\refno{\def\r@fn@##1}
\def\keyright{\def\r@fit@m{\hang\textindent}}
\def\keyflat{\def\r@fit@m##1{\setbox\z@\hbox{\rm ##1\enspace}\hang\noindent
 \ifnum\wd\z@<\parindent\indent\hglue-\wd\z@\fi\unhbox\z@}}

\def\R@fb@x{\global\setbox\r@f@b@x} \def\K@yb@x{\global\setbox\k@yb@x}
\def\ref{\par\b@gr@\rm\R@fb@x\box\voidb@x\K@yb@x\box\voidb@x\@fn@mfalse
 \@fl@bfalse\b@g@nr@f}
\def\c@nc@t#1{\setbox\z@\lastbox
 \setbox\adjb@x\hbox{\unhbox\adjb@x\unhbox\z@\unskip\unskip\unpenalty#1}}
\def\adjust#1{\relax\ifmmode\penalty-\@M\null\hfil$\clubpenalty\z@
 \widowpenalty\z@\interlinepenalty\z@\offinterlineskip\endgraf
 \setbox\z@\lastbox\unskip\unpenalty\c@nc@t{#1}\nt$\hfil\penalty-\@M
 \else\endgraf\c@nc@t{#1}\nt\fi}
\def\adjustnext#1{\P@nct\hbox{#1}\ignore}
\def\adjustend#1{\def\@djp@{#1}\ignore}
\def\punct#1{\adjustend{#1\space}\nopunct}
\def\cl@s@{\adjust{\@djp@}\endgraf\setbox\z@\lastbox
 \global\setbox\@ne\hbox{\unhbox\adjb@x\ifvoid\z@\else\unhbox\z@\unskip\unskip
 \unpenalty\fi}\egroup\ifnum\c@rr@nt=\k@yb@x\global\fi
 \setbox\c@rr@nt\hbox{\unhbox\@ne\box\p@nct@}\P@nct\null}
\def\@p@n#1{\def\c@rr@nt{#1}\setbox\c@rr@nt\vbox\bgroup\let\@djp@\relax
 \hsize\maxdimen\nt}
\def\b@g@nr@f{\bgroup\@p@n\z@}
\def\key{\cl@s@\ifvoid\k@yb@x\@p@n\k@yb@x\else\@p@n\z@\fi}
\def\label{\cl@s@\ifvoid\k@yb@x\global\@fl@btrue\@p@n\k@yb@x\else\@p@n\z@\fi}
\def\no{\cl@s@\ifvoid\k@yb@x\gad\rcount\global\@fn@mtrue
 \K@yb@x\hbox{\the\rcount}\fi\@p@n\z@}
\def\labelno{\cl@s@\ifvoid\k@yb@x\gad\rcount\@fl@btrue\@p@n\k@yb@x\the\rcount
 \else\@p@n\z@\fi}
\def\by{\cl@s@\@p@n\b@b@x} \def\paper{\cl@s@\@p@n\p@p@rb@x\it\ignore}
\def\jour{\cl@s@\@p@n\j@@rb@x} \def\yr{\cl@s@\@p@n\y@@rb@x}
\def\vol{\cl@s@\@p@n\v@lb@x\bf\ignore} \def\issue{\cl@s@\@p@n\is@b@x}
\def\page{\cl@s@\ifp@g@s\@p@n\z@\else\p@g@true\@p@n\p@g@b@x\fi}
\def\pages{\cl@s@\ifp@g@\@p@n\z@\else\p@g@strue\@p@n\p@g@b@x\fi}
\def\inbook{\cl@s@\@p@n\inb@@kb@x} \def\book{\cl@s@\@p@n\b@@kb@x\it\ignore}
\def\publ{\cl@s@\@p@n\p@blb@x} \def\publaddr{\cl@s@\@p@n\p@bl@db@x}
\def\ed{\cl@s@\ifed@s\@p@n\z@\else\ed@true\@p@n\ed@b@x\fi}
\def\eds{\cl@s@\ifed@\@p@n\z@\else\ed@strue\@p@n\ed@b@x\fi}
\def\info{\cl@s@\@p@n\inf@b@x} \def\paperinfo{\cl@s@\@p@n\p@p@nf@b@x}
\def\bookinfo{\cl@s@\@p@n\b@@nf@b@x} \let\finalinfo\info
\def\P@nct{\global\setbox\p@nct@} \def\nopunct{\P@nct\box\voidb@x}
\def\p@@@t#1#2{\ifvoid\p@nct@\else#1\unhbox\p@nct@#2\fi}
\def\sp@@{\penalty-50 \space\hskip\z@ plus.1em}
\def\c@mm@{\p@@@t,\sp@@} \def\sp@c@{\p@@@t\empty\sp@@} \def\p@@nt{.\kern.3em}
\def\p@tb@x#1#2{\ifvoid#1\else#2\@nb@x#1\fi}
\def\@nb@x#1{\unhbox#1\P@nct\lastbox}
\def\endr@f@{\cl@s@\nopunct
 \R@fb@x\hbox{\unhbox\r@f@b@x \p@tb@x\b@b@x\empty
 \ifvoid\j@@rb@x\ifvoid\inb@@kb@x\ifvoid\p@p@rb@x\ifvoid\b@@kb@x
  \ifvoid\p@p@nf@b@x\ifvoid\b@@nf@b@x
  \p@tb@x\v@lb@x\c@mm@ \ifvoid\y@@rb@x\else\sp@c@(\@nb@x\y@@rb@x)\fi
  \p@tb@x\is@b@x{\c@mm@ no\p@@nt}\p@tb@x\p@g@b@x\c@mm@ \p@tb@x\inf@b@x\c@mm@
  \else\p@tb@x \b@@nf@b@x\c@mm@ \p@tb@x\v@lb@x\c@mm@
  \p@tb@x\is@b@x{\sp@c@ no\p@@nt}%
  \ifvoid\ed@b@x\else\sp@c@(\@nb@x\ed@b@x,\space\ifed@ ed.\else eds.\fi)\fi
  \p@tb@x\p@blb@x\c@mm@ \p@tb@x\p@bl@db@x\c@mm@ \p@tb@x\y@@rb@x\c@mm@
  \p@tb@x\p@g@b@x{\c@mm@\ifp@g@ p\p@@nt\else pp\p@@nt\fi}%
  \p@tb@x\inf@b@x\c@mm@\fi
  \else \p@tb@x\p@p@nf@b@x\c@mm@ \p@tb@x\v@lb@x\c@mm@
  \ifvoid\y@@rb@x\else\sp@c@(\@nb@x\y@@rb@x)\fi
  \p@tb@x\is@b@x{\c@mm@ no\p@@nt}\p@tb@x\p@g@b@x\c@mm@ \p@tb@x\inf@b@x\c@mm@\fi
  \else \p@tb@x\b@@kb@x\c@mm@
  \p@tb@x\b@@nf@b@x\c@mm@ \p@tb@x\p@blb@x\c@mm@
  \p@tb@x\p@bl@db@x\c@mm@ \p@tb@x\y@@rb@x\c@mm@
  \ifvoid\p@g@b@x\else\c@mm@\@nb@x\p@g@b@x p\fi \p@tb@x\inf@b@x\c@mm@ \fi
  \else \c@mm@\@nb@x\p@p@rb@x\ic@\p@tb@x\p@p@nf@b@x\c@mm@
  \p@tb@x\v@lb@x\sp@c@ \ifvoid\y@@rb@x\else\sp@c@(\@nb@x\y@@rb@x)\fi
  \p@tb@x\is@b@x{\c@mm@ no\p@@nt}\p@tb@x\p@g@b@x\c@mm@\p@tb@x\inf@b@x\c@mm@\fi
  \else \p@tb@x\p@p@rb@x\c@mm@\ic@\p@tb@x\p@p@nf@b@x\c@mm@
  \c@mm@\@nb@x\inb@@kb@x \p@tb@x\b@@nf@b@x\c@mm@ \p@tb@x\v@lb@x\sp@c@
  \p@tb@x\is@b@x{\sp@c@ no\p@@nt}%
  \ifvoid\ed@b@x\else\sp@c@(\@nb@x\ed@b@x,\space\ifed@ ed.\else eds.\fi)\fi
  \p@tb@x\p@blb@x\c@mm@ \p@tb@x\p@bl@db@x\c@mm@ \p@tb@x\y@@rb@x\c@mm@
  \p@tb@x\p@g@b@x{\c@mm@\ifp@g@ p\p@@nt\else pp\p@@nt\fi}%
  \p@tb@x\inf@b@x\c@mm@\fi
  \else\p@tb@x\p@p@rb@x\c@mm@\ic@\p@tb@x\p@p@nf@b@x\c@mm@\p@tb@x\j@@rb@x\c@mm@
  \p@tb@x\v@lb@x\sp@c@ \ifvoid\y@@rb@x\else\sp@c@(\@nb@x\y@@rb@x)\fi
  \p@tb@x\is@b@x{\c@mm@ no\p@@nt}\p@tb@x\p@g@b@x\c@mm@ \p@tb@x\inf@b@x\c@mm@
 \fi}}
\def\m@r@f#1#2{\endr@f@\ifvoid\p@nct@\else\R@fb@x\hbox{\unhbox\r@f@b@x
 #1\unhbox\p@nct@\penalty-200\enskip#2}\fi\egroup\b@g@nr@f}
\def\endref{\endr@f@\ifvoid\p@nct@\else\R@fb@x\hbox{\unhbox\r@f@b@x.}\fi
 \parindent\r@f@nd
 \r@fit@m{\ifvoid\k@yb@x\else\if@fn@m\r@fn@{\unhbox\k@yb@x}\else
 \if@fl@b\unhbox\k@yb@x\else\r@fk@y{\unhbox\k@yb@x}\fi\fi\fi}\unhbox\r@f@b@x
 \endgraf\egroup\endgroup}
\def\moreref{\m@r@f;\empty}
\def\transl{\m@r@f;{\unskip\space
 {\sl English translation\ic@}:\penalty-66 \space}}
\def\endRefs{\endgraf\goodbreak\endgroup}

\hyphenation{acad-e-my acad-e-mies af-ter-thought anom-aly anom-alies
an-ti-deriv-a-tive an-tin-o-my an-tin-o-mies apoth-e-o-ses
apoth-e-o-sis ap-pen-dix ar-che-typ-al as-sign-a-ble as-sist-ant-ship
as-ymp-tot-ic asyn-chro-nous at-trib-uted at-trib-ut-able bank-rupt
bank-rupt-cy bi-dif-fer-en-tial blue-print busier busiest
cat-a-stroph-ic cat-a-stroph-i-cally con-gress cross-hatched data-base
de-fin-i-tive de-riv-a-tive dis-trib-ute dri-ver dri-vers eco-nom-ics
econ-o-mist elit-ist equi-vari-ant ex-quis-ite ex-tra-or-di-nary
flow-chart for-mi-da-ble forth-right friv-o-lous ge-o-des-ic
ge-o-det-ic geo-met-ric griev-ance griev-ous griev-ous-ly
hexa-dec-i-mal ho-lo-no-my ho-mo-thetic ideals idio-syn-crasy
in-fin-ite-ly in-fin-i-tes-i-mal ir-rev-o-ca-ble key-stroke
lam-en-ta-ble light-weight mal-a-prop-ism man-u-script mar-gin-al
meta-bol-ic me-tab-o-lism meta-lan-guage me-trop-o-lis
met-ro-pol-i-tan mi-nut-est mol-e-cule mono-chrome mono-pole
mo-nop-oly mono-spline mo-not-o-nous mul-ti-fac-eted mul-ti-plic-able
non-euclid-ean non-iso-mor-phic non-smooth par-a-digm par-a-bol-ic
pa-rab-o-loid pa-ram-e-trize para-mount pen-ta-gon phe-nom-e-non
post-script pre-am-ble pro-ce-dur-al pro-hib-i-tive pro-hib-i-tive-ly
pseu-do-dif-fer-en-tial pseu-do-fi-nite pseu-do-nym qua-drat-ic
quad-ra-ture qua-si-smooth qua-si-sta-tion-ary qua-si-tri-an-gu-lar
quin-tes-sence quin-tes-sen-tial re-arrange-ment rec-tan-gle
ret-ri-bu-tion retro-fit retro-fit-ted right-eous right-eous-ness
ro-bot ro-bot-ics sched-ul-ing se-mes-ter semi-def-i-nite
semi-ho-mo-thet-ic set-up se-vere-ly side-step sov-er-eign spe-cious
spher-oid spher-oid-al star-tling star-tling-ly sta-tis-tics
sto-chas-tic straight-est strange-ness strat-a-gem strong-hold
sum-ma-ble symp-to-matic syn-chro-nous topo-graph-i-cal tra-vers-a-ble
tra-ver-sal tra-ver-sals treach-ery turn-around un-at-tached
un-err-ing-ly white-space wide-spread wing-spread wretch-ed
wretch-ed-ly Brown-ian Eng-lish Euler-ian Feb-ru-ary Gauss-ian
Grothen-dieck Hamil-ton-ian Her-mit-ian Jan-u-ary Japan-ese Kor-te-weg
Le-gendre Lip-schitz Lip-schitz-ian Mar-kov-ian Noe-ther-ian
No-vem-ber Rie-mann-ian Schwarz-schild Sep-tem-ber}

\def\leftheadtext#1{} \def\rightheadtext#1{}

\let\nopagenumber\p@gen@false \let\putpagenumber\p@gen@true
\let\pagefirst\nopagenumber \let\pagenext\putpagenumber

\else

\amsppttrue

\let\twelvepoint\relax \let\Twelvepoint\relax \let\putpagenumber\relax
\let\logo@\relax \let\pagefirst\firstpage@true \let\pagenext\firstpage@false
\def\nopagenumber{\let\f@li@ld\folio\def\folio{\global\let\folio\f@li@ld}}

\def\ftext#1{\footnotetext""{\vsk-.8>\nt #1}}

\def\m@R@f@[#1]{\Refs\nofrills{}\m@th\tenpoint
 {
 \s@ct{}{#1}}\wr@@c{\string\Refcd{#1}{\the\pageno}}
 \def\k@yf@##1{\hss[##1]\enspace} \let\keyformat\k@yf@
 \def\widest##1{\s@twd@\refindentwd{\tenpoint\k@yf@{##1}}}
 \let\Key\key \def\refin{\kern\refindentwd}}
\let\info\finalinfo \r@R@fs\Refs
\def\adjust#1{#1} \let\adjustend\relax
\let\adjustnext\adjust \let\punct\adjust

\fi

\outer\def\myRefs{\myR@fs} \r@st@re\proclaim
\def\bye{\par\vfill\supereject\cl@selbl\cl@secd\b@e} \r@endd@\b@e
\let\Cite\cite \let\Key\key \def\endpro{\par\endproclaim}
\let\d@c@\document \def\document{\d@c@\tenpoint}
\hyphenation{ortho-gon-al}

\newtoks\@@tp@t \@@tp@t\output
\output=\@ft@{\let\{\noexpand\the\@@tp@t}
\let\{\relax

\newif\ifVersion

\def\s@ct#1#2{\ifVersion
 \skip@\lastskip\ifdim\skip@<1.5\bls\vskip-\skip@\p@n@l{-200}\vsk.5>%
 \p@n@l{-200}\vsk.5>\p@n@l{-200}\vsk.5>\p@n@l{-200}\vsk-1.5>\else
 \p@n@l{-200}\fi\ifdim\skip@<.9\bls\vsk.9>\else
 \ifdim\skip@<1.5\bls\vskip\skip@\fi\fi
 \vtop{\twelvepoint\raggedright\bf\vp1\vsk->\vskip.16ex\s@twd@\parindent{#1}%
 \ifdim\parindent>\z@\adv\parindent.5em\fi\hang\textindent{#1}#2\strut}
 \else
 \p@sk@p{-200}{.8\bls}\vtop{\bf\s@twd@\parindent{#1}%
 \ifdim\parindent>\z@\adv\parindent.5em\fi\hang\textindent{#1}#2\strut}\fi
 \nointerlineskip\nobreak\vtop{\strut}\nobreak\vskip-.6\bls\nobreak}

\def\p@n@l#1{\ifnum#1=\z@\else\penalty#1\relax\fi}

\def\s@bs@ct#1#2{\ifVersion
 \skip@\lastskip\ifdim\skip@<1.5\bls\vskip-\skip@\p@n@l{-200}\vsk.5>%
 \p@n@l{-200}\vsk.5>\p@n@l{-200}\vsk.5>\p@n@l{-200}\vsk-1.5>\else
 \p@n@l{-200}\fi\ifdim\skip@<.9\bls\vsk.9>\else
 \ifdim\skip@<1.5\bls\vskip\skip@\fi\fi
 \vtop{\elevenpoint\raggedright\it\vp1\vsk->\vskip.16ex%
 \s@twd@\parindent{#1}\ifdim\parindent>\z@\adv\parindent.5em\fi
 \hang\textindent{#1}#2\strut}
 \else
 \p@sk@p{-200}{.6\bls}\vtop{\it\s@twd@\parindent{#1}%
 \ifdim\parindent>\z@\adv\parindent.5em\fi\hang\textindent{#1}#2\strut}\fi
 \nointerlineskip\nobreak\vtop{\strut}\nobreak\vskip-.8\bls\nobreak}

\def\gadv{\global\adv} \def\gad#1{\gadv#1\@ne} \def\gadneg#1{\gadv#1-\@ne}

\newcount\t@@n \t@@n=10 \newbox\testbox

\newcount\Sno \newcount\Lno \newcount\Fno

\def\pr@cl#1{\r@st@re\pr@c@\pr@c@{#1}\global\let\pr@c@\relax}

\def\tagg#1{\tag"\rlap{\rm(#1)}\kern.01\p@"}
\def\l@L#1{\l@bel{#1}L} \def\l@F#1{\l@bel{#1}F} \def\<#1>{\l@b@l{#1}F}
\def\Tag#1{\tag{\l@F{#1}}} \def\Tagg#1{\tagg{\l@F{#1}}}
\def\Rem{\demo{\sl Remark}} 
\def\Pf#1.{\demo{Proof #1}} \def\epf{\qed\enddemo}
\def\Ap@x{Appendix}
\def\Appendix{\Sno=64 \t@@n\@ne \wr@@c{\string\Appencd}
 \def\sf@rm{\char\the\Sno} \def\sf@rm@{\Ap@x\space\sf@rm} \def\sf@rm@@{\Ap@x}
 \def\s@ct@n##1##2{\s@ct\empty{\setbox\z@\hbox{##1}\ifdim\wd\z@=\z@
 \if##2*\sf@rm@@\else\if##2.\sf@rm@@.\else##2\fi\fi\else
 \if##2*\sf@rm@\else\if##2.\sf@rm@.\else\sf@rm@.\enspace##2\fi\fi\fi}}}
\def\Appcd#1#2#3{\def\Ap@@{\hglue-\l@ftcd\Ap@x}\ifx\@ppl@ne\empty
 \def\l@@b{\@fwd@@{#1}{\space#1}{}}\if*#2\entcd{}{\Ap@@\l@@b}{#3}\else
 \if.#2\entcd{}{\Ap@@\l@@b.}{#3}\else\entcd{}{\Ap@@\l@@b.\enspace#2}{#3}\fi\fi
 \else\def\l@@b{\@fwd@@{#1}{\c@l@b{#1}}{}}\if*#2\entcd{\l@@b}{\Ap@x}{#3}\else
 \if.#2\entcd{\l@@b}{\Ap@x.}{#3}\else\entcd{\l@@b}{#2}{#3}\fi\fi\fi}

\let\s@ct@n\s@ct
\def\s@ct@@[#1]#2{\@ft@\xdef\csname @#1@S@\endcsname{\sf@rm}\wr@@x{}%
 \wr@@x{\string\labeldef{S}\space{\?#1@S?}\space{#1}}%
 {
 \s@ct@n{\sf@rm@}{#2}}\wr@@c{\string\Entcd{\?#1@S?}{#2}{\the\pageno}}}
\def\s@ct@#1{\wr@@x{}{
 \s@ct@n{\sf@rm@}{#1}}\wr@@c{\string\Entcd{\sf@rm}{#1}{\the\pageno}}}
\def\s@ct@e[#1]#2{\@ft@\xdef\csname @#1@S@\endcsname{\sf@rm}\wr@@x{}%
 \wr@@x{\string\labeldef{S}\space{\?#1@S?}\space{#1}}%
 {
 \s@ct@n\empty{#2}}\wr@@c{\string\Entcd{}{#2}{\the\pageno}}}
\def\s@cte#1{\wr@@x{}{
 \s@ct@n\empty{#1}}\wr@@c{\string\Entcd{}{#1}{\the\pageno}}}
\def\theSno#1#2{\dff\?#1@S?{#2}%
 \wr@@x{\string\labeldef{S}\space{#2}\space{#1}}\fi}

\newif\ifd@bn@\d@bn@true
\def\Section{\gad\Sno\ifd@bn@\Fno\z@\Lno\z@\fi\@fn@xt[\s@ct@@\s@ct@}
\def\section{\gad\Sno\ifd@bn@\Fno\z@\Lno\z@\fi\@fn@xt[\s@ct@e\s@cte}
\let\Sect\Section 
\def\subsection{\@fn@xt*\subs@ct@\subs@ct}
\def\subs@ct#1{{
 \s@bs@ct\empty{#1}}\wr@@c{\string\subcd{#1}{\the\pageno}}}
\def\subs@ct@*#1{\vsk->\vsk>{
 \s@bs@ct\empty{#1}}\wr@@c{\string\subcd{#1}{\the\pageno}}}
\let\subsect\subsection \def\Snodef#1{\Sno #1}

\def\l@b@l#1#2{\def\n@@{\csname #2no\endcsname}%
 \if*#1\gad\n@@ \@ft@\xdef\csname @#1@#2@\endcsname{\l@f@rm}\else\def\t@st{#1}%
 \ifx\t@st\empty\gad\n@@ \@ft@\xdef\csname @#1@#2@\endcsname{\l@f@rm}%
 \else\@ft@\ifx\csname @#1@#2@mark\endcsname\relax\gad\n@@
 \@ft@\xdef\csname @#1@#2@\endcsname{\l@f@rm}%
 \@ft@\gdef\csname @#1@#2@mark\endcsname{}%
 \wr@@x{\string\labeldef{#2}\space{\?#1@#2?}\space\ifnum\n@@<10 \space\fi{#1}}%
 \fi\fi\fi}
\def\labeldef#1#2#3{\dff\?#3@#1?{#2}}
\def\Labeldef#1#2#3{\dff\?#3@#1?{#2}\@ft@\gdef\csname @#3@#1@mark\endcsname{}}

\def\l@bel#1#2{\l@b@l{#1}{#2}\?#1@#2?}

\newcount\c@cite
\def\?#1?{\csname @#1@\endcsname}
\def\[{\@fn@xt:\c@t@sect\c@t@}
\def\c@t@#1]{{\c@cite\z@\@fwd@@{\?#1@L?}{\adv\c@cite1}{}%
 \@fwd@@{\?#1@F?}{\adv\c@cite1}{}\@fwd@@{\?#1?}{\adv\c@cite1}{}%
 \relax\ifnum\c@cite=\z@{\bf ???}\wrs@x{No label [#1]}\else
 \ifnum\c@cite=1\let\@@PS\relax\let\@@@\relax\else\let\@@PS\underbar
 \def\@@@{{\rm<}}\fi\@@PS{\?#1?\@@@\?#1@L?\@@@\?#1@F?}\fi}}
\def\(#1){{\rm(\c@t@#1])}}
\def\c@t@s@ct#1{\@fwd@@{\?#1@S?}{\?#1@S?\relax}%
 {{\bf ???}\wrs@x{No section label {#1}}}}
\def\c@t@sect:#1]{\c@t@s@ct{#1}} \let\SNo\c@t@s@ct

\newdimen\l@ftcd \newdimen\r@ghtcd \let\nlc\relax

\def\d@tt@d{\leaders\hbox to 1em{\kern.1em.\hfil}\hfill}
\def\entcd#1#2#3{\item{#1}{#2}\alb\kern.9em\hbox{}\kern-.9em\d@tt@d
 \kern-.36em{#3}\kern-\r@ghtcd\hbox{}\par}
\def\Entcd#1#2#3{\def\l@@b{\@fwd@@{#1}{\c@l@b{#1}}{}}\vsk.2>%
 \entcd{\l@@b}{#2}{#3}}
\def\subcd#1#2{{\adv\leftskip.333em\entcd{}{\it #1}{#2}}}
\def\Refcd#1#2{\def\t@@st{#1}\ifx\t@@st\empty\ifx\r@fl@ne\empty\relax\else
 \R@fcd{\r@fl@ne}{#2}\fi\else\R@fcd{#1}{#2}\fi}
\def\R@fcd#1#2{\sk@@p{.6\bls}\entcd{}{\hglue-\l@ftcd\bf #1}{#2}}
\def\Refline{\def\r@fl@ne} \def\Refempty{\let\r@fl@ne\empty}
\def\Appencd{\par\adv\leftskip-\l@ftcd\adv\rightskip-\r@ghtcd\@ppl@ne
 \adv\leftskip\l@ftcd\adv\rightskip\r@ghtcd\let\Entcd\Appcd}
\def\appline{\def\@ppl@ne} \def\Appempty{\let\@ppl@ne\empty}
\def\Appline#1{\def\@ppl@ne{\s@bs@ct{}{\bf#1}}}
\def\leftcd#1{\adv\leftskip-\l@ftcd\s@twd@\l@ftcd{\c@l@b{#1}\enspace}
 \adv\leftskip\l@ftcd}
\def\rightcd#1{\adv\rightskip-\r@ghtcd\s@twd@\r@ghtcd{#1\enspace}
 \adv\rightskip\r@ghtcd}
\def\C@nt{Contents} \def\Ap@s{Appendices} \def\R@fcs{References}
\def\contents{\@fn@xt*\cont@@\cont@}
\def\cont@{\@fn@xt[\cnt@{\cnt@[\C@nt]}}
\def\cont@@*{\@fn@xt[\cnt@@{\cnt@@[\C@nt]}}
\def\cnt@[#1]{\c@nt@{M}{#1}{44}{\s@bs@ct{}{\bf\Ap@s}}}
\def\cnt@@[#1]{\c@nt@{M}{#1}{44}{}}
\def\endco{\par\penalty-500\vsk>\vskip\z@\endgroup}
\def\readcd{\@np@t{\jobname.cd}}
\def\Cde{\@fn@xt*\Cde@@\Cde@}
\def\Cde@{\@fn@xt[\Cd@{\Cd@[\C@nt]}}
\def\Cde@@*{\@fn@xt[\Cd@@{\Cd@@[\C@nt]}}
\def\Cd@[#1]{\cnt@[#1]\readcd\endco}
\def\Cd@@[#1]{\cnt@@[#1]\readcd\endco}
\def\contlabeldef{\def\c@l@b}

\long\def\c@nt@#1#2#3#4{\s@twd@\l@ftcd{\c@l@b{#1}\enspace}
 \s@twd@\r@ghtcd{#3\enspace}\adv\r@ghtcd1.333em
 \def\@ppl@ne{#4}\def\r@fl@ne{\R@fcs}\s@ct{}{#2}\B@gr@\parindent\z@\let\nlc\nl
 \let\nl\relax\parskip.2\bls\adv\leftskip\l@ftcd\adv\rightskip\r@ghtcd}

\def\writecd{\immediate\openout\@@cd\jobname.cd \def\wr@@c{\write\@@cd}
 \def\cl@secd{\immediate\write\@@cd{\string\endinput}\immediate\closeout\@@cd}
 \def\closecd{\cl@secd\global\let\cl@secd\relax}}
\let\cl@secd\relax \def\wr@@c#1{} \let\closecd\relax

\def\dff{\@ft@\d@f} \def\d@f{\@ft@\def}
\def\edff{\@ft@\ed@f} \def\ed@f{\@ft@\edef}
\def\defi#1#2{\def#1{#2}\wr@@x{\string\def\string#1{#2}}}

\def\qed{\hbox{}\nobreak\hfill\nobreak{\m@th$\,\square$}}
\def\back#1 {\strut\kern-.33em #1\enspace\ignore} 
\def\Text#1{\crcr\noalign{\alb\vsk>\normalbaselines\vsk->\vbox{\nt #1\strut}%
 \nobreak\nointerlineskip\vbox{\strut}\nobreak\vsk->\nobreak}}

\def\hcor#1{\advance\hoffset by #1}
\def\vcor#1{\advance\voffset by #1}
\let\bls\baselineskip \let\ignore\ignorespaces
\ifx\ic@\undefined \let\ic@\/\fi
\def\vsk#1>{\vskip#1\bls} \let\adv\advance
\def\vv#1>{\vadjust{\vsk#1>}\ignore}
\def\vvn#1>{\vadjust{\nobreak\vsk#1>\nobreak}\ignore}
\def\vvv#1>{\vskip\z@\vsk#1>\nt\ignore}
\def\vvgood{\vadjust{\penalty-500}} \def\vgood{\strut\vvn->\vvgood\nl\ignore}
\def\Par{\vsk.5>} \def\setparindent{\edef\Parindent{\the\parindent}}
\def\Type{\vsk.5>\bgroup\parindent\z@\tt\rightskip\z@ plus1em minus1em%
 \spaceskip.3333em \xspaceskip.5em\relax}
\def\endType{\vsk.5>\egroup\nt} 

 \let\Tilde\widetilde \let\dollar\$ \let\ampersand\&
\let\sss\scriptscriptstyle  
\let\vp\vphantom \let\hp\hphantom \let\nt\noindent
\let\cline\centerline \let\lline\leftline \let\rline\rightline
\def\nn#1>{\noalign{\vskip#1\p@@}} \def\NN#1>{\openup#1\p@@}
\def\cnn#1>{\noalign{\vsk#1>}}
 
\let\Lim\lim \def\lim{\Lim\limits} \let\Sum\sum \def\sum{\Sum\limits}
\def\Plus{\bigoplus\limits} 
\let\Prod\prod \def\prod{\Prod\limits} \let\Int\int \def\int{\Int\limits}

\def\tprod{\mathop{\tsize\Prod}\limits}
\def\&{.\kern.1em} \def\>{{\!\;}} \def\]{{\!\!\;}} \def\){\>\]} \def\}{\]\]}
\def\nl{\leavevmode\hfill\break} \def\~{\leavevmode\@fn@xt~\m@n@s\@md@sh}
\def\m@n@s~{\raise.15ex\mbox{-}} \def\@md@sh{\raise.13ex\hbox{--}}
\let\procent\% \def\%#1{\ifmmode\mathop{#1}\limits\else\procent#1\fi}
\let\@ml@t\" \def\"#1{\ifmmode ^{(#1)}\else\@ml@t#1\fi}
\let\@c@t@\' \def\'#1{\ifmmode _{(#1)}\else\@c@t@#1\fi}
\let\colon\: \def\:{^{\vp|}}

\let\texspace\ \def\ {\ifmmode\alb\fi\texspace}

\let\n@wp@ge\newpage \def\newpage{\endgraf\n@wp@ge}
\let\=\m@th \def\mbox#1{\hbox{\m@th$#1$}}
\def\mtext#1{\text{\m@th$#1$}} \def\^#1{\text{\m@th#1}}
\def\Line#1{\kern-.5\hsize\line{\m@th$\dsize#1$}\kern-.5\hsize}
\def\Lline#1{\kern-.5\hsize\lline{\m@th$\dsize#1$}\kern-.5\hsize}
\def\Cline#1{\kern-.5\hsize\cline{\m@th$\dsize#1$}\kern-.5\hsize}
\def\Rline#1{\kern-.5\hsize\rline{\m@th$\dsize#1$}\kern-.5\hsize}

\def\Ll@p#1{\llap{\m@th$#1$}} \def\Rl@p#1{\rlap{\m@th$#1$}}
 \def\Cl@p#1{\llap{\m@th$#1$\hss}}
\def\Llap#1{\mathchoice{\Ll@p{\dsize#1}}{\Ll@p{\tsize#1}}{\Ll@p{\ssize#1}}%
 {\Ll@p{\sss#1}}}
\def\Clap#1{\mathchoice{\Cl@p{\dsize#1}}{\Cl@p{\tsize#1}}{\Cl@p{\ssize#1}}%
 {\Cl@p{\sss#1}}}
\def\Rlap#1{\mathchoice{\Rl@p{\dsize#1}}{\Rl@p{\tsize#1}}{\Rl@p{\ssize#1}}%
 {\Rl@p{\sss#1}}}
 
\def\LRtph#1#2{\setbox\z@\hbox{#1}\dimen\z@\wd\z@\hbox{\hbox to\dimen\z@{#2}}}
\def\LRph#1#2{\LRtph{\m@th$#1$}{\m@th$#2$}}

\def\CCph#1#2{\LRph{#1}{\hss#2\hss}}
 \def\RRph#1#2{\LRph{#1}{#2\hss}}
\def\Cph#1#2{\mathchoice{\CCph{\dsize#1}{\dsize#2}}{\CCph{\tsize#1}{\tsize#2}}
 {\CCph{\ssize#1}{\ssize#2}}{\CCph{\sss#1}{\sss#2}}}

\def\Rph#1#2{\mathchoice{\RRph{\dsize#1}{\dsize#2}}{\RRph{\tsize#1}{\tsize#2}}
 {\RRph{\ssize#1}{\ssize#2}}{\RRph{\sss#1}{\sss#2}}}
\def\Lto#1{\setbox\z@\mbox{\tsize{#1}}%
 \mathrel{\mathop{\hbox to\wd\z@{\rightarrowfill}}\limits#1}}
\def\Lgets#1{\setbox\z@\mbox{\tsize{#1}}%
 \mathrel{\mathop{\hbox to\wd\z@{\leftarrowfill}}\limits#1}}
\def\vpb#1{{\vp{\big(}}^{\]#1}} \def\vpp#1{{\vp{\big]}}_{#1}}
\def\lbc{\mathopen{[\![}} \def\rbc{\mathclose{]\!]}}
 
\def\LBc{\mathopen{\Big[\!\}\Big[}} \def\RBc{\mathclose{\Big]\!\}\Big]}}
\def\lpt{\mathopen{(\!(}} \def\rpt{\mathclose{)\!)}}

\let\alb\allowbreak 
\def\ald{\noalign{\alb}} \let\alds\allowdisplaybreaks

 \let\x\times \let\ox\otimes 
  \let\tabs\+
\let\le\leqslant \let\ge\geqslant
 \let\8\infty \let\*\star
\let\bra\langle \let\ket\rangle
 
\let\map\mapsto  
 \let\rto\rightrightarrows
 \def\vert{\ |\ } \def\nin{\not\in}

\let\lb\lbrace \let\rb\rbrace
 
 \let\tright\triangleright

\def\lsym#1{#1\alb\ldots\relax#1\alb}
\def\lc{\lsym,}   
\def\llc{\,,\alb\ {\ldots\ ,}\alb\ }
 
\def\End{\mathop{\roman{End}\>}}

\def\Res{\mathop{\roman{Res}\>}\limits}

  \def\for{\text{for }\,}
\def\1{^{-1}} \def\_#1{_{\Rlap{#1}}}
\def\vst#1{{\lower1.9\p@@\mbox{\bigr|_{\raise.5\p@@\mbox{\ssize#1}}}}}
\def\vrp#1:#2>{{\vrule height#1 depth#2 width\z@}}
\def\vru#1>{\vrp#1:\z@>} \def\vrd#1>{\vrp\z@:#1>}
\def\qqq{\qquad\quad} 
\def\sscr#1{\raise.3ex\mbox{\sss#1}} \def\@@PS{\bold{OOPS!!!}}

\def\intcl{\mathop
 {\Rlap{\raise.3ex\mbox{\kern.12em\curvearrowleft}}\int}\limits}
\def\intcr{\mathop
 {\Rlap{\raise.3ex\mbox{\kern.24em\curvearrowright}}\int}\limits}

\def\pms{\raise.25ex\mbox{\ssize\pm}\>}
\def\mps{\raise.25ex\mbox{\ssize\mp}\>}

\let\al\alpha
\let\bt\beta
\let\gm\gamma \let\Gm\Gamma 
\let\dl\delta  
 \let\eps\varepsilon \let\epsilon\eps

\let\zt\zeta
\let\tht\theta \let\Tht\Theta
\let\thi\vartheta 

\let\ka\kappa

\let\si\sigma 
 
 \let\phi\varphi

\let\om\omega  

\def\C{\Bbb C}

\def\Z{\Bbb Z}

\def\AA{\Bbb A}

\def\TT{\Bbb T}

\def\Zp{\Z_{\ge 0}} 
\def\Zn{\Z_{\le 0}} 

\def\difl/{differential} \def\dif/{difference}
\def\cf.{cf.\ \ignore} \def\Cf.{Cf.\ \ignore}
\def\egv/{eigenvector} \def\eva/{eigenvalue} \def\eq/{equation}
\def\lhs/{the left hand side} \def\rhs/{the right hand side}
\def\Lhs/{The left hand side} \def\Rhs/{The right hand side}
\def\gby/{generated by} \def\wrt/{with respect to} \def\st/{such that}
\def\resp/{respectively} \def\off/{offdiagonal} \def\wt/{weight}
\def\pol/{polynomial} \def\rat/{rational} \def\tri/{trigonometric}
\def\fn/{function} \def\var/{variable} \def\raf/{\rat/ \fn/}
\def\inv/{invariant} \def\hol/{holomorphic} \def\hof/{\hol/ \fn/}
\def\mer/{meromorphic} \def\mef/{\mer/ \fn/} \def\mult/{multiplicity}
\def\sym/{symmetric} \def\perm/{permutation} \def\fd/{finite-dimensional}
\def\rep/{representation} \def\irr/{irreducible} \def\irrep/{\irr/ \rep/}
\def\hom/{homomorphism} \def\aut/{automorphism} \def\iso/{isomorphism}
\def\lex/{lexicographical} \def\as/{asymptotic} \def\asex/{\as/ expansion}
\def\ndeg/{nondegenerate} \def\neib/{neighbourhood} \def\deq/{\dif/ \eq/}
\def\hw/{highest \wt/} \def\gv/{generating vector} \def\eqv/{equivalent}
\def\msd/{method of steepest descend} \def\pd/{pairwise distinct}
\def\wlg/{without loss of generality} \def\Wlg/{Without loss of generality}
\def\onedim/{one-dimensional} \def\qcl/{quasiclassical}
\def\hgeom/{hyper\-geometric} \def\hint/{\hgeom/ integral}
\def\hwm/{\hw/ module} \def\emod/{evaluation module} \def\Vmod/{Verma module}
\def\symg/{\sym/ group} \def\sol/{solution} \def\eval/{evaluation}
\def\anf/{analytic \fn/} \def\anco/{analytic continuation}
\def\qg/{quantum group} \def\qaff/{quantum affine algebra}

\def\Rm/{\^{$R$-}matrix} \def\Rms/{\^{$R$-}matrices} \def\YB/{Yang-Baxter \eq/}
\def\Ba/{Bethe ansatz} \def\Bv/{Bethe vector} \def\Bae/{\Ba/ \eq/}
\def\KZv/{Knizh\-nik-Zamo\-lod\-chi\-kov} \def\KZvB/{\KZv/-Bernard}
\def\KZ/{{\sl KZ\/}} \def\qKZ/{{\sl qKZ\/}}
\def\KZB/{{\sl KZB\/}} \def\qKZB/{{\sl qKZB\/}}
\def\qKZo/{\qKZ/ operator} \def\qKZc/{\qKZ/ connection}
\def\KZe/{\KZ/ \eq/} \def\qKZe/{\qKZ/ \eq/} \def\qKZBe/{\qKZB/ \eq/}

\def\h@ph{\discretionary{}{}{-}} \def\$#1$-{\,\^{$#1$}\h@ph}

\def\TFT/{Research Insitute for Theoretical Physics}
\def\HY/{University of Helsinki} \def\AoF/{the Academy of Finland}
\def\CNRS/{Supported in part by MAE\~MICECO\~CNRS Fellowship}
\def\LPT/{Laboratoire de Physique Th\'eorique ENSLAPP}
\def\ENSLyon/{\'Ecole Normale Sup\'erieure de Lyon}
\def\LPTaddr/{46, All\'ee d'Italie, 69364 Lyon Cedex 07, France}
\def\enslapp/{URA 14\~36 du CNRS, associ\'ee \`a l'E.N.S.\ de Lyon,
au LAPP d'Annecy et \`a l'Universit\`e de Savoie}
\def\ensemail/{vtarasov\@ enslapp.ens-lyon.fr}
\def\DMS/{Department of Mathematics, Faculty of Science}
\def\DMO/{\DMS/, Osaka University}
\def\DMOaddr/{Toyonaka, Osaka 560, Japan}
\def\dmoemail/{vt\@ math.sci.osaka-u.ac.jp}
\def\SPb/{St\&Peters\-burg}
\def\home/{\SPb/ Branch of Steklov Mathematical Institute}
\def\homeaddr/{Fontanka 27, \SPb/ \,191011, Russia}
\def\homemail/{vt\@ pdmi.ras.ru}
\def\absence/{On leave of absence from \home/}
\def\UNC/{Department of Mathematics, University of North Carolina}
\def\ChH/{Chapel Hill}
\def\UNCaddr/{\ChH/, NC 27599, USA} \def\avemail/{av\@ math.unc.edu}
\def\grant/{NSF grant DMS\~9501290}	
\def\Grant/{Supported in part by \grant/}

\def\Aomoto/{K\&Aomoto}
\def\Dri/{V\]\&G\&Drin\-feld}
\def\Fadd/{L\&D\&Fad\-deev}
\def\Feld/{G\&Felder}
\def\Fre/{I\&B\&Fren\-kel}
\def\Gustaf/{R\&A\&Gustafson}
\def\Kazh/{D\&Kazhdan} \def\Kir/{A\&N\&Kiril\-lov}
\def\Kor/{V\]\&E\&Kore\-pin}
\def\Lusz/{G\&Lusztig}
\def\MN/{M\&Naza\-rov}
\def\Resh/{N\&Reshe\-ti\-khin} \def\Reshy/{N\&\]Yu\&Reshe\-ti\-khin}
\def\SchV/{V\]\&\]V\]\&Schecht\-man} \def\Sch/{V\]\&Schecht\-man}
\def\Skl/{E\&K\&Sklya\-nin}
\def\Smirn/{F\]\&Smirnov} \def\Smirnov/{F\]\&A\&Smirnov}
\def\Takh/{L\&A\&Takh\-tajan}
\def\VT/{V\]\&Ta\-ra\-sov} \def\VoT/{V\]\&O\&Ta\-ra\-sov}
\def\Varch/{A\&\]Var\-chenko} \def\Varn/{A\&N\&\]Var\-chenko}

\def\AMS/{Amer.\ Math.\ Society}
\def\CMP/{Comm.\ Math.\ Phys.{}}
\def\DMJ/{Duke.\ Math.\ J.{}}
\def\Inv/{Invent.\ Math.{}} 
\def\IMRN/{Int.\ Math.\ Res.\ Notices}
\def\JPA/{J.\ Phys.\ A{}}
\def\JSM/{J.\ Soviet\ Math.{}}
\def\LMP/{Lett.\ Math.\ Phys.{}}
\def\LMJ/{Leningrad Math.\ J.{}}
\def\LpMJ/{\SPb/ Math.\ J.{}}
\def\SIAM/{SIAM J.\ Math.\ Anal.{}}
\def\SMNS/{Selecta Math., New Series}
\def\TMP/{Theor.\ Math.\ Phys.{}}
\def\ZNS/{Zap.\ nauch.\ semin. LOMI}

\def\ASMP/{Advanced Series in Math.\ Phys.{}}

\def\AMSa/{AMS \publaddr Providence}
\def\Birk/{Birkh\"auser}
\def\CUP/{Cambridge University Press} \def\CUPa/{\CUP/ \publaddr Cambridge}
\def\Spri/{Springer-Verlag} \def\Spria/{\Spri/ \publaddr Berlin}
\def\WS/{World Scientific} \def\WSa/{\WS/ \publaddr Singapore}

\newbox\lefthbox \newbox\righthbox

\let\sectsep. \let\labelsep. \let\contsep. \let\labelspace\relax
\let\sectpre\relax \let\contpre\relax
\def\sf@rm{\the\Sno} \def\sf@rm@{\sectpre\sf@rm\sectsep}
\def\c@l@b#1{\contpre#1\contsep}
\def\l@f@rm{\ifd@bn@\sf@rm\labelsep\fi\labelspace\the\n@@}

\def\sectformdef{\def\sf@rm}

\let\DoubleNum\d@bn@true \let\SingleNum\d@bn@false

\def\NoNewNum{\let\writeldf\relax\def\l@b@l##1##2{\if*##1%
 \@ft@\xdef\csname @##1@##2@\endcsname{\mbox{*{*}*}}\fi}}
\def\NoNewTime{\def\todaydef##1{\def\today{##1}}
 \def\nowtimedef##1{\def\nowtime{##1}}}
\def\NoInput{\let\Input\input\let\writeldf\relax}
\def\Fixed{\NoNewTime\NoInput}

\tenpoint

\csname beta.def\endcsname

\Fixed

\expandafter\ifx\csname hgeom.def\endcsname\relax \else\endinput\fi
\expandafter\edef\csname hgeom.def\endcsname{%
 \catcode`\noexpand\@=\the\catcode`\@\space}
\catcode`\@=11

\font\Brm=cmr12 scaled 1200
\font\Bit=cmti12 scaled 1200

\ifx\twelvemsa\undefined\font\twelvemsa=msam10 scaled 1200\fi
\ifMag\else\fi

\def\Th#1{\pr@cl{Theorem \l@L{#1}}\ignore}
\def\Lm#1{\pr@cl{Lemma \l@L{#1}}\ignore}
\def\Cr#1{\pr@cl{Corollary \l@L{#1}}\ignore}
\def\Df#1{\pr@cl{Definition \l@L{#1}}\ignore}
\def\Cj#1{\pr@cl{Conjecture \l@L{#1}}\ignore}
\def\Prop#1{\pr@cl{Proposition \l@L{#1}}\ignore}

\newif\ifosaka \newcount\rc@@nt
\def\Pf#1.{\demo{\ifosaka\rm\fi Proof #1}}
\def\Rem{\demo{\ifosaka\smc\else\sl\fi Remark}}
\def\osakaRefs{\wr@@x{}\let\by@@\by \def\by{\by@@\punct:}
 \def\Key##1 {\gad\rc@@nt\wr@@x{\string\labeldef{R}\space{\the\rc@@nt}\space
 \ifnum\rc@@nt<10 \space\fi{##1}}\no##1}}
\def\Ref#1{\ifosaka\[#1@R]\else#1\fi}
\def\Cite#1{\alb\@fwd@@{\Ref{#1}}{\cite{\Ref{#1}}}{\cite{\bf???}}}

 \def\pii{\pi i}
 \def\nepsr{\ne{\)p^s\eta^r}} 

\def\Sb{\bold S}
\def\sb{\bold s}

\def\Ec{\Cal E}
\def\Fc{\Cal F}

\def\Zc{\Cal Z}

\def\eg{\frak e}
\def\hg{\frak h}
\def\lg{\frak l}
\def\mg{\frak m}
\def\ng{\frak n}

\def\Ib{\bar I}

\def\Zlnb{\Rlap{\;\overline{\!\}\phantom\Zc\}}\>}\Zln}

\def\alt{\Tilde\al}

\def\Pht{\Tilde\Phi}
\def\TTt{\Tilde\TT}

\def\zti{\tilde z}

\def\Cn{\C^{\>n}} \def\Cl{\C^{\,\ell}} 
 
\def\Cx{\C^\x}  \def\Cxl{\C^{\x\ell}}

\def\Bq{B_{\]\sss e\]l\}l}}
\def\Fq{\Fc_{\!\!\]\sss e\]l\}l}}
\def\Sq{S_{\]\sss e\]l\}l}}

\def\IH{\mathop{\italic{Int}\)}\limits}
\def\Resi{\mathop{\italic{Res}\)}\limits}
\def\limA{\lim\nolimits_{\)\AA}\:}
\def\Val{\mathop{\italic{Value}}}

\def\Sl{\Sb_\ell}  \def\Zpn{\Zp^n}
\def\Zln{\Zc_\ell^n}  

\def\pZ{p^{\)\Z}} \def\ve[#1]{v^{\)[#1]}}
\def\phif#1#2{{\dsize{}_{#1}\:\phi_{#2}\:}} \def\phifn{\phif{n}{n-1}}
\def\Fifn{{\dsize{}_nF_{n-1}}}

\def\tho{\tht'(1)}

\def\l@inf{\lower.21ex\mbox{\ssize\8}} \def\9{_{\kern-.02em\l@inf}}

\def\smm{\sum_{m=1}} 
\def\sun{\smm^n} 
 
\def\susi{\sum_{\>\si\in\Sl\}}}
\def\smiZ{\sum_{\,\mg\in\Zln\!}^{\vp.}}
\def\smiZb{\sum_{\,\mg\in\Zlnb\!}^{\vp.}}
\def\Pm@Z{\Plus_{\,\mg\in\Zln\!}}
\def\PmiZ{\mathchoice{{\Pm@Z}}{\!\Pm@Z\}}{\Pm@Z}{\@@PS}}

\def\prodp{\prod\nolimits^{\sss\]\bullet\]}}
\def\prm{\prod_{m=1}} \def\tprm{\tprod_{m=1}}
\def\pron{\prm^n} \def\tpron{\tprm^n}
\def\prmn{\prm^{n-1}} 
\def\prmk{\prod_{\tsize{m=1\atop m\ne k}}^n}
\def\plmn{\prod_{1\le l<m\le n}}
 
\def\prlm{\prod_{1\le l<m}} \def\prllm{\prod_{1\le l\le m}}
 
\def\prjm{\prod_{1\le j<m}} \def\tprjjm{\!\!\tprod_{1\le j\le m}\!}

\def\pros{\prod_{s=0}^{\ell-1}} 
\def\prsl{\prod_{s=1-\ell}^{\ell-1}}
\def\pral{\prod_{a=1}^\ell} \def\tpral{\tprod_{a=1}^\ell}

\def\prab{\prod_{1\le a<b\le\ell}\!\!}

\def\prabsi{\prod_{\tsize{1\le a<b\le\ell\atop\si_a>\si_b}}\!\!}

\def\an{a_1\lc a_n} \def\bn{b_1\lc b_{n-1}}
\def\ani{a_1\1\}\lc a_n\1} \def\bni{b_1\1\}\lc b_{n-1}\1}
\def\anpi{p\)a_1\1\}\lc p\)a_n\1} \def\anp{p\)a_1\lc p\)a_n}
\def\anb{b_1a_1\1\}\lc b_1a_n\1} \def\anbi{p\)b_1a_1\1\}\lc p\)b_1a_n\1}
\def\anpb{p\)b_1\1a_1\lc p\)b_1\1a_n} \def\anpbm{p\)b_m\1a_1\lc p\)b_m\1a_n}
\def\anbmi{p\)b_ma_1\1\}\lc p\)b_ma_n\1}

\def\kn{k=1\lc n}  \def\lcn{l=1\lc n} \def\mn{m=1\lc n}
\def\xn{x_1\lc x_n} \def\yn{y_1\lc y_n} \def\zn{z_1\lc z_n}
\def\xin{\xi_1\lc\xi_n} \def\elln{\ell_1\lc\ell_n}
\def\xy@{\!\pron\!x_m/y_m} \def\yx@{\!\pron\!y_m/x_m}
\def\xy{\mathchoice{{\tsize\xy@}}{\xy@}{\xy@}{\xy@}}
\def\yx{\mathchoice{{\tsize\yx@}}{\yx@}{\yx@}{\yx@}}
\def\alti{\eta^{2-2\ell}\xy} 
\def\aell{a=1\lc\ell}
  \def\roll{r=0\lc\ell-1}
\def\rolt{r=0\lc 2\ell-2}
\def\tell{t_1\lc t_\ell} \def\tpell{t_1\lc p\)t_a\lc t_\ell}
 
\def\tsill{t_{\si_1}\lc t_{\si_\ell}}
\def\xt{x\mathbin{\smash\tright\]}}
\def\yt{y\mathbin{\smash\tright\]}}

 \def\dtt{(dt/t)^\ell} 
\def\ts{t^{\sss\star}} 

\def\gl{\frak{gl}_2} \def\gsl{\frak{sl}_2} \def\Eqg{E_{\rho,\gm}(\gsl)}
\def\Uu{U_q(\gsl)} 
 \def\Ugg{U_q'(\Tilde{\gl})}

\def\azo/{\as/ zone} \def\asol/{\as/ \sol/}
\def\conn/{connection coefficient} \def\sconn/{system of \conn/s}
\def\loc/{local system} \def\dloc/{discrete local system}
\def\prim/{primitive factor} \def\hform/{\hgeom/ form}
\def\hpair/{\hgeom/ pairing} \def\hgf/{\hgeom/ space}
\def\hmap/{\hgeom/ map} \def\hsol/{\hgeom/ \sol/}
\def\bid/{bilinear identity} \def\hRid/{\hgeom/ Riemann identity}
\def\thgf/{\tri/ \hgf/} \def\ehgf/{elliptic \hgf/}
\def\hcg/{\hgeom/ cohomology group} \def\hhg/{\hgeom/ homology group}
\def\wtd/{\wt/ decomposition} \def\diag/{diagonalizable}
\def\qbeta/{\^{$q$-}beta} \def\qSelberg/{\^{$q$-}Selberg}
\def\dyne/{dynamical elliptic} \def\ecurve/{elliptic curve}
\def\Spair/{Shapo\-va\-lov pairing} \def\Sform/{Shapo\-va\-lov form}
\def\eqg/{elliptic \qg/}

\def\GM/{Gauss-Manin connection} \def\MB/{Mellin-Barnes}

\def\stype/{\^{$\gsl$-}type}
\def\Umod/{\^{$\Uu$-}module} \def\Uhmod/{\^{$\Ugg$-}module}
\def\Emod/{\^{$\Eqg$-}module} \def\hmod/{\^{$\hg$-}module}
\def\dhmod/{\diag/ \^{$\hg$-}module}
\def\eVmod/{evaluation \Vmod/} \def\epoint/{evaluation point}

\def\tenco/{tensor coordinates} \def\qlo/{quantum loop algebra}
\def\phf/{phase \fn/} \def\wtf/{\wt/ \fn/}
\def\twf/{\tri/ \wtf/} \def\ewf/{elliptic \wtf/} \def\Atype/{{\sl A\/}-type}
\def\cosub/{coboundary subspace}\def\bosub/{boundary subspace}

\def\Vval/{\^{$V\!\!\;$-}valued} 

\def\qhgeom/{\^{$q$-}\hgeom/} \def\qhint/{\^{$q$-}\hint/}
\def\bhs/{basic \hgeom/ series}

\let\goodbm\relax  \let\mmgood\relax \let\mmline\relax
\let\nline\nl \def\vvm#1>{\relax} \let\qqm\qqq

\ifMag\let\goodbm\goodbreak  \let\mmgood\vvgood \let\mmline\nl
 \let\vvm\vv \let\qqm\qquad
 \let\goodbreak\relax \let\vgood\relax \let\vvgood\relax \let\nline\relax \fi

\setparindent

\let\fixedpage\newpage

\whattime\readldf\writeldf

\csname hgeom.def\endcsname

\labeldef{S} {2} {2}
\labeldef{F} {2\labelsep \labelspace 1}  {gener}
\labeldef{F} {2\labelsep \labelspace 2}  {act}
\labeldef{F} {2\labelsep \labelspace 3}  {ect}
\labeldef{F} {2\labelsep \labelspace 4}  {order}
\labeldef{F} {2\labelsep \labelspace 5}  {xtry}
\labeldef{F} {2\labelsep \labelspace 6}  {Rests}
\labeldef{F} {2\labelsep \labelspace 7}  {xyz}
\labeldef{F} {2\labelsep \labelspace 8}  {hint0}
\labeldef{L} {2\labelsep \labelspace 1}  {well}
\labeldef{F} {2\labelsep \labelspace 9}  {list}
\labeldef{F} {2\labelsep \labelspace 10} {hint1}
\labeldef{F} {2\labelsep \labelspace 11} {Phi}
\labeldef{F} {2\labelsep \labelspace 12} {hpair}
\labeldef{F} {2\labelsep \labelspace 13} {Imap}
\labeldef{F} {2\labelsep \labelspace 14} {CD}
\labeldef{L} {2\labelsep \labelspace 5}  {Resint}
\labeldef{F} {2\labelsep \labelspace 15} {Spair}
\labeldef{F} {2\labelsep \labelspace 16} {Sh}
\labeldef{L} {2\labelsep \labelspace 7}  {Sp}
\labeldef{L} {2\labelsep \labelspace 8}  {Sq}
\labeldef{L} {2\labelsep \labelspace 9}  {main}

\labeldef{S} {3} {3}
\labeldef{F} {3\labelsep \labelspace 1}  {wW}
\labeldef{F} {3\labelsep \labelspace 2}  {W'}
\labeldef{L} {3\labelsep \labelspace 1}  {wbasis}
\labeldef{L} {3\labelsep \labelspace 2}  {wbasis'}
\labeldef{L} {3\labelsep \labelspace 3}  {Wbasis}
\labeldef{F} {3\labelsep \labelspace 3}  {alW}
\labeldef{L} {3\labelsep \labelspace 5}  {ResW}
\labeldef{F} {3\labelsep \labelspace 4}  {vv}
\labeldef{F} {3\labelsep \labelspace 5}  {vv*}
\labeldef{L} {3\labelsep \labelspace 7}  {main2}

\labeldef{S} {4} {4}
\labeldef{F} {4\labelsep \labelspace 1}  {KL}
\labeldef{F} {4\labelsep \labelspace 2}  {Mm}
\labeldef{L} {4\labelsep \labelspace 1}  {qKZ}
\labeldef{L} {4\labelsep \labelspace 2}  {qKZ'}
\labeldef{L} {4\labelsep \labelspace 3}  {limI}
\labeldef{L} {4\labelsep \labelspace 4}  {limI'}
\labeldef{L} {4\labelsep \labelspace 5}  {limIb}

\labeldef{S} {5} {P}
\labeldef{F} {5\labelsep \labelspace 1}  {Flm}
\labeldef{L} {5\labelsep \labelspace 1}  {mu}
\labeldef{F} {5\labelsep \labelspace 2}  {KK}
\labeldef{F} {5\labelsep \labelspace 3}  {TF}

\labeldef{S} {6} {R}
\labeldef{F} {6\labelsep \labelspace 1}  {MIIN}
\labeldef{F} {6\labelsep \labelspace 2}  {gene}
\labeldef{L} {6\labelsep \labelspace 1}  {welr}
\labeldef{F} {6\labelsep \labelspace 3}  {lisr}
\labeldef{L} {6\labelsep \labelspace 2}  {main3}
\labeldef{F} {6\labelsep \labelspace 4}  {MIINr}
\labeldef{F} {6\labelsep \labelspace 5}  {NIIM}

\labeldef{S} {7} {E}
\labeldef{F} {7\labelsep \labelspace 1}  {II}
\labeldef{L} {7\labelsep \labelspace 1}  {detIF}

\labeldef{S} {\char 65} {S}
\labeldef{L} {\char 65\labelsep \labelspace 1}  {thi}
\labeldef{F} {\char 65\labelsep \labelspace 1}  {Glg}
\labeldef{F} {\char 65\labelsep \labelspace 2}  {dnm}
\labeldef{L} {\char 65\labelsep \labelspace 3}  {detQ}
\labeldef{F} {\char 65\labelsep \labelspace 3}  {Xi}
\labeldef{L} {\char 65\labelsep \labelspace 4}  {detS}
\labeldef{L} {\char 65\labelsep \labelspace 5}  {combi}

\labeldef{S} {\char 66} {D}
\labeldef{F} {\char 66\labelsep \labelspace 1}  {glg}
\labeldef{L} {\char 66\labelsep \labelspace 1}  {detX}
\labeldef{L} {\char 66\labelsep \labelspace 2}  {detIW}
\labeldef{L} {\char 66\labelsep \labelspace 3}  {detIG}

\labeldef{S} {\char 67} {J}
\labeldef{F} {\char 67\labelsep \labelspace 1}  {fPT}
\labeldef{F} {\char 67\labelsep \labelspace 2}  {TA}
\labeldef{L} {\char 67\labelsep \labelspace 1}  {AJ}

\ifamsppt\NoRunningHeads\fi

\document

\ifMag\ifamsppt\else\osakatrue\fi\fi

\Sno 0

\ifosaka
\hfuzz 6pt
\center
\=
\vp1
\vsk-.5>
{\bls 16pt
\Bbf
Bilinear identity for {\Bmmi q}-\hint/s
\par}
\vsk2>
{\smc Vitaly TARASOV}
\vsk>
\endcenter
\else
\vp1\vsk->{\ifMag\eightpoint\fi\vsk-3>
\rline{{\sl Osaka J.\ Math.}\ (1998)}\rline{\sl to appear\/}\vsk>}
\center
\=
{\bls 16pt
\Bbf
Bilinear Identity for {\Bmmi q}-Hypergeometric Integrals
\par}
\vsk1.5>
\VT/
\vsk1.5>
{\it
\ifMag
\DMO/\\ \DMOaddr/
\else
\DMO/, \DMOaddr/
\fi}
\vsk1.8>
{\sl August \,19, 1997}
\endcenter
\ftext{\=\bls11pt \absence/\vv-.05>\nl
{\tenpoint\sl E-mail\/{\rm:} \dmoemail/\,, \,\homemail/}}
\vsk.2>
\fi

\vsk0>

\Sect{Introduction}

In this paper we describe a \bid/ satisfied by certain multidimensional
\qhint/s.  We give two \eqv/ forms of the identity, \cf. Theorems \[main]
and \[main2].  We call this identity the \em{\hRid/}. In the \onedim/ case
the \qhint/s can be expressed in terms of the \bhs/ $\phifn$ and the \hRid/
reduces to bilinear identities for these series, see Section~\[:E].
\par
There are two interpretations of the \hRid/ which make it a subject of our
interest. First, the \hRid/ is an analogue of the Riemann bilinear relation
for the twisted (co)homology groups defined by the reciprocally dual \loc/s,
see \Cite{CM}, \Cite{M}. An appropriate form of the \hRid/ is given by
Theorem~\[main]. From this point of view, the \tri/ and \ehgf/s correspond
\resp/ to the top cohomology and homology groups, the \hpair/s provide natural
duality between them and the \Spair/s play the role of the intersection forms.
This analogy can be seen directly from the explicit formulae. But in fact,
there is much deeper similarity, including the dual discrete \loc/s, the \dif/
twisted de~Rham complex etc., see \Cite{A}, \Cite{TV1}. The deformation of the
Riemann bilinear relation for the hyperelliptic Riemann surfaces was obtained
in \Cite{S}.
\par
The second interpretation of the \hRid/ comes from the \rep/ theory of \qaff/s,
namely, via the quantum \KZv/ (\qKZ/) \eq/. The \qKZe/ is a remarkable system
of \deq/s introduced in \Cite{FR}. It is a natural deformation of the famous
\difl/ \KZv/ \eq/ inheriting many of its nice properties.
\par
In \Cite{TV1} Varchenko and the author constructed all \sol/s of the \qKZe/
with values in a tensor product of $\Uu\!$ \Vmod/s with generic \hw/s
using the \qhint/s. The space of \sol/s of the \qKZe/ was indentified in
\Cite{TV1} with the tensor product of the corresponding \eVmod/s over the \eqg/
$\Eqg$.  The dual version of that construction is also available: the system
of \deq/s is dual to the \qKZe/, its \sol/s take values in the dual space to
the tensor product of \Umod/s and the space of \sol/s can be identified with
the dual space to the tensor product of \Emod/s. In this context the \hRid/
means that the \qhgeom/ \sol/s of the \qKZ/ and dual \qKZe/s transform
the natural pairing of the spaces of \sol/s to the natural pairing of
the target spaces; namely, given respective \sol/s $\Psi$ and $\Psi^*$
of the \qKZ/ and dual \qKZe/s we have that
$$
\bra\Val\Psi^*\},\Val\Psi\ket_{\,\text{target spaces}}\;=\,
\bra\Psi^*\},\Psi\ket_{\,\text{spaces of solutions}}\,.
$$
In particular, we can say that the \hRid/ is a deformation of
Gaudin-Korepin's formula for norms of the \Bv/s \Cite{KBI}, \cf. \Cite{TV3}.
\par
In the main part of the paper we make no references to the deformed
(co)homology theory and limit ourselves to only a few remarks about
the \rep/ theory; for more conceptual point of view see \Cite{TV1}, \Cite{TV2}.
\par
The results of \Cite{TV1} are crucial for our proof of the \hRid/.
It is shown in \Cite{TV1} that a matrix formed by certain \hint/s satisfies
a system of linear \deq/s, \cf. Theorem~\[qKZ] in this paper, which is \eqv/
to the \qKZe/. Moreover, this matrix has a finite limit in a suitable \azo/
and the limit is a triangular matrix, \cf. Propositions~\[limI], \[limI'].
To prove the \hRid/ we first use the system of \deq/s and its dual system,
and show that it is enough to check the identity as the parameters tends to
limit in the \azo/. The remaining verification is quite straightforward
since all the matrices involved have constant triangular \as/s in the \azo/.
\par
In this paper we consider the so-called \em{\tri/} or \em{multiplicative} case,
which involves the \qhint/s. There is the \rat/ version of the story which
involves the multidimensional \hint/s of \MB/ type, formulae being written in
terms of the gamma \fn/, \rat/ and \tri/ \fn/s. The \rat/ case can be
considered as the logarithmic deformation taking the intermediate place between
the classical and \tri/ deformed cases. The \rat/ case of the \qKZe/ and
related deformed (co)homologies were studied in \Cite{TV2}. The \rat/ version
of the \hRid/ can be obtained similarly to the \tri/ one using results of
\Cite{TV2}. It will appear separately.
\par
The paper is organized as follows. In Section~\[:2] we give definitions and
formulate the main result of the paper -- the \hRid/, \cf.  Theorem \[main].
The \eqv/ form of the identity is given in Section~\[:3],
\cf. Theorem~\[main2]. In Section~\[:4] we describe a system of \deq/s
satisfied by the \qhint/s and study their behaviour in a suitable \azo/.
The final step of the proof of the \hRid/ is made in Section~\[:P].
We give the restricted version of the \hRid/ in Section~\[:R] and
consider the \onedim/ example in Section~\[:E].
\par
There are three Appendices in the paper. Appendix~\[:S] contains the proof of
Proposition \[Sq]. Some of the determinant formulae from \Cite{TV1} relevant to
this paper are reproduced in Appendix~\[:D]. In Appendix~\[:J] we explain that
the \qhint/s which we are using in the paper essentially coincide with
the \sym/ \Atype/ Jackson integrals.
\Par
The author thanks K\&Aomoto, A\&Nakayashiki and \Varch/ for helpful
discussions.

\ifosaka\else\ifMag\fixedpage\fi\fi
\Sect[2]{The \hRid/}
\subsect*{Basic notations}
Let ${\Cx\!=\C\setminus\lb 0\rb}$. Fix a nonzero complex number $p$
\st/ ${|p|<1}$. Set
\ifMag
$$
\pZ\>=\>\lb\)p^s\vert s\in\Z\)\rb\,.
$$
\else
\vv.5>
${\pZ=\lb\)p^s\vert s\in\Z\)\rb}$.
\fi
Let $(u)\9=(u;p)\9=\topsmash{\prod_{s=0}^\8(1-p^su)}$ and let
\vvm-.35>
${\tht(u)=(u)\9(p\)u\1)\9(p)\9}$ be the Jacobi theta-\fn/.
\Par
Fix a nonnegative integer $\ell$. Take nonzero complex numbers $\eta$, $\xn$,
$\yn$ called \em{parameters}. Say that the parameters are generic if for any
$\roll$ and any $k,\mn$, we have
\ifMag
$$
\gather
\qquad\eta^{r+1}\nin\pZ,\qqq \eta^{\pm r}x_k/x_m\nin\pZ,\ \quad
\eta^{\pm r}y_k/y_m\nin\pZ\;\ \quad \for\ k\ne m\,,
\Tag{gener}
\\
\nn4>
\eta^{\pm r}x_k/y_m\nin\pZ.
\endgather
$$
\else
$$
\ \eta^{r+1}\nin\pZ,\qqq \eta^{\pm r}x_k/x_m\nin\pZ,\ \quad
\eta^{\pm r}y_k/y_m\nin\pZ\;\ \quad \for\ k\ne m\,,\qqq
\eta^{\pm r}x_k/y_m\nin\pZ.
\Tag{gener}
$$
\fi
All over the paper we assume that the parameters are generic, unless otherwise
stated.
\Par
For any \fn/ $f(\tell)$ and any \perm/ $\si\in\Sl$ set
$$
\align
[\)f\)]\vpp\si(\tell)\, &{}=\,f(\tsill)
\prabsi\;{t_{\si_b}-\eta\)t_{\si_a}\over\eta\)t_{\si_b}-t_{\si_a}}\;,
\Tag{act}
\\
\nn2>
\lbc\)f\)\rbc\vpp\si(\tell)\, &{}=\,f(\tsill)
\prabsi{\eta\)\tht(\eta\1t_{\si_b}/t_{\si_a})\]\over
\tht(\eta\)t_{\si_b}/t_{\si_a})}\;.
\Tag{ect}
\endalign
$$
Each of the formulae defines an action of the \symg/ $\Sl$.
\Par
For any ${\lg=(\lg_1\lc\lg_n)\in\Zpn}\,$ set $\,{\lg^m=\lg_1\lsym+\lg_m}$,
$\mn$. Set
$$
\Zln\,=\,\lb\>\lg\in\Zpn\vert\lg^n=\ell\>\rb\,.
$$
For any $\lg,\mg\in\Zln$, $\lg\ne\mg$, say that
$$
\lg\ll\mg\qquad\ \text{if}\qquad \lg^k\le\mg^k\qquad
\text{for any}\quad \kn-1\,.
\Tag{order}
$$
\Rem
We can identify $\lg\in\Zpn$ with a partition ${\lg^n\}\lsym\ge\lg^1}$.
The introduced above ordering on $\Zln$ coincides with the inverse dominance
ordering for the corresponding partitions.
\enddemo
For any $x\in\Cn\!$ and $\lg\in\Zln$ we define the point
${\xt\lg\>[\eta\)]\in\Cxl}\!$ as follows:
$$
\xt\lg\>[\eta\)]\,=\,(\eta^{1-\lg_1}x_1,\>\eta^{\)2-\lg_1}x_1\lc x_1,\>
\eta^{1-\lg_2}x_2\lc x_2,\,\ldots\,,\>\eta^{1-\lg_n}x_n\lc x_n)\,,
\ifMag\kern-1.3em\fi
\Tag{xtry}
$$
\par
For any \fn/ $f(\tell)$ and a point ${\ts\}=(\ts_1\lc\ts_\ell)}$ we define
the multiple residue $\Res f(t)\vst{t=\ts}$ by the formula
$$
\Res f(t)\vst{t=\ts}\)=\,\Res\bigl(\,\ldots\,
\Res f(\tell)\vst{t_\ell=\ts_\ell}\ldots\,\bigr)\vst{t_1=\ts_1}\,.
\Tag{Rests}
$$
We often use in the paper the following compact notations:
$$
t=(\tell)\,,\qqq x=(\xn)\,, \qqq y=(\yn)\,.
$$
\par
For any vector space $V$ we denote by ${V^*\]}$ the dual vector space,
and for a linear operator $A$ we denote by ${A^*\]}$ the dual operator.
\par
In this paper we extensively use results from \Cite{TV1}. We have the following
correspondence of parameters $\xn$, $\yn$ in this paper and parameters $\xin$,
$\zn$ in \Cite{TV1}\):
$$
x_m=\xi_mz_m\,,\qqq y_m=\xi_m\1z_m\,,\qqq\mn\,.
\Tag{xyz}
$$

\subsect{The \hint/}
Let $\Pht(t;x;y;\eta)$ be the following \fn/:
$$
\Pht(t;x;y;\eta)\,=\,\pron\,\pral\,{1\over(x_m/t_a)\9\>(t_a/y_m)\9\!}\;
\prod_{\tsize{a,b=1\atop a\ne b}}^\ell{1\;\over(\eta\1t_a/t_b)\9\!}\;.
$$
For any \fn/ $f(\tell)$ \hol/ in $\Cxl$ we define below the \em{\hint/}
${\IH(f\)\Pht)}$.
\par
Assume that $|\eta|>1$ and $|x_m|<1$, $|y_m|>1$, $\mn$. Then we set
$$
\IH[x;y;\eta\)](f\)\Pht)\,=\;{1\over(2\pii)^\ell\!}\,
\int_{\){\TT^\ell}\]}f(t)\>\Pht(t;x;y;\eta)\;\dtt
\Tag{hint0}
$$
where $\,\topsmash{\dtt=\}\pral dt_a/t_a}\,$ and
$\,{\TT^\ell=\lb\)t\in\Cl\vert\ |t_1|=1\llc|t_\ell|=1\)\rb}$. We define
${\IH[x;y;\eta\)](f\)\Pht)}$ for arbitrary values of the parameters $\eta$,
$\xn$, $\yn$ by the \anco/ \wrt/ the parameters.
\Prop{well}
For generic values of the parameters $\eta$, $\xn$, $\yn$, see \(gener),
the \hint/ ${\IH[x;y;\eta\)](f\)\Pht)}$ is well defined and is a \hof/ of
the parameters.
\endpro
\Pf.
For generic values of the parameters $\eta$, $\xn$, $\yn$ singularities
of the integrand ${f(t)\>\Pht(t;x;y;\eta)}$ are at most at the following
hyperplanes:
$$
t_a=0\,,\qquad t_a=p^sx_m\,,\qquad t_a=p^sy_m\,,\qquad t_a=p^{-s}\eta\)t_b\,,
\Tag{list}
$$
${a,b=1\lc\ell}$, ${a\ne b}$, $\,{\mn}$, $\,{s\in\Zp}$. The number of edges
(nonempty intersections of the hyperplanes) of configuration \(list) and
dimensions of the edges are always the same for nonzero generic values of the
parameters.  Therefore, the topology of the complement in $\Cl$ of the union of
the hyperplanes \(list) does not change if the parameters are nonzero generic.
\par
The rest of the proof is similar to the proof of Theorem 5.7 in \Cite{TV2}.
\epf
It is clear from the proof of Proposition \[well] that for generic value of
the parameters the \hint/ ${\IH[x;y;\eta\)](f\)\Pht)}$ can be represented as
an integral
$$
\IH[x;y;\eta\)](f\)\Pht)\,=\;{1\over(2\pii)^\ell\!}\,
\int_{\,\TTt^\ell[x;y;\eta\)]\!}\!\!\!f(t)\>\Pht(t;x;y;\eta)\;\dtt
\Tag{hint1}
$$
where ${\TTt^\ell[x;y;\eta\)]}$ is a suitable deformation of the torus
$\TT^\ell$ which does not depend on $f$.
\Rem
In what follows we are using the \hint/s ${\IH(f\)\Pht)}$ only for \sym/ \fn/s
$f$ which have a certain particular form. In this case the \hint/s coincide
with the \sym/ \Atype/ Jackson integrals, \cf. Appendix~\[:J].
\enddemo

\subsect{The \hgf/s and the \hpair/}
Let ${\Fc[x;\eta\);\ell\>]}$ be the space of \raf/s $f(\tell)$ \st/ the product
$$
f(\tell)\,\tpral t_a\1\,\pron\,\pral\,(t_a-x_m)
\prab\,{\eta\)t_a-t_b\over t_a-t_b}
$$
is a \sym/ \pol/ of degree less than $n$ in each of the \var/s $\tell$.
Elements of the space ${\Fc[x;\eta\);\ell\>]}$ are invariant \wrt/
action \(act) of the \symg/ $\Sl$. Set
$$
\Fc'[x;\eta\);\ell\>]\,=\,\bigl\lb\)f(\tell)\vert
t_1\ldots t_\ell\>f(\tell)\in \Fc[x;\eta\1;\ell\>]\)\bigr\rb\,.
$$
The spaces $\Fc$ and $\Fc'$ are called the \em{\thgf/s}.
\Rem
We have quite a few motivations to call the spaces $\Fc$ and $\Fc'$ \tri/,
though no \tri/ \fn/s will appear actually in the paper.
\enddemo
Fix $\al\in\Cx\!$. Let ${\Fq[\al;x;\eta\);\ell\>]}$ be the space of \fn/s
$g(\tell)$ \st/
$$
g(\tell)\,\pron\,\pral\,\tht(t_a/x_m)
\prab\,{\tht(\eta\)t_a/t_b)\over\tht(t_a/t_b)}
$$
is a \sym/ \hof/ of $\tell$ in $\Cxl\!$ and
$$
g(\tpell)\,=\,\al\)\eta^{2-2a}\>g(\tell)\,.
$$
Elements of the space ${\Fq[\al;x;\eta\);\ell\>]}$ are invariant \wrt/
action \(ect) of the \symg/ $\Sl$. Set
$$
\Fq'[\al;x;\eta\);\ell\>]\,=\,\Fq[\al\1;x;\eta\1;\ell\>]\,.
$$
The spaces $\Fq$ and $\Fq'$ are called the \em{\ehgf/s}.
\Rem
The parameter $\al$ here is related to the parameter $\ka$ in \Cite{TV1}\):
\mmline
$\,\al=\ka\>\eta^{\ell-1}\!\!\pron\!\xi_m\1\!$, \ \cf. \(xyz).
\enddemo
In what follows we do not indicate explicitly all arguments for the \hgf/s and
related maps if it causes no confusion. The suppressed arguments are supposed
to be the same for all the spaces and maps involved.
\Rem
The \thgf/s can be considered as degenerations of the \ehgf/s as $p\to 0$ and
then $\al\to 0$. In this limit the spaces $\Fq[\al]$ and $\Fq'[\al]$ degenerate
into the spaces $\Fc$ and $\Fc'\!$, \resp/. Because of this correspondence we
use two slightly different versions of the \thgf/s.
\enddemo
\Prop{}
\back\Cite{TV1}
For any $\al,\,\eta,\,\xn$ we have that
$$
\dim\Fc[x;\eta\);\ell\>]\,=\,\dim\Fc'[x;\eta\);\ell\>]\,=\,
\dim\Fq[\al;x;\eta\);\ell\>]\,=\,{n+\ell-1\choose n-1}\,.
$$
\endpro
Let $\Phi(t;x;y;\eta)$ be the following \fn/:
$$
\Phi(t;x;y;\eta)\,=\,\pron\,\pral\,{(t_a/x_m)\9\!\over(t_a/y_m)\9}\,
\prab\;{(\eta\)t_a/t_b)\9\!\over(\eta\1t_a/t_b)\9\!}\;.
\Tag{Phi}
$$
We call the \fn/ $\Phi(t;x;y;\eta)$ the \em{\phf/}. Notice that
\ifMag
$$
\Phi(t;y;x;\eta\1)\,=\,\bigl(\Phi(t;x;y;\eta)\bigr)\1.
$$
\par
\else
$\;{\Phi(t;y;x;\eta\1)\,=\,\bigl(\Phi(t;x;y;\eta)\bigr)\1}\!$.
\Par
\fi
The \hint/ \(hint1) induces the \em{\hpair/s} of the \tri/ and \ehgf/s:
\ifMag
$$
\alignat2
\Llap{I\)[\al;x;y;\eta\)]\):\)\Fq} &
\Rlap{[\al;x;\eta\)]\ox\Fc[x;\eta\)]\,\to\,\C\,,} &&
\Tag{hpair}
\\
\nn4>
\Llap{f\ox g\>\map\,{1\over\ell\)!}\IH{}} &
\Rlap{[x;y;\eta\)]\bigl(fg\>\Phi(\cdot;x;y;\eta)\bigr)\,,} &&
\\
\nn8>
\ \ I'[\al;x;y;\eta\)]\):\)\Fq' &
[\al;y;\eta\)] && {}\ox\Fc'[y;\eta\)]\,\to\,\C\,,
\\
\nn4>
&&& \Llap{f\ox g\>\map\,{1\over\ell\)!}\IH[y;x;\eta\1]}
\bigl(fg\>\Phi(\cdot;y;x;\eta\1)\bigr)\,.
\endalignat
$$
\else
$$
\alignat2
\quad I\)[\al;x;y;\eta\)]\):\)
\Fq &[\al;x;\eta\)]\ox\Fc[x;\eta\)]\,\to\,\C\,, \qqq &
I'[\al;x;y;\eta\)]\):\)\Fq'[\al;y;\eta\)] &{}\ox\Fc'[y;\eta\)]\,\to\,\C\,,
\kern-1em
\Tag{hpair}
\\
\nn6>
f\ox g\>\map\,{1\over\ell\)!}
\IH &[x;y;\eta\)]\bigl(fg\>\Phi(\cdot;x;y;\eta)\bigr)\,, \qqq &
f\ox g\>\map\,{1\over\ell\)!}
\IH[y;x;\eta\1] & \bigl(fg\>\Phi(\cdot;y;x;\eta\1)\bigr)\,.
\kern-1em
\endalignat
$$
\fi
\vgood
We also consider these pairings as linear maps from the \ehgf/s to the dual
spaces of the \thgf/s, denoting them by the same letters:
$$
\align
I\)&[\al;x;y;\eta\)]\):\)\Fq[\al;x;\eta\)]\,\to\,\bigl(\Fc[x;\eta\)]\bigr)^*,
\Tag{Imap}
\\
\nn4>
I'&[\al;x;y;\eta\)]\):\)\Fq'[\al;y;\eta\)]\,\to\,\bigl(\Fc'[y;\eta\)]\bigr)^*.
\endalign
$$
\Rem
In this paper we multiply the \hpair/s by an additional factor $\dsize\ifMag
{1\over(2\pii)^\ell\>\ell\)!}\else\smash{1\over(2\pii)^\ell\>\ell\)!}\fi$
compared with the \hpair/ in \Cite{TV1}.
\enddemo
\Prop{}
\back\Cite{TV1}
Let the parameters $\eta$, $\xn$, $\yn$ be generic.
\mmline
Assume that
$$
\al\nepsr\,,\qqm \al\>\alti\ne\)p^{-s-1}\eta^{-r}\,,
\qqm \roll\,,\quad s\in\Zp\,.
$$
Then the \hpair/ $\,{I\)[\al;x;y;\eta\)]\):\)
\Fq[\al;x;\eta\)]\,\to\,\bigl(\Fc[x;\eta\)]\bigr)^*}$ is \ndeg/.
\endpro
\nt
The statement follows from Corollary \[detIG].
\Cr{}
Let the parameters $\eta$, $\xn$, $\yn$ be generic. Assume that
$$
\al\ne\)p^{-s-1}\eta^r\,,\qqm \al\>\alti\ne\)p^s\eta^{-r}\,,
\qqm \roll\,,\quad s\in\Zp\,.
$$
Then the \hpair/ $\,{I'[\al;x;y;\eta\)]\):\)
\Fq'[\al;y;\eta\)]\,\to\,\bigl(\Fc'[y;\eta\)]\bigr)^*}$ is \ndeg/.
\endpro
\Pf.
Let $\pi$ be the following map:
$\,{\pi:f(\tell)\map t_1\ldots t_\ell\>f(\tell)}$. Then the next diagram
is commutative:
$$
\CD
\Fq[p\1\al\1;y;\eta\1]\; @>{\tsize\ \;I\)[p\1\al\1;y;x;\eta\1]\ \;}>>
\,\bigl(\Fc[y;\eta\1]\bigr)^*
\\
@V{\tsize \pi\ }VV
@VV{\tsize\ \pi^*}V
\\
\Fq'[\al;y;\eta\)]\;
@>>{\tsize\ \;\Cph{I\)[p\1\al\1;y;x;\eta\1]}{I'[\al;x;y;\eta\)]}\ \;}>
\,\bigl(\Fc'[y;\eta\)]\bigr)^*
\endCD
\Tag{CD}
$$
and the vertical arrows are invertible, which proves the statement.
\epf

\subsect{The \Spair/}
Let the points ${\xt\lg\>[\eta\)]\in\Cxl}\!$, $\lg\in\Zln$, be defined by
\(xtry). For any \fn/ $f(\tell)$ set
$$
\Resi[x;\eta\)](f)\,=
\smiZ\Res\bigl(t_1\1\!\ldots t_\ell\1f(\tell)\bigr)\vst{t=\xt\mg[\eta]}\,.
$$
\Lm{Resint}
\back\Cite{TV1}
For any $f\in\Fc[x;\eta;\ell]$ and $g\in\Fc'[y;\eta;\ell]$ we have
$$
\Resi[x;\eta\)](fg)\,=\,(-1)^\ell\}\Resi[y;\eta\1](fg)\,=\;{1\over(2\pii)^\ell
\>\ell\)!}\,\int_{\,\TTt^\ell[x;y;\eta\)]\!}\!\!\!f(t)\>g(t)\;\dtt
$$
where ${\TTt^\ell[x;y;\eta\)]}$ is the deformation of the torus $\TT^\ell$
defined by \(hint1).
\endpro
\Lm{}
\back\Cite{TV1}
For any $f\in\Fq[\al;x;\eta;\ell]$ and $g\in\Fq'[\al;y;\eta;\ell]$ we have
$$
\Resi[x;\eta\)](fg)\,=\,(-1)^\ell\}\Resi[y;\eta\1](fg)\,.
$$
\endpro
We define the \em{\Spair/s} of the \tri/ and \ehgf/s as follows:
\ifMag
$$
\align
S\)[x;y;\eta\)]\):\)\Fc' & [y;\eta;\ell]\ox\Fc[x;\eta;\ell]\,\to\,\C\,,
\Tag{Spair}
\\
\nn6>
f\ox g\;& \]{}\map\,\Resi[x;\eta\)](fg)\,,
\\
\nn8>
\Sq[\al;x;y;\eta\)]\):\)\Fq' &[\al;y;\eta;\ell]\ox\Fq[\al;x;\eta;\ell]
\,\to\,\C\,,
\\
\nn6>
f\ox g\;& \]{}\map\,\Resi[x;\eta\)](fg)\,.
\endalign
$$
\else
$$
\alignat2
\quad S\)[x;y;\eta\)]\):\)\Fc' &[y;\eta;\ell]\ox\Fc[x;\eta;\ell]\,\to\,\C\,,
\qquad &
\Sq[\al;x;y;\eta\)]\):\)\Fq' &[\al;y;\eta;\ell]\ox\Fq[\al;x;\eta;\ell]
\,\to\,\C\,,\kern-2em
\Tag{Spair}
\\
\nn6>
f\ox g\;& \]{}\map\,\Resi[x;\eta\)](fg)\,, &
f\ox g\;& \]{}\map\,\Resi[x;\eta\)](fg)\,.
\endalignat
$$
\fi
\vgood
We also consider these pairings as linear maps, denoting them by
the same letters:
\ifMag
$$
\gather
S\)[x;y;\eta\)]\):\)\Fc'[y;\eta\)]\,\to\,\bigl(\Fc[x;\eta\)]\bigr)^*,
\Tag{Sh}
\\
\nn4>
\Sq[\al;x;y;\eta\)]\):\)\Fq'[\al;y;\eta\)]\,\to\,
\bigl(\Fq[\al;x;\eta\)]\bigr)^*.
\endgather
$$
\else
$$
S\)[x;y;\eta\)]\):\)\Fc'[y;\eta\)]\,\to\,\bigl(\Fc[x;\eta\)]\bigr)^*,
\qqq
\Sq[\al;x;y;\eta\)]\):\)\Fq'[\al;y;\eta\)]\,\to\,
\bigl(\Fq[\al;x;\eta\)]\bigr)^*.
\Tag{Sh}
$$
\fi
\Prop{Sp}
Let the parameters $\eta$, $\xn$, $\yn$ be generic. Then the \Spair/
${S\)[x;y;\eta\)]\):\)\Fc'[y;\eta\)]\,\to\,\bigl(\Fc[x;\eta\)]\bigr)^*}$
is \ndeg/.
\endpro
\nt
The statement follows from Lemma \[ResW].
\Prop{Sq}
Let the parameters $\eta$, $\xn$, $\yn$ be generic. Assume that
$$
\al\>\eta^{-r}\nin\pZ,\qqq \al\>\eta^{r+2-2\ell}\xy\nin\pZ,\qqq\roll\,.
$$
Then the \Spair/ ${\Sq[\al;x;y;\eta\)]\):\)\Fq'[\al;y;\eta\)]\,\to\,
\bigl(\Fq[\al;x;\eta\)]\bigr)^*}$ is \ndeg/.
\endpro
\nt
The statement follows from Proposition \[detS].

\subsect{The \hRid/}
In this section we formulate the main result of the paper, the \hRid/ which
involves both the \hgeom/ and \Spair/s, see Theorem \[main]. We prove this
result in Section~\[:P].
\Th{main}
Let the parameters $\eta$, $\xn$, $\yn$ be generic.
Then the following diagram is commutative:
$$
\CD
\Fq[\al;x;\eta\)]\;
@>{\tsize\ \;\Cph{\bigl(I'[\al;x;y;\eta\)]\bigr)^*}{I\)[\al;x;y;\eta\)]}\ \;}>>
\,\bigl(\Fc[x;\eta\)]\bigr)^*
\\
@V{\tsize{\vrp18pt:12pt>}\bigl(\Sq[\al;x;y;\eta\)]\bigr)^*\ }VV
@VV{\tsize\ (-1)^\ell\>\bigl(S\)[x;y;\eta\)]\bigr)\1{\vrp18pt:12pt>}}V
\\
\bigl(\Fq'[\al;y;\eta\)]\bigr)^*\;
@<<{\tsize\ \;\bigl(I'[\al;x;y;\eta\)]\bigr)^*\ \;}<
\,\Fc'[y;\eta\)]
\endCD
\ifMag\ \;\fi
$$
\endpro
\Rem
Given bases of the \hgf/s, Theorem \[main] translates into bilinear relations
for the corresponding \hint/s. In the next section we describe an important
example of the bases --- the bases given by the \wtf/s \(wW), \(W').
\enddemo

\Sect[3]{Tensor coordinates on the \hgf/s}
In this section we give an \eqv/ form of the \hRid/, see Theorem \[main2].

\subsect{Bases of the \hgf/s}
For any $\lg\in\Zln$ define the \fn/s $w_\lg$ and $W_\lg$ by the formulae:
\ifMag
$$
\gather
w_\lg(t;x;y;\eta)\,=\>\pron\,\prod_{s=1}^{\lg_m}\>{1-\eta\over1-\eta^s}\,
\susi\,\Bigl[\,\pron\>\prod_{a\in\Gm_{\Rph lm}}\Bigl(\,
{t_a\over t_a-x_m}\ \prlm\,{t_a-y_l\over t_a-x_l}\,\Bigr)\>\Bigr]_\si \,,
\kern-2em
\Tag{wW}
\\
\nn10>
\aligned
\ W_\lg & (t;\al;x;y;\eta)\,={}
\\
\nn6>
& {}=\>\pron\,\prod_{s=1}^{\lg_m}\,{\tht(\eta)\over\tht(\eta^s)}\,
\susi\,\LBc\,\pron\,\prod_{a\in\Gm_{\Rph lm}}
\Bigl(\,{\tht(\eta^{\)2a-2}\al_m\1t_a/x_m)\over\tht(t_a/x_m)}
\prlm\,{\tht(t_a/y_l)\over\tht(t_a/x_l)}\,\Bigr)\>\RBc_\si\,\hp,\kern-2em
\endaligned
\endgather
$$
\else
$$
\gather
w_\lg(t;x;y;\eta)\,=\>\pron\,\prod_{s=1}^{\lg_m}\>{1-\eta\over1-\eta^s}\,
\susi\,\Bigl[\,\pron\>\prod_{a\in\Gm_{\Rph lm}}\Bigl(\,
{t_a\over t_a-x_m}\ \prlm\,{t_a-y_l\over t_a-x_l}\,\Bigr)\>\Bigr]_\si\,,
\kern-1em
\Tag{wW}
\\
\nn6>
W_\lg(t;\al;x;y;\eta)\,=\>
\pron\,\prod_{s=1}^{\lg_m}\,{\tht(\eta)\over\tht(\eta^s)}\,
\susi\,\LBc\,\pron\,\prod_{a\in\Gm_{\Rph lm}}
\Bigl(\,{\tht(\eta^{\)2a-2}\al_m\1t_a/x_m)\over\tht(t_a/x_m)}
\prlm\,{\tht(t_a/y_l)\over\tht(t_a/x_l)}\,\Bigr)\>\RBc_\si
\endgather
$$
\fi
where $\,{\Gm_m=\lb 1+\lg^{m-1}\,\lc\lg^m\rb}\,$ and
$\,\botsmash{\al_m=\al\!\}\prlm\!\!\}x_l/y_l}$, $\ \mn$. Set
\ifMag
$$
\gather
w'_\lg(t;x;y;\eta)\,=\)\tpron y_m^{\lg_m}\]
\tpral\}t_a\1\)w_\lg(t;y;x;\eta\1)\,,
\Tag{W'}
\\
\nn6>
W'_\lg(t;\al;x;y;\eta)\,=\,W_\lg(t;\al\1;y;x;\eta\1)\,.
\endgather
$$
\else
$$
w'_\lg(t;x;y;\eta)\,=\)\tpron y_m^{\lg_m}\]
\tpral\}t_a\1\)w_\lg(t;y;x;\eta\1)\,,
\qqq W'_\lg(t;\al;x;y;\eta)\,=\,W_\lg(t;\al\1;y;x;\eta\1)\,.
\Tag{W'}
$$
\fi
The \fn/s $w_\lg,\,w'_\lg$ and $W_\lg,\,W'_\lg$ are called the \em{\tri/} and
\em{\ewf/s}, \resp/.
\Rem
In this paper we use a slightly different normalization of the \twf/s compared
with \Cite{TV1}.
\enddemo
\Prop{wbasis}
\back\Cite{TV1}
Let the parameters $\eta$, $\xn$, $\yn$ be generic. Then the \fn/s
${\lb\)w_\lg(t;x;y;\eta)\)\rb\vpp{\lg\in\Zln}}$ form a basis in the \thgf/
${\Fc[x;\eta\)]}$.
\endpro
\nt
The statement follows from Proposition \[detX].
\Cr{wbasis'}
Under the above assumptions the \fn/s
${\lb\)w'_\lg(t;x;y;\eta)\)\rb\vpp{\lg\in\Zln}}$ form a basis in the \thgf/
${\Fc'[y;\eta\)]}$.
\endpro
\Prop{Wbasis}
\back\Cite{TV1}
Let the parameters $\eta$, $\xn$, $\yn$ be generic. Assume that
$$
\al\>\eta^{-r}\!\!\tsize\prllm\!x_l/y_l\nin\pZ,\qqq\mn-1\,,\quad\rolt\,.
\ifMag\kern-2.2em\fi
\Tag{alW}
$$
Then the \fn/s ${\lb\)W_\lg(t;\al;x;y;\eta)\)\rb\vpp{\lg\in\Zln}}$ form a basis
in the \ehgf/ ${\Fq[\al;x;\eta\)]}$.
\endpro
\nt
The statement follows from Proposition \[detQ].
\Cr{}
Under the above assumptions the \fn/s
${\lb\)W'_\lg(t;\al;x;y;\eta)\)\rb\vpp{\lg\in\Zln}}$ form a basis
in the \ehgf/ ${\Fq'[\al;x;\eta\)]}$.
\endpro
The next lemma shows that the bases ${\lb\)w_\lg\)\rb}$ and
${\lb\)w'_\lg\)\rb}$ of the \thgf/ are biorthogonal \wrt/ the \Spair/ \(Spair),
and the same holds for the bases ${\lb\)W_\lg\)\rb}$, ${\lb\)W'_\lg\)\rb}$ of
the \ehgf/s.
\Lm{ResW}
\vv->
$$
\gather
S\)(w'_\lg,w_\mg)\,=\,\dl_{\lg\)\mg}\>\pron\prod_{s=0}^{\lg_m-1}\,
{(1-\eta)\>\eta^s\)y_m\over(1-\eta^{s+1})\>(x_m-\eta^sy_m)}\;,
\\
\nn8>
\Sq(W'_\lg,W_\mg)\,=\,\dl_{\lg\)\mg}\>
\pron\prod_{s=0}^{\lg_m-1}\,{\eta^s\)\tht(\eta)\>
\tht(\eta^s\al_{\lg,m}\1)\>\tht(\eta^{1-s-\lg_m}\al_{\lg,m}\)x_m/y_m)
\over\tho\>\tht(\eta^{s+1})\>\tht(\eta^{-s}x_m/y_m)}
\endgather
$$
where $\;{\al_{\lg,m}=\al\!\prjm\!\eta^{-2\lg_j}x_j/y_j}\;$ and
$\;\dsize\tho=\){d\over du}\tht(u)\vst{u=1}\]=-\)(p)^3_\8$.
\endpro
\Pf.
The formulae are \resp/ \eqv/ to formulae (C.9) and (C.4) in \Cite{TV1}.
\epf

\subsect{The \tenco/ and the \hmap/s}
Let ${V=\PmiZ\C\>v_\mg}$ and let ${V^*=\PmiZ\C\>v^*_\mg}\ifMag\vp{\Plus^{}}\fi$
be the dual space. Denote by ${\bra\>{,}\>\ket}$ the canonical pairing:
${\bra v^*_\lg,v_\mg\ket=\dl_{\lg\)\mg}}$.
\Par
Introduce the \em{\tenco/} on the \hgf/s, \cf. \Cite{TV1}, \Cite{V}.
They are the following linear maps:
$$
\alignat2
B[{}& x;y;\eta\)]:V^*\]\to\>\Fc[x;\eta\)]\,, &
\Bq[{}& \al;x;y;\eta\)]:V^*\]\to\>\Fq[\al;x;\eta\)]\,,
\\
\nn4>
& v^*_\mg\,\map\,w_\mg(t;x;y;\eta)\,, &
& v^*_\mg\,\map\,W_\mg(t;\al;x;y;\eta)\,,
\\
\nn8>
B'[{}& x;y;\eta\)]:V^*\]\to\>\Fc'[y;\eta\)]\,, \qquad\qquad &
\Bq'[{}& \al;x;y;\eta\)]:V^*\]\to\>\Fq'[\al;y;\eta\)]\,,
\\
\nn4>
& v^*_\mg\,\map\,w'_\mg(t;x;y;\eta)\,, &
& v^*_\mg\,\map\,W'_\mg(t;\al;x;y;\eta)\,.\vrd2.75ex>
\endalignat
$$
Under the assumptions of Propositions \[wbasis] and \[Wbasis] the \tenco/
are \iso/s of the respective vector spaces.
\Rem
The \tenco/ used in this paper differ from the \tenco/ in \Cite{TV1} by
normalization factors.
\enddemo
The \tenco/ and the \Spair/s \(Sh) induce bilinear forms
$$
\gathered
(\>{,}\>)[x;y;\eta\)]\):V\ox V\)\to\>\C\,,
\\
\nn4>
(u,v)=\bigl\bra B'\vpb{-1}S\vpb{-1}B^*\vpb{-1} u,v\bigr\ket\,,
\endgathered
\qquad\qquad
\gathered
\lpt\>{,}\>\rpt[\al;x;y;\eta\)]\):V^*\]\ox V^*\]\to\>\C\,,
\\
\nn4>
\lpt u,v\rpt=(-1)^\ell\bigl\bra u,(\Bq)^*\Sq\)\Bq'\>v\bigr\ket\,.
\endgathered
$$
We omit the common arguments in the second line.
The explicit formulae for these pairings are:
\vvm->
$$
\gather
\ifMag\qquad\fi
(v_\lg,v_\mg)\,=\,\dl_{\lg\)\mg}\>\pron\prod_{s=0}^{\lg_m-1}\,
{(1-\eta^{s+1})\>(x_m-\eta^sy_m)\over(1-\eta)\>\eta^s\)y_m}\;,
\Tag{vv}
\\
\nn8>
\ifMag\qquad\fi
\lpt v^*_\lg,v^*_\mg\rpt\,=\,
\dl_{\lg\)\mg}\pron\prod_{s=0}^{\lg_m-1}\,{\eta^s\)\tht(\eta)\>
\tht(\eta^s\al_{\lg,m}\1)\>\tht(\eta^{1-s-\lg_m}\al_{\lg,m}\)x_m/y_m)
\over\>\tht(\eta^{s+1})\>\tht(\eta^{-s}x_m/y_m)\>(p)\9^3}
\Tag{vv*}
\endgather
$$
where $\;{\al_{\lg,m}=\al\!\prjm\!\eta^{-2\lg_j}x_j/y_j}$,
\;\cf. Lemma \[ResW]. These formulae imply the next proposition.
\Prop{}
Let the parameters $\eta$, $\xn$, $\yn$ be generic. Then the form
${(\>{,}\>)}$ is \ndeg/. The form ${\lpt\>{,}\>\rpt}$ is \ndeg/ provided that
$$
\gather
\ifMag\qquad\fi
\al\>\eta^{-r}\nin\pZ,\qqq \al\>\eta^{r+2-2\ell}\xy\nin\pZ,\qqq\roll\,,
\\
\nn-3>
\Text{and}
\nn-2>
\ifMag\qquad\fi
\al\>\eta^{-r}\!\!\tsize\prllm\!x_l/y_l\nin\pZ,\qqq\mn-1\,,\quad\rolt\,.
\endgather
$$
\endpro
\Rem
The space $V$ can be identified with a weight subspace in a tensor product of
$\Uu\!$ \Vmod/s, the form ${(\>{,}\>)}$ coinciding with the tensor product of
the corresponding $\Uu$ \Sform/s. The space $V$ also can be identified with
a weight subspace in a tensor product of \eVmod/s over the \eqg/ $\Eqg$.
The last space has a certain natural bilinear form which is an elliptic
analogue of the tensor product of the \Sform/s. The form ${\lpt\>{,}\>\rpt}$
on $V^*\!$ and the ``elliptic \Sform/'' on $V$ correspond to each other.
\enddemo
Consider the following linear maps:
$$
\gathered
\Ib\)[\al;x;y;\eta\)]\):V^*\]\to\>V\,,
\\
\nn4>
\Ib=B^*\)I\>\Bq\,,
\endgathered
\qquad\qquad
\gathered
\Ib'[\al;x;y;\eta\)]\):V^*\]\to\>V\,,
\\
\nn4>
\Ib'=(B')^*\)I'\)\Bq'\,,
\endgathered
$$
where $I$ and $I'$ are given by \(Imap) and we omit the common arguments in
the second line. We call $\Ib$ and $\Ib'$ the \em{\hmap/s}.
\Par
Theorem \[main] is \eqv/ to the following statement.
\Th{main2}
Let the parameters $\eta$, $\xn$, $\yn$ be generic. Then the \hmap/s
${\Ib\)[\al;x;y;\eta\)]}$, ${\Ib'[\al;x;y;\eta\)]}$ respect the forms
${(\>{,}\>)[x;y;\eta\)]}$,
\mmline
${\lpt\>{,}\>\rpt[\al;x;y;\eta\)]}$. That is,
for any $u,v\in V^*\!$ we have
$$
\lpt u,v\rpt[\al]\,=\,
\bigl(\Ib'[\al]\>u,\Ib\)[\al]\>v\bigr)\,.
$$
\endpro
\vsk-.5>
\vsk0>

\Sect[4]{Difference \eq/s and asymptotics}
In this section we describe a system of \deq/s satisfied by the \hmap/s and
\as/s of the \hmap/s in a suitable \azo/ of the parameters $\xn$, $\yn$.

\subsect{Difference \eq/s for the \hmap/s}
Let ${L_k[x;y;\eta\)]}$ and ${L'_k[x;y;\eta\)]}$, $\kn$, be linear operators
acting in the \thgf/s ${\Fc[x;y;\eta\)]}$ and ${\Fc'[x;y;\eta\)]}$, \resp/.
The ope\-rators are defined by their actions on the bases of the \twf/s:
$$
\align
L_k[\al;x;y;\eta\)]\>w_\lg(\cdot;x;y;\eta)\, &{}=\,
\bigl(\al\)\eta^{1-\ell}\xy\bigr)^{\lg^k}
w_{{\vp\lg}^k\}\lg}(\cdot;{}^k\]x;{}^k\]y;\eta)\,,
\\
\nn4>
L'_k[\al;x;y;\eta\)]\>w'_\lg(\cdot;x;y;\eta)\, &{}=\,
\bigl(\al\)\eta^{1-\ell}\xy\bigr)^{-\lg^k}
w'_{{\vp\lg}^k\}\lg}(\cdot;{}^k\]x;{}^k\]y;\eta)\,,
\endalign
$$
where ${^k\lg=(\lg_{k+1}\lc\lg_n,\lg_1\lc\lg_k)}$,
${^k\]x=(x_{k+1}\lc x_n,x_1\lc x_k)}$,
${^k\]y=(y_{k+1}\lc y_n,{}}\alb y_1\lc y_k)$.
\Par
Using the \tenco/ we introduce operators ${K_m,K'_m\in\End(V)}$,
\ifMag $\mn$\else ${\mn}$\fi, by the formulae:
$$
\align
K_m[\al;x;y;\eta\)]\, &{}=\,\bigl(\]\bigl(B[x;y;\eta\)]\bigr)\vpb{-1}
L_m[\al;x;y;\eta\)]\>B[x;y;\eta\)]\)\bigr)^*,
\Tag{KL}
\\
\nn4>
K'_m[\al;x;y;\eta\)]\, &{}=\,\bigl(\]\bigl(B'[x;y;\eta\)]\bigr)\1
L'_m[\al;x;y;\eta\)]\>B'[x;y;\eta\)]\)\bigr)^*.
\endalign
$$
We also define operators ${M_m[\al;x;y;\eta\)]\in\End(V^*)}$, $\mn$:
$$
M_m[\al;x;y;\eta\)]\>v^*_\lg\,=\,\mu_{\lg,m}[\al;x;y;\eta\)]\>v^*_\lg\,,\qqq
\mu_{\lg,m}\,=\,\bigl(\al\)\eta^{1-\lg^m}\!\!
\tprjjm\!x_j/y_j\bigr)^{-\lg^m}\!\!. \ifMag\kern-1em\fi
\Tag{Mm}
$$
\vgood
Let $T_m^h$, $\mn$, be the multiplicative shift operators acting on \fn/s of
$\xn$, $\yn$:
\ifMag
$$
\align
(T_m^h f) & (\xn;\yn)\,={}
\\
\nn4>
& {}=\,f(hx_1\lc hx_m,x_{m+1}\lc x_n;hy_1\lc hy_m,y_{m+1}\lc y_n)\,.
\endalign
$$
\else
$$
(T_m^h f)(\xn;\yn)\,=\,
f(hx_1\lc hx_m,x_{m+1}\lc x_n;hy_1\lc hy_m,y_{m+1}\lc y_n)\,.
$$
\fi
Set \,$T_m=T_m^p$, \,$\mn$.
\Th{qKZ}
\back\Cite{TV1}
The \hmap/ ${\Ib\)[\al;x;y;\eta\)]}$ satisfies the following system of \deq/s:
$$
T_m\Ib\)[\al;x;y;\eta\)]\,=\,K_m[\al;x;y;\eta\)]\>
\Ib\)[\al;x;y;\eta\)]\>M_m[\al;x;y;\eta\)]\,,\qqq\mn\,.
$$
\endpro
\Cr{qKZ'}
The \hmap/ ${\Ib'[\al;x;y;\eta\)]}$ satisfies the following system of \deq/s:
$$
T_m\Ib'[\al;x;y;\eta\)]\,=\,K'_m[\al;x;y;\eta\)]\>
\Ib'[\al;x;y;\eta\)]\>\bigl(M_m[\al;x;y;\eta\)]\bigr)\1\,,\qqq\mn\,.
$$
\endpro
\nt
The last claim results from the commutativity of diagram \(CD) and
formulae \(W').
\Rem
The numbers $\mu_{\lg,m}$ are related to the transformation properties of
the \ewf/s:
$$
T_m W_{\lg}\>=\>\mu_{\lg,m}\;\tprjjm(x_j/y_j)^\ell\,W_{\lg}\,,
\qqq T_m W'_{\lg}\>=\>
\mu_{\lg,m}\1\;\tprjjm(x_j/y_j)^{-\ell}\,W'_{\lg}\,.
$$
\enddemo
\Rem
The system of \deq/s ${T_m\Psi\)=K_m\Psi}$, $\mn$, can be identified with
the \qKZe/ with values in a weight subspace ${({}\}=V\))}$ of a tensor product
of $\Uu\!$ \Vmod/s. Its \sol/s have the form ${\Psi=\Ib\>Y}$, where
${Y\in\End(V^*)}$ solves the system of \deq/s ${T_mY=M_m\1\)Y}$. Notice that
the operators $M_1\lc M_n$ are invariant \wrt/ the shift operators
$T_1^h\lc T_n^h$ for any nonzero $h$. The factor $Y$ plays the role
of an adjusting map in \Cite{TV1}.
\par
The system ${T_m\Psi'\)=K_m'\Psi'}$, $\mn$, corresponds to the dual \qKZe/,
if we identify the spaces $V$ and $V^*\!$ using the \Sform/ ${(\>{,}\>)}$.
\enddemo

\subsect{Asymptotics of the \hmap/s}
Let $\AA$ be the following \azo/ of the parameters $\xn$, $\yn$:
$$
\AA\,=\,\left\lb\,
\alignedat2
& |x_m/x_{m+1}|\ll 1\,,\quad && \mn-1
\\
& |x_m/y_m|\simeq 1\,, && \mn
\endalignedat
\,\right\rb\,.
$$
We say that $(x;y)$ tends to limit in $\AA$ and write $(x;y)\rto\AA$ if
$$
x_m/x_{m+1}\>\to\,0\,,\qqq\mn-1\,,
$$
and the ratios $x_m/y_m$ and $y_m/x_m$ remain bounded for any $\mn$.
If a \fn/ $f(x;y)$ has a finite limit as $(x;y)\rto\AA$, we denote this limit
by $\limA f$. Notice that the limit $\limA f$ can depend on $\xn$, $\yn$, but
it is invariant \wrt/ the shift operators $T_1^h\lc T_n^h$ for any nonzero $h$.
\Rem
The operators $K_1\lc K_n$, $K'_1\lc K'_n$, \cf. \(KL), have finite limits as
\alb\ $(x;y)\rto\AA$, and the limits are \resp/ lower and upper triangular
\wrt/ the basis ${\lb\)v_\lg\)\rb_{\lg\in\Zln}}$ and ordering \(order). 
The diagonal parts of the limits $\limA K_m\1$ and $\limA K'_m$ are equal
to $M^*_m$.
\enddemo
\goodbreak
Define the \fn/s $I_{\lg\mg}(x,y)$ and $I'_{\lg\mg}(x,y)$ by the formulae:
$$
\align
I_{\lg\mg}(x,y)\,=\,
I\)& [\al;x;y]\bigl(W_\lg(\cdot;\al;x;y),w_\mg(\cdot;x;y)\bigr)\,,
\\
\nn6>
I'_{\lg\mg}(x,y)\,=\,
I'& [\al;x;y]\bigl(W'_\lg(\cdot;\al;x;y),w'_\mg(\cdot;x;y)\bigr)\,.
\endalign
$$
\Prop{limI}
\back\Cite{TV1}
For any $\lg,\mg\in\Zln$ the \hint/ $I_{\lg\mg}(x;y)$ has a finite limit as
$(x;y)\rto\AA$. Moreover, $\limA I_{\lg\mg}=0$ unless $\,{\lg\ll\mg}\,$ or
$\,{\lg=\mg}$, \cf. \(order), and
\vvm.3>
$$
\limA I_{\lg\)\lg}\,=\pron\prod_{s=0}^{\lg_m-1}\,{\eta^{-s}\>(\eta\1)\9\>
(\eta^s\al_{\lg,m}\1)\9\>(p\)\eta^{1-s-\lg_m}\al_{\lg,m}x_m/y_m)\9
\over(\eta^{-s-1})\9\>(\eta^{-s}x_m/y_m)\9\>(p)\9}\;.
$$
\endpro
\nt
Recall that
$\;{\al_{\lg,m}=\al\!\prjm\!\eta^{-2\lg_j}x_j/y_j}\ifMag\vp{\prod^.}\fi$.
\Prop{limI'}
For any $\lg,\mg\in\Zln$ the \hint/ $I'_{\lg\mg}(x;y)$ has a finite limit as
$(x;y)\rto\AA$. Moreover, $\limA I'_{\lg\mg}=0$ unless $\,{\lg\gg\mg}\,$ or
$\,{\lg=\mg}$, \cf. \(order), and
\vvm-.5>
$$
\limA I'_{\lg\)\lg}\,=\)\pron\prod_{s=0}^{\lg_m-1}\,{-\eta^{-s}\al_{\lg,m}\>
(\eta)\9\>(p\)\eta^{-s}\al_{\lg,m})\9\>
(\eta^{s-1+\lg_m}\al_{\lg,m}\1\)y_m/x_m)\9
\over(\eta^{s+1})\9\>(\eta^sy_m/x_m)\9\>(p)\9}\;.
$$
\endpro
\nt
The proof is similar to the proof of Proposition \[limI].
\Cr{limIb}
For any $\al,\eta$ the \hmap/s ${\Ib\)[\al;x;y;\eta\)]}$ and
${\Ib'[\al;x;y;\eta\)]}$ have finite limits as $(x;y)\rto\AA$. Moreover,
\ifMag
$$
\gather
\limA\Ib v^*_\lg\,=\,v_\lg\>\limA I_{\lg\)\lg}\;+\>
\sum_{\mg\gg\lg}v_\mg\>\limA I_{\lg\)\mg}
\\
\nn-5>
\Text{and}
\nn-3>
\limA\Ib' v^*_\lg\,=\,v_\lg\>\limA I'_{\lg\)\lg}\;+\>
\sum_{\mg\ll\lg}v_\mg\>\limA I'_{\lg\)\mg}\,.
\endgather
$$
\else
$$
\limA\Ib v^*_\lg\,=\,v_\lg\>\limA I_{\lg\)\lg}\;+\>
\sum_{\mg\gg\lg}v_\mg\>\limA I_{\lg\)\mg}\qquad\text{and\/}\qquad
\limA\Ib' v^*_\lg\,=\,v_\lg\>\limA I'_{\lg\)\lg}\;+\>
\sum_{\mg\ll\lg}v_\mg\>\limA I'_{\lg\)\mg}\,.
$$
\fi
\endpro

\vsk-.75>
\vsk0>
\Sect[P]{Proof of the \hRid/}
In this section we prove the \bid/ for the \hint/s. Its \eqv/ forms are given
by Theorems \[main] and \[main2]. We will prove the latter theorem.
\Pf of Theorem \[main2].
Let ${G_{\lg\mg}(\al;x;y;\eta)=\bigl(\Ib'[\al;x;y;\eta\)]\>v^*_\lg,
\Ib\)[\al;x;y;\eta\)]\>v^*_\mg\bigr)[x;y;\eta\)]}$.
\mmline
We have to prove that
$$
G_{\lg\mg}(\al;x;y;\eta)\,=\,\lpt v^*_\lg,v^*_\mg\rpt[\al;x;y;\eta\)]\,,\qqq
\lg,\mg\in\Zln\,.
\Tag{Flm}
$$
Since both sides of the above equality are analytic \fn/s of $\al$, we can
assume that $\al$ is generic. In particular, we will use the next statement.
\Lm{mu}
Let $\al$ be generic. Let $\mu_{\lg,k}$ be defined by \(Mm).
If $\,{\lg,\mg\in\Zln}$ are \st/ $\mu_{\lg,k}=\mu_{\mg,k}$ for any $\kn$,
then $\lg=\mg$.
\endpro
 From the definitions of operators $L_m$, $L'_m$ and the \Spair/ $S$
it is easy to see that
$$
S\)[x;y;\eta\)]\,=\,\bigl(L_m[\al;x;y;\eta\)]\bigr)^*S\)[x;y;\eta\)]\>
L'_m[\al;x;y;\eta\)]\,,\qqm \mn\,.
\ifMag\kern-2em\fi
$$
Therefore, for any $u,v\in V$ we have
$$
\bigl(K'_m[\al;x;y;\eta\)]\>u, K_m[\al;x;y;\eta\)]\>v\bigr)[x;y;\eta\)]\,=\,
(u,v)[x;y;\eta\)]\,,\qqm \mn\,.\ifMag\kern-2.5em\fi
\Tagg{KK}
$$
Hence, the \fn/ $G_{\lg\mg}(\al;x;y;\eta)$ satisfies a system of \deq/s
$$
T_k\)G_{\lg\mg}\,=\,\mu_{\lg,k}\1\>\mu_{\mg,k}\,G_{\lg\mg}\,,\qqq \kn\,,
\Tag{TF}
$$
see Theorem \[qKZ], Corollary \[qKZ'] and formulae \(Mm), \(KK).
\vvgood
On the other hand, Corollary \[limIb] shows that the \fn/ $G_{\lg\mg}(x;y)$
has a finite limit as $(x;y)\rto\AA$, and equations \(TF) imply that
$$
\smash{\limA G_{\lg\mg}\,=\,\mu_{\lg,k}\1\>\mu_{\mg,k}\,\limA G_{\lg\mg}\,,
\qqq\kn\,.}\strut
$$
Therefore, $\,{\limA G_{\lg\mg}=0}\,$ for ${\lg\ne\mg}$ by Lemma \[mu].
\par
Using once again \eq/s \(TF) we obtain that ${G_{\lg\mg}(x;y)\>=\)0}\,$
for $\lg\ne\mg$, and ${G_{\lg\lg}(x;y)\>=\)\limA G_{\lg\lg}}$. In particular,
the \fn/s ${G_{\lg\mg}(x;y)}$ are invariant \wrt/ the shift operators
$T_1^h\lc T_n^h$ for any nonzero $h$.
\par
Obviously, \rhs/ of \(Flm) enjoys the same properties:
$\lpt v^*_\lg,v^*_\mg\rpt=0\,$ for ${\lg\ne\mg}$, and
${\lpt v^*_\lg,v^*_\lg\rpt}$ is invariant \wrt/ $T_1^h\lc T_n^h$.
Hence, it remains to show that
$$
\limA G_{\lg\lg}\,=\,\lpt v^*_\lg,v^*_\lg\rpt\,,\qqq\lg\in\Zln\,,\vp{\Big|}
$$
which is a straightforward calculation using formulae \(vv), \(vv*),
Propositions \[limI], \[limI'] and Corollary \[limIb].
Theorem \[main2] is proved.
\epf

\Sect[R]{The restricted \hRid/}
Let us write down the \hRid/ using the bases of the \wtf/s in the \hgf/s:
$$
\gather
\smiZ M_\mg\>I'(W'_\lg,w'_\mg)\>I(W_\ng,w_\mg)\,=\,\dl_{\lg\ng}\>N_\lg\,,\qqq
\lg,\ng\in\Zln\,,
\Tag{MIIN}
\\
\nn4>
M_\mg\,=\,\prod_{i=1}^n\prod_{s=0}^{\mg_i-1}\,
{(1-\eta^{s+1})\>(x_i-\eta^sy_i)\over(1-\eta)\>\eta^s\)y_i}\;,
\\
\nn8>
N_\lg\,=\,\pron\prod_{s=0}^{\lg_m-1}\,{\eta^s\)\tht(\eta)\>
\tht(\eta^s\al_{\lg,m}\1)\>\tht(\eta^{1-s-\lg_m}\al_{\lg,m}\)x_m/y_m)
\over\tho\>\tht(\eta^{s+1})\>\tht(\eta^{-s}x_m/y_m)}
\endgather
$$
where $\,{\al_{\lg,m}=\al\!\prjm^{\vp.}\!\eta^{-2\lg_j}x_j/y_j}$,
\mmgood
\,\cf. Theorem \[main] and Lemma \[ResW]. Formula \(MIIN) holds for generic
values of the parameters $\eta$, $\xn$, $\yn$, \cf. \(gener), and all
the coefficients $N_\lg$, ${\lg\in\Zln}\!$, \,are clearly regular in this case.
Suppose now that ${x_1=\eta^ry_1}$ for some nonnegative integer $r$ less
than $\ell$, thus violating \(gener), and all other assumptions \(gener) hold.
Then the coefficients $N_\lg$ with $\lg_1\le r$ remain regular while
the coefficients $N_\lg$ with $\lg_1> r$ have a pole at $x_1=\eta^ry_1$.
This suggests that if the parameters $\eta$, $\xn$, $\yn$ are slightly
nongeneric, then some of the \hint/s can survive and still satisfy a certain
version of the \hRid/. We study such a possibility in this section.
\Par
Fix integers $\elln$ \st/ $1\le\ell_m\le\ell$, $\mn$. Set
$$
\Zlnb\,=\,\lb\)\lg\in\Zln\vert\lg_m\le\ell_m\,,\quad\mn\)\rb\,.
$$
Assume that for any $\roll$ and any $k,\mn$, we have
\ifMag
$$
\gather
\eta^{r+1}\nin\pZ,
\Tag{gene}
\\
\nn4>
\eta^{\pm r}x_k/x_m\nin\pZ,\qquad
\eta^{\pm r}y_k/y_m\nin\pZ,\qquad \eta^{-r}x_k/y_m\nin\pZ, \qquad k\ne m\,,
\\
\nn4>
\eta^{-s}x_m/y_m\nin\pZ,\qqq s=0\lc \ell_m-1\,.
\endgather
$$
\else
$$
\gather
\qquad\eta^{r+1}\nin\pZ,\qqq \eta^{\pm r}x_k/x_m\nin\pZ,\qquad
\eta^{\pm r}y_k/y_m\nin\pZ,\qquad \eta^{-r}x_k/y_m\nin\pZ, \qquad k\ne m\,,
\Tag{gene}
\\
\nn8>
\eta^{-s}x_m/y_m\nin\pZ,\qqq s=0\lc \ell_m-1\,.
\endgather
$$
\fi
Comparing with \(gener), here we impose weaker conditions for the ratios
$x_k/y_m$.
\Prop{welr}
Assume that nonzero parameters $\eta$, $\xn$, $\yn$ satisfy conditions \(gene).
Then for any $\lg\in\Zlnb$, $\mg\in\Zln$ the \hint/s $I(W_\lg,w_\mg)$ and
$I'(W'_\lg,w'_\mg)$ are well defined and are \hof/s of the parameters.
\endpro
\goodbreak
\Pf.
Assume that $|\eta|>1$ and $|x_m|<1$, $|y_m|>1$, $\mn$. Then
$$
I(W_\lg,w_\mg)\,=\;{1\over(2\pii)^\ell\>\ell\)!}\,
\int_{\){\TT^\ell}\]}W_\lg(t)\>w_\mg(t)\>\Phi(t)\;\dtt\;,
$$
\cf. \(hint0) and \(hpair). Observe that the integrand
${W_\lg(t)\>w_\mg(t)\>\Phi(t)}$ is a \sym/ \fn/ of $\tell$. Since the
integration contour $\TT^\ell$ is invariant \wrt/ \perm/s of the \var/s
$\tell$, we can drop the summation in the definition of the \fn/ $W_\lg$,
\cf. \(wW), multiplying the result of the integration by ${\ell\)!}$:
$$
\gather
I(W_\lg,w_\mg)\,=\;{1\over(2\pii)^\ell\!}\,
\int_{\){\TT^\ell}\]}P_\lg(t)\>w_\mg(t)\>\Phi(t)\;\dtt\;,
\\
\nn6>
P_\lg(t)\,=\,\pron\,\prod_{s=1}^{\lg_m}\,{\tht(\eta)\over\tht(\eta^s)}\,
\pron\,\prod_{a\in\Gm_{\Rph lm}}
\Bigl(\,{\tht(\eta^{\)2a-2}\al_m\1t_a/x_m)\over\tht(t_a/x_m)}
\prlm\,{\tht(t_a/y_l)\over\tht(t_a/x_l)}\,\Bigr)
\endgather
$$
where $\,{\Gm_m=\lb 1+\lg^{m-1}\,\lc\lg^m\rb}\,$ and
$\,{\al_m=\al\!\!\prlm\!\!\}x_l/y_l}$, $\ \mn$.
\goodbm
\vsk.2>
Under assumptions \(gene) singularities of the integrand
${P_\lg(t)\>w_\mg(t)\>\Phi(t)}$ are at most at the following hyperplanes:
\ifMag
$$
\gather
t_a=0\,,\qqq t_a=p^{-s}\eta\)t_b\,,\qquad 1\le a<b\le\ell\,,
\Tag{lisr}
\\
\nn4>
t_a=p^sx_j\,,\qquad t_a=p^sy_k\,,\qquad a\in\lb 1+\lg^{m-1}\,\lc\lg^m\rb\,,
\quad 1\le j\le m\le k\le\ell\,,
\endgather
$$
\else
$$
\alignat2
t_a &{}=0\,,\qqq  & t_a &{}=p^{-s}\eta\)t_b\,,\qquad 1\le a<b\le\ell\,,
\Tag{lisr}
\\
\nn6>
t_a=p^sx_j\,,\qquad t_a &{}=p^sy_k\,,\qquad &
a &{}\in\lb 1+\lg^{m-1}\,\lc\lg^m\rb\,,\quad 1\le j\le m\le k\le\ell\,,
\endalignat
$$
\fi
${s\in\Zp}$. The number of edges of configuration \(lisr) and dimensions of
the edges are always the same provided that the parameters are nonzero and
assumptions \(gene) hold. Therefore, the topology of the complement in $\Cl$
of the union of the hyperplanes \(lisr) does not change under the assumptions
of Proposition \[welr]. The rest of the proof is similar to the proof of
Theorem 5.7 in \Cite{TV2}.
\par
The proof for the \hint/ $I'(W'_\lg,w'_\mg)$ is similar to the proof for the
\hint/ $I(W_\lg,w_\mg)$ given above. Proposition \[welr] is proved.
\epf
\Rem
More detailed results for a similar problem concerning the multidimensional
\hint/s of \MB/ type are obtained in \Cite{MV}.
\enddemo
\Th{main3}
Assume that nonzero parameters $\eta$, $\xn$, $\yn$ satisfy conditions \(gene)
and $\,x_m=\eta^{\ell_m}y_m$ \,if $\,\ell_m<\ell$. Then
\vv-.3>
$$
\gather
\smiZb M_\mg\>I'(W'_\lg,w'_\mg)\>I(W_\ng,w_\mg)\,=\,\dl_{\lg\ng}\>N_\lg\,,
\qqq \lg,\ng\in\Zlnb\,,
\Tag{MIINr}
\\
\nn4>
\Text{where the coefficients $M_\mg,N_\lg$ are defined by \(MIIN).
Moreover, if $N_\mg\ne 0$ for all $\,\mg\in\Zlnb$, then}
\nn4>
\smiZb N_\mg\1\>I'(W'_\mg,w'_\lg)\>I(W_\mg,w_\ng)\,=\,\dl_{\lg\ng}\>M_\lg\1\,,
\qqq \lg,\ng\in\Zlnb.
\Tag{NIIM}
\endgather
$$
\vvv->
\endpro
\Pf.
Formula \(MIINr) follows from formula \(MIIN), because all the terms in \(MIIN)
are well defined under the assumptions of the theorem, see Proposition~\[welr],
and $M_\mg=0$ unless $\mg\in\Zlnb$.
\par
Writing down relation \(MIINr) in the matrix form: ${I'\)M\>I^t=N}$, for
matrices $I,I'\!,M,N$ with entries
$$
I_{\lg\mg}\)=\)I(W_\lg,w_\mg)\,,\qquad
I'_{\lg\mg}\)=\)I'(W'_\lg,w'_\mg)\,,\qquad
M_{\lg\mg}\)=\)\dl_{\lg\mg}M_\lg\,,\qquad N_{\lg\mg}\)=\)\dl_{\lg\mg}N_\lg\,,
$$
we immediately get that ${I^t\>N\1\)I'=M\1}\!$, which is the same as formula
\(NIIM).
\epf
\nt
We call relations \(MIINr) and \(NIIM) the \em{restricted \hgeom/ Riemann
identities}.
\Par
It is possible to introduce restricted versions of the \hgf/s, the \Spair/s and
the \hpair/s, and reformulate Theorem \[main3] similarly to Theorem \[main].
This will be done elsewhere.


\Sect[E]{Bilinear identities for \bhs/}
In this section we consider the \hRid/ in the one-dimen\-sional case.
That is, all over this section we assume that $\ell=1$.
\par
Let $\phifn(a_1\lc a_n;b_1\lc b_{n-1};z)$ be the \bhs/ \Cite{GR}:
\ifMag
$$
\gather
\phifn(\an;\bn;z)\,=\,\sum_{k=0}^\8\,{(a_1)_k\:\ldots(a_n)_k\:\over
(b_1)_k\:\ldots(b_{n-1})_k\:\)(p)_k\:\!}\;z^k\,,
\\
(u)_k\:=\>\prod_{s=0}^{k-1}(1-p^su)\,.
\endgather
$$
\else
$$
\phifn(\an;\bn;z)\,=\,\sum_{k=0}^\8\,{(a_1)_k\:\ldots(a_n)_k\:\over
(b_1)_k\:\ldots(b_{n-1})_k\:\)(p)_k\:\!}\;z^k\,,
\qqq (u)_k\:=\>\prod_{s=0}^{k-1}(1-p^su)\,.
$$
\fi
For any $\kn$ we define \fn/s $f_k$, $f'_k$, $F_k$, $F'_k$ by the formulae:
$$
\gather
f_k(t)\,=\,{t\over y_k\!}\,\pron{y_k-x_m\over t-x_m}\,
\prmk\>{t-y_m\over y_k-y_m}\;,\qquad\qquad f'_k(t)\,=\;{y_k\over t-y_k}\;,
\\
\ald
\nn4>
\ifMag
\Line{\ F_k(t)\,=\,(p)\9^2\,{\tht(\alt\1t/y_k)\over\tht(\alt\1)}\,\pron
{(p\)x_m/y_k)\9\over\tht(t/x_m)}\,\prmk\>{\tht(t/y_m)\over(p\)y_m/y_k)\9\!}\;,
\hfil \alt=\al\,\)\xy\,,\;}
\else
F_k(t)\,=\,(p)\9^2\,{\tht(\alt\1t/y_k)\over\tht(\alt\1)}\,\pron
{(p\)x_m/y_k)\9\over\tht(t/x_m)}\,\prmk\>{\tht(t/y_m)\over(p\)y_m/y_k)\9\!}\;,
\qqq \alt=\al\,\)\xy\,,
\fi
\\
\nn4>
F'_k(t)\,=\,(p)\9^2\,{\tht(\al t/y_k)\over\tht(\al)\>\tht(t/y_k)}\,
\pron(y_k/x_m)\9\,\prmk(y_k/y_m)\9\1\,.
\endgather
$$
The \fn/s ${\lb\)f_m\)\rb_{m=1}^n}$, ${\lb\)f'_m\)\rb_{m=1}^n}$,
\vv.2>
${\lb\)F_m\)\rb_{m=1}^n}$, ${\lb\)F'_m\)\rb_{m=1}^n}$ form bases in
the respective \hgf/s $\Fc$, $\Fc'\!$, $\Fq[\al]$, $\Fq'[\al]$,
and these bases are biorthonormal \wrt/ the \Spair/s:
$$
S(f'_l,f_m)\,=\,\dl_{lm}\,,\qqq\Sq(F'_l,F_m)\,=\,-\)\dl_{lm}\,.
$$
The \hint/s $I(F_l,f_m)$ and $I'(F'_l,f'_m)$ can be expressed via the \bhs/
$\phifn$. For instance,
\ifMag
$$
\alds
\align
& \aligned
I(F_1,f_1)\, &{}=\,
\phifn(x_1y_1\1\}\lc x_ny_1\1;y_2y_1\1\}\lc y_ny_1\1;\alt\1)\,,
\\
\nn6>
I(F_1,f_k)\, &{}=\,-\alt\1\>y_1\>y_k\1\)(y_1-p\)y_k)\1\prmk(y_m-y_k)\1
\pron(x_m-y_k)\,\x{}
\\
\nn2>
&\>{}\x\,\phifn(p\)x_1y_1\1\}\lc p\)x_ny_1\1;
p\)y_2y_1\1\}\lc p^2y_ky_1\1\}\lc p\)y_{n-1}y_1\1;\alt\1)\,,
\endaligned
\\
\nn12>
& \aligned
I'(F'_1,f'_1)\, &{}=\,
\phifn(y_1x_1\1\}\lc y_1x_n\1;y_1y_2\1\}\lc y_1y_n\1;\al\1)\,,
\\
\nn6>
I'(F'_1,f'_k)\,&{}=\,\alt\1\>y_k\>y_1\1\)(y_k-p\)y_1)\1
\prod_{m=2}^n(y_m-y_1)\1\pron(x_m-y_1)\,\x{}
\\
\nn2>
&\>{}\x\,\phifn(p\)y_1x_1\1\}\lc p\)y_1x_n\1;
p\)y_1y_2\1\}\lc p^2y_1y_k\1\}\lc p\)y_1y_{n-1}\1;\al\1)\,,\vrd2.7ex>
\endaligned
\endalign
$$
\else
$$
\alds
\gather
I(F_1,f_1)\,=\,\phifn(x_1y_1\1\}\lc x_ny_1\1;y_2y_1\1\}\lc y_ny_1\1;\alt\1)\,,
\\
\nn6>
\aligned
I(F_1,f_k)\, &{}=\,-\alt\1\>y_1\>y_k\1\)(y_1-p\)y_k)\1\prmk(y_m-y_k)\1
\pron(x_m-y_k)\,\x{}
\\
\nn2>
&\>{}\x\,\phifn(p\)x_1y_1\1\}\lc p\)x_ny_1\1;
p\)y_2y_1\1\}\lc p^2y_ky_1\1\}\lc p\)y_{n-1}y_1\1;\alt\1)\,,
\endaligned
\\
\nn12>
I'(F'_1,f'_1)\,=\,
\phifn(y_1x_1\1\}\lc y_1x_n\1;y_1y_2\1\}\lc y_1y_n\1;\al\1)\,,
\\
\nn6>
\aligned
I'(F'_1,f'_k)\,&{}=\,\alt\1\>y_k\>y_1\1\)(y_k-p\)y_1)\1
\prod_{m=2}^n(y_m-y_1)\1\pron(x_m-y_1)\,\x{}
\\
\nn2>
&\>{}\x\,\phifn(p\)y_1x_1\1\}\lc p\)y_1x_n\1;
p\)y_1y_2\1\}\lc p^2y_1y_k\1\}\lc p\)y_1y_{n-1}\1;\al\1)\,,\vrd2.7ex>
\endaligned
\endgather
$$
\fi
$k=2\lc n$. General formulae for $I(F_l,f_m)$, $I'(F'_l,f'_m)$ can be obtained
by a suitable change of indices.
\par
In the \onedim/ example in question, Theorem \[main] is \eqv/ to each of
the next formulae:
$$
\sun\,I(F_k,f_m)\>I'(F'_l,f'_m)\,=\,\dl_{kl}\,,\qqm\quad
\sun\,I(F_m,f_k)\>I'(F'_m,f'_l)\,=\,\dl_{kl}\,.\kern-1em
\Tag{II}
$$
These formulae deliver bilinear identities for the \bhs/. They read as follows:
\ifMag
$$
\alds
\gather
\Lline{\phifn(\an;\bn;z)\,\phifn(\ani;\bni;\zti)\,={}}
\\
\nn4>
\Rline{\aligned
\!\!\!{}=\,{}1+\smm^{n-1}\,{z^2A_0\>A_m\over(1-p\)b_m)\>(p-b_m)}\;
& \phifn(\anp;p\)b_1\lc p^2\]b_m\lc p\)b_{n-1};z)\,\x{}
\\
{}\x\,{} & \phifn(\anpi;p\)b_1\1\}\lc p^2\]b_m\1\}\lc p\)b_{n-1}\1;\zti)\,,
\endaligned}
\\
\nn8>
\Lline{\aligned
& \phifn(\an;\bn;z)\,\x{}
\\
\nn4>
&{}\,\x\,\phifn(\anbi;p^2\]b_1,p\)b_1b_2\1\}\lc p\)b_1b_{n-1}\1;\zti)\,={}
\\
\nn4>
&\kern1.7em{}=\,\phifn(\anp;p^2\]b_1,p\)b_2\lc p\)b_{n-1};z)\,\x{}
\\
\nn4>
&\kern3.4em{}\x\,\phifn(\anb;b_1,b_1b_2\1\}\lc b_1b_{n-1}\1;\zti)\,+{}
\endaligned}
\\
\nn5>
\Rline{\aligned
{}+\sum_{m=2}^{n-1}\,{z\>A_m\>(1-p\)b_1)\over(1-p\)b_m)\>(b_m-p\)b_1)}\;
\phifn(\anp;p\)b_1\lc p^2\]b_m\lc p\)b_{n-1};z)\x{} &
\\
{}\x\,\phifn(\anbi;
p\)b_1,p\)b_1b_2\1\}\lc p^2\]b_1b_m\1\}\lc p\)b_1b_{n-1}\1;\zti) & \,,
\endaligned}
\\
\nn10>
\Lline{\phifn(\an;\bn;z)\,\phifn(\ani;\bni;\zti)\,={}}
\\
\nn4>
\Rline{\aligned
{}=\,{}1+\smm^{n-1}\, & {z^2A_0\>A_m\over(1-p\)b_m)\>(p-b_m)}\ \x{}
\\
\nn2>
&{}\x\,
\phifn(\anpbm;p\)b_m\1b_1\lc p^2\]b_m\1\}\lc p\)b_m\1b_{n-1};z)\,\x{}
\\
\nn4>
&\kern1.6em{}\x\,
\phifn(\anbmi;p\)b_mb_1\1\}\lc p^2\]b_m\lc p\)b_mb_{n-1}\1;\zti)\,,
\endaligned}
\\
\nn8>
\gathered
\Lline{\phifn(\an;\bn;z)\,
\phifn(\anpi;p^2\]b_1\1\},p\)b_2\1\}\lc p\)b_{n-1}\1;\zti)\,={}}
\\
\Lline{\aligned
\hp{\phifn}\Llap{{}=\,{}}
\phifn(\anpb;p^2\]b_1\1,p\)b_1\1b_2\lc p\)b_1\1b_{n-1};z)\,\x{}&
\\
\nn3>
{}\x\,\phifn(\anb;b_1,b_1b_2\1\}\lc b_1b_{n-1}\1;\zti)\,+{}&
\endaligned}
\endgathered
\\
\nn5>
\Rline{\aligned
{}+\sum_{m=2}^{n-1}\, & {z\>A_m\>(p-b_1)\over(p-b_m)\>(b_1-p\)b_m)}\ \x{}
\\
\nn2>
&{}\x\,
\phifn(\anpbm;p\)b_m\1b_1\lc p^2\]b_m\1\}\lc p\)b_m\1b_{n-1};z)\x{}
\\
\nn4>
&\kern1.6em{}\x\,
\phifn(\anbmi;p\)b_m,p^2\]b_mb_1\1\},p\)b_mb_2\1\}\lc p\)b_mb_{n-1}\1;\zti)
\endaligned}
\\
\nn2>
\Text{where}
\zti\)=\)z\;{\prod a_m\over\prodp\]b_m}\;,\qqq
A_0\)=\,{\prod\>(1-a_m)\over\prodp(1-b_m)}\;,\qquad
A_k\)=\,{\prod\>(a_m-b_k)\over(1-b_k)\>
\mathop{\prodp}\limits_{\!m\ne k}\}(b_m-b_k)}\;,
\endgather
$$
$\kn-1$. \;Here and below $\prod$ stands for $\!\topsmash{\pron}\!$ \;and
$\,\prodp$ stands for $\topsmash{\prmn}$.
\else
$$
\alds
\gather
\aligned
\phifn( & \an;\bn;z)\,\phifn(\ani;\bni;\zti)\,={}
\\
\nn4>
& \aligned
\!\!\!{}=\,{}1+\smm^{n-1}\,{z^2A_0\>A_m\over(1-p\)b_m)\>(p-b_m)}\;
& \phifn(\anp;p\)b_1\lc p^2\]b_m\lc p\)b_{n-1};z)\,\x{}
\\
{}\x\,{} & \phifn(\anpi;p\)b_1\1\}\lc p^2\]b_m\1\}\lc p\)b_{n-1}\1;\zti)\,,
\endaligned
\endaligned
\\
\nn8>
\gathered
\Lline{\phifn(\an;\bn;z)\,
\phifn(\anbi;p^2\]b_1,p\)b_1b_2\1\}\lc p\)b_1b_{n-1}\1;\zti)\,={}}
\\
\nn4>
\Rline{{}=\,\phifn(\anp;p^2\]b_1,p\)b_2\lc p\)b_{n-1};z)\,
\phifn(\anb;b_1,b_1b_2\1\}\lc b_1b_{n-1}\1;\zti)\,+{}}
\\
\nn5>
\Rline{\aligned
{}+\sum_{m=2}^{n-1}\,{z\>A_m\>(1-p\)b_1)\over(1-p\)b_m)\>(b_m-p\)b_1)}\;
\phifn(\anp;{}& p\)b_1\lc p^2\]b_m\lc p\)b_{n-1};z)\x{}
\\
{}\x\,\phifn(\anbi; {}&
p\)b_1,p\)b_1b_2\1\}\lc p^2\]b_1b_m\1\}\lc p\)b_1b_{n-1}\1;\zti)\,,
\endaligned}
\endgathered
\\
\nn10>
\gathered
\Lline{\phifn(\an;\bn;z)\,\phifn(\ani;\bni;\zti)\,={}}
\\
\nn4>
\Rline{{}=\,{}1+\smm^{n-1}\,{z^2A_0\>A_m\over(1-p\)b_m)\>(p-b_m)}\;
\phifn(\anpbm;p\)b_m\1b_1\lc p^2\]b_m\1\}\lc p\)b_m\1b_{n-1};z)\,\x{}\hp{\,,}}
\\
\Rline{{}\x\,
\phifn(\anbmi;p\)b_mb_1\1\}\lc p^2\]b_m\lc p\)b_mb_{n-1}\1;\zti)\,,}
\endgathered
\\
\nn8>
\gathered
\Lline{\phifn(\an;\bn;z)\,
\phifn(\anpi;p^2\]b_1\1\},p\)b_2\1\}\lc p\)b_{n-1}\1;\zti)\,={}}
\\
\nn5>
\Lline{\aligned
\hp{\phifn}\Llap{{}=\,{}}
\phifn(\anpb;p^2\]b_1\1,p\)b_1\1b_2\lc p\)b_1\1b_{n-1};z)\,\x{}&
\\
\nn3>
{}\x\,\phifn(\anb;b_1,b_1b_2\1\}\lc b_1b_{n-1}\1;\zti)\,+{}&
\endaligned}
\\
\nn5>
\Rline{\aligned
{}+\sum_{m=2}^{n-1}\,{z\>A_m\>(p-b_1)\over(p-b_m)\>(b_1-p\)b_m)}\;
\phifn(\anpbm;{}& p\)b_m\1b_1\lc p^2\]b_m\1\}\lc p\)b_m\1b_{n-1};z)\x{}
\\
{}\x\,\phifn(\anbmi; {}&
p\)b_m,p^2\]b_mb_1\1\},p\)b_mb_2\1\}\lc p\)b_mb_{n-1}\1;\zti)\,,
\endaligned}
\endgathered
\\
\Text{where}
\nn4>
\zti\)=\)z\;{\prod a_m\over\prodp\]b_m}\;,\qqq
A_0\)=\,{\prod\>(1-a_m)\over\prodp(1-b_m)}\;,\qquad
A_k\)=\,{\prod\>(a_m-b_k)\over(1-b_k)\>
\mathop{\prodp}\limits_{\!m\ne k}\}(b_m-b_k)}\;,\qquad\kn-1\,.
\endgather
$$
Here and below $\prod$ stands for $\!\topsmash{\pron}\!$ \;and $\,\prodp$
stands for $\topsmash{\prmn}$.
\fi
\Par
\Prop{detIF}
$\ \dsize\det\bigl[I(F_l,f_m)\bigr]_{l,m=1}^n\,=\,
{(\al\1)\9\over(\al\1\prod y_m/x_m)\9\!}$.
\endpro
\Pf.
Let $\eg(m)=(0\lc\%1_{\sss\^{$m$-th}}\lc 0)$, $\mn$. Let
$w_m(t)=w_{\eg(m)}(t;x;y;\eta)$ and $W_m(t)=W_{\eg(m)}(t;\al;x;y;\eta)$ be
the \wtf/s. We have that
$$
\alds
\alignat3
w_m(y_l)\) &{}=\)0\,,& \qqq W_m(y_l)\) &{}=\)0\,, & \qqq & 1\le l<m\le n\,,
\\
\nn\ifMag-4\else-2\fi>
\Text{and}
\nn\ifMag-4\else-2\fi>
f_m(y_l)\) &{}=\)0\,,& \qqq F_m(y_l)\) &{}=\)0\,, &\qqq & 1\le l\ne m\le n\,.
\endalignat
$$
Therefore,
$$
\gather
\\
\nn-28>
\det\bigl[I(F_l,f_m)\bigr]_{l,m=1}^n\,=\,\pron\,
{f_m(y_m)\>F_m(y_m)\over w_m(y_m)\>W_m(y_m)}\;
\det\bigl[I(W_l,w_m)\bigr]_{l,m=1}^n
\endgather
$$
and the claim follows from Proposition \[detIW] for $\ell=1$.
\epf
Calculating a matrix inverse to the matrix
$\bigl[I'(F'_l,f'_m)\bigr]_{l,m=1}^n$ in two different ways using
either formulae \(II) or Proposition \[detIF], we obtain that
$$
I(F_l,f_m)\,=\;{(-1)^{l+m}\>(\al\1)\9\over(\al\1\prod y_m/x_m)\9}\;
\det\bigl[I'(F'_j,f'_k)\bigr]_{\vtop{\bls.63\bls
\mbox{\ssize j,k=1}\mbox{\ssize j\ne l}\mbox{\ssize k\ne m}}}^n\,,
\qqq l,\mn\,.
$$
For $n=2$ these relations are \eqv/ to
$$
\phif21(a_1,a_2;b\);z)\>(z)\9\,=\,
\phif21(b\)a_1\1\},b\)a_2\1;b\);z\)a_1a_2\)b\1)\>(z\)a_1a_2\)b\1)\9
$$
and for $n=3$ they give
\ifMag
$$
\gather
\\
\nn-24>
\Lline{\phif32(a_1,a_2,a_3;b_1,b_2;z)\,{(z)\9\over(\zti)\9}\;={}}
\\
\nn6>
\aligned
{}=\,{}\phif32(b_1a_1\1\},b_1a_2\1\},b_1a_3\1;b_1,b_1b_2\1;\zti)\,
\phif32(b_2a_1\1\},b_2a_2\1\},b_2a_3\1;b_2,b_2b_1\1;\zti)\,-{} &
\\
\nn6>
{}-\,z^2\,{(a_1-b_1)\>(a_2-b_1)\>(a_3-b_1)\>(a_1-b_2)\>(a_2-b_2)\>(a_3-b_2)
\over(1-b_1)\>(1-b_2)\>(b_1-b_2)^2\>(b_1-p\)b_2)\>(p\)b_1-b_2)}\ \x{}
\kern.7em &
\\
\nn5>
{}\x\,
\phif32(p\)b_1a_1\1\},p\)b_1a_2\1\},p\)b_1a_3\1;p\)b_1,p^2b_1b_2\1;\zti)\,\x{}
\kern1.4em &
\endaligned
\\
\nn5>
\Rline{{}\x\,
\phif32(p\)b_2a_1\1\},p\)b_2a_2\1\},p\)b_2a_3\1;p\)b_2,p^2b_2b_1\1;\zti)\,,}
\endgather
$$
$$
\gather
\Lline{
\phif32(p\)a_1,p\)a_2,p\)a_3;p^2b_1,p\)b_2;z)\,{(z)\9\over(\zti)\9}\;={}}
\\
\nn6>
\Lline{{}\]=\,{}
\phif32(p\)b_1a_1\1\},p\)b_1a_2\1\},p\)b_1a_3\1;p^2b_1,p\)b_1b_2\1;\zti)\,
\phif32(b_2a_1\1\},b_2a_2\1\},b_2a_3\1;b_2,b_2b_1\1;\zti)\,+{}}
\\
\nn6>
\Rline{\aligned
& {}+\,z\,{(1-p\)b_1)\>(a_1-b_2)\>(a_2-b_2)\>(a_3-b_2)\over
(1-b_2)\>(1-p\)b_2)\>(b_1-b_2)\>(p\)b_1-b_2)}\ \x{}
\\
\nn5>
& \kern4em {}\x\,
\phif32(p\)b_1a_1\1\},p\)b_1a_2\1\},p\)b_1a_3\1;p\)b_1,p^2b_1b_2\1;\zti)\,\x{}
\\
\nn5>
& \kern8em {}\x\,
\phif32(p\)b_2a_1\1\},p\)b_2a_2\1\},p\)b_2a_3\1;p^2b_2,p\)b_2b_1\1\};\zti)
\endaligned}
\endgather
$$
\else
$$
\align
&\) \phif32(a_1,a_2,a_3;b_1,b_2;z)\,{(z)\9\over(\zti)\9}\;={}
\\
\nn6>
& \aligned
{}=\,{} & \phif32(b_1a_1\1\},b_1a_2\1\},b_1a_3\1;b_1,b_1b_2\1;\zti)\,
\phif32(b_2a_1\1\},b_2a_2\1\},b_2a_3\1;b_2,b_2b_1\1;\zti)\,-{}
\\
\nn6>
& {}\>-\,z^2\,{(a_1-b_1)\>(a_2-b_1)\>(a_3-b_1)\>(a_1-b_2)\>(a_2-b_2)\>(a_3-b_2)
\over(1-b_1)\>(1-b_2)\>(b_1-b_2)^2\>(b_1-p\)b_2)\>(p\)b_1-b_2)}\ \x{}
\\
\nn4>
& {}\>\x\,
\phif32(p\)b_1a_1\1\},p\)b_1a_2\1\},p\)b_1a_3\1;p\)b_1,p^2b_1b_2\1;\zti)\,
\phif32(p\)b_2a_1\1\},p\)b_2a_2\1\},p\)b_2a_3\1;p\)b_2,p^2b_2b_1\1;\zti)\,,
\endaligned
\\
\nn8>
&\) \phif32(p\)a_1,p\)a_2,p\)a_3;p^2b_1,p\)b_2;z)\,{(z)\9\over(\zti)\9}\;={}
\\
\nn6>
& \aligned
{}=\,{}&
\phif32(p\)b_1a_1\1\},p\)b_1a_2\1\},p\)b_1a_3\1;p^2b_1,p\)b_1b_2\1;\zti)\,
\phif32(b_2a_1\1\},b_2a_2\1\},b_2a_3\1;b_2,b_2b_1\1;\zti)\,+{}
\\
\nn6>
& {}\>+\,z\,{(1-p\)b_1)\>(a_1-b_2)\>(a_2-b_2)\>(a_3-b_2)\over
(1-b_2)\>(1-p\)b_2)\>(b_1-b_2)\>(p\)b_1-b_2)}\ \x{}
\\
\nn4>
& {}\>\x\,
\phif32(p\)b_1a_1\1\},p\)b_1a_2\1\},p\)b_1a_3\1;p\)b_1,p^2b_1b_2\1;\zti)\,
\phif32(p\)b_2a_1\1\},p\)b_2a_2\1\},p\)b_2a_3\1;p^2b_2,p\)b_2b_1\1;\zti)
\endaligned
\endalign
$$
\fi
where $\zti=z\)a_1a_2a_3\)b_1\1b_2\1\!$. \,Notice that
$(z)\9/(\zti)\9=\phif10(z/\zti;\zti)=\phif10(\zti/z;z)\vpb{-1}\!$.
\Par
The \rat/ version of the \hRid/ in the \onedim/ case gives similar formulae for
the generalized \hgeom/ \fn/ $\Fifn$. They can be obtained from the formulae
for the \bhs/ after the standard substitution ${a_m\}=p^{\)\al_m}}$,
${b_m\}=p^{\>\bt_m}}$ in the limit ${p\to 1}$ which degenerates
$\phifn(\an;\bn;z)\,$ to $\Fifn(\al_1\lc\al_n;\bt_1\lc\bt_{n-1};z)$.

\Appendix

\Sect[S]{Nondegeneracy of the elliptic Shapovalov pairing}
Let $\topsmash{A=\al\)\eta^{1-\ell}\!\pron\!\]x_m}\vp{\prod^{}}$.
\vv.1>
Let $\Ec[A]$ be the space of \hof/s on $\Cx\!$ \st/ $f({p\)u})=A(-u)^{-n}f(u)$.
It is easy to see that $\dim\)\Ec[A]=n$, say by Fourier series.
\Par
Let $\om=\exp(2\pii/n)$. Fix complex numbers $\xi$ and $\zt$ \st/
${\xi^n\]=p}$ and ${\zt^n\]=-A\1\!}$. Set
\vvm-.4>
$$
\smash{\thi_l(u)\,=\,u^{l-1}\pron\tht(-\zt\)\xi^{l-1}\om^mu)\,,\qqq\lcn\,.}
\vp{\prod_{}}
$$
\Lm{thi}
The \fn/s $\thi_1\lc\thi_n$ form a basis in the space $\Ec[A]$.
\endpro
\Pf.
Clearly, $\thi_l\in\Ec[A]$ for any $\lcn$. Moreover,
${\thi_l(\om\)u)=\om^{l-1}\thi_l(u)}$, that is the \fn/s $\thi_1\lc\thi_n$
are eigen\fn/s of the translation operator with distinct \eva/s.
Hence, they are linearly independent.
\epf
For any $\lg\in\Zln$ let $G_\lg(t;\al;x;y;\eta)$ be the following \fn/:
\ifMag
$$
\align
\qquad & \,G_\lg(t;\al;x;y;\eta)\,={}
\Tagg{Glg}
\\
\nn3>
& {}=\;{1\over\lg_1!\ldots\lg_n!}\,
\pron\,\pral\,{1\over\tht(t_a/x_m)}\;
\prab\ {\tht(t_a/t_b)\over\tht(\eta\)t_a/t_b)}\;
\susi\,\pron\,\prod_{a\in\Gm_{\Rph lm}}\thi_m(t_{\si_a})\,.
\endalign
$$
\else
$$
G_\lg(t;\al;x;y;\eta)\,=\;{1\over\lg_1!\ldots\lg_n!}\,
\pron\,\pral\,{1\over\tht(t_a/x_m)}\;
\prab\ {\tht(t_a/t_b)\over\tht(\eta\)t_a/t_b)}\;
\susi\,\pron\,\prod_{a\in\Gm_{\Rph lm}}\thi_m(t_{\si_a})\,.
\Tag{Glg}
$$
\fi
Here $\Gm_m=\lb 1+\lg^{m-1}\,\lc\lg^m\rb$, $\mn$.
\Lm{}
The \fn/s $G_\lg(t;\al;x;\eta)$, ${\lg\in\Zln}$, form a basis in the \ehgf/
${\Fq[\al;x;\eta\);\ell\>]}$.
\endpro
\goodbm
\Pf.
The \ehgf/ ${\Fq[\al;x;\eta\);\ell\>]}$ is naturally isomorphic to the
\^{$\ell$-}th \sym/ power of the space $\Ec[A]$ --- the space of \sym/ \fn/s
in $\tell$ which considered as \fn/s of one \var/ $t_a$ belong to $\Ec[A]$ for
any $\aell$. The \iso/ reads as follows:
\ifMag
\vv.5>
$$
\Line{f(\tell)\,\map\,f(\tell)\,\pron\,\pral\,\tht(t_a/x_m)
\prab\,{\tht(\eta\)t_a/t_b)\over\tht(t_a/t_b)}\;,\hfil
f\in\Fq[\al;x;\eta\)]\,.}
$$
\vvv.5>
\else
$$
f(\tell)\,\map\,f(\tell)\,\pron\,\pral\,\tht(t_a/x_m)
\prab\,{\tht(\eta\)t_a/t_b)\over\tht(t_a/t_b)}\;,\qqq
f\in\Fq[\al;x;\eta\)]\,.
$$
\fi
Now the proposition follows from Lemma \[thi].
\epf
Let $W_\lg$, $\lg\in\Zln$, be the \ewf/s. Define a matrix $Q(\al;x;y;\eta)$ by
the rule:
\vvm->
$$
\gather
W_\lg(t;\al;x;y;\eta)\,=\smiZ Q_{\lg\mg}(\al;x;y;\eta)\>G_\mg(t;\al;x;\eta)\,,
\qqq \lg\in\Zln\,.
\\
\nn-4>
\Text{Set}
\nn-4>
d(n,m,\ell,s)\,=
\sum_{\tsize{\vp( i,j\ge 0\atop\tsize{\vp( i+j<\ell\atop i-j=s\vp(}}}
{m-1+i\choose m-1}\,{n-m-1+j\choose n-m-1}
\Tag{dnm}
\endgather
$$
\ifMag\vsk->\vsk0>\fi
\Prop{detQ}
\back\Cite{TV1}
$$
\gather
\\
\nn-24>
\ifMag\nn-4>\fi
\aligned
\det Q(\al;x;y;\eta)\,& {}=\;\Xi\>\prsl\,\prmn\>\tht\bigl(
{\tsize\eta^{s+\ell-1}\al\1\!\!\!\prllm\!\!y_l/x_l}\bigr)^{d(n,m,\ell,s)}\ \x
\\
\nn6>
& \>{}\x\,\pron y_m^{(m-n){\tsize{n+\ell-1\choose n}}}\;
\pros\,\plmn\>\tht(\eta^s y_l/x_m)\vpb{\tsize{n+\ell-s-2\choose n-1}}
\endaligned
\\
\nn4>
\Text{where}
\nn-6>
\Xi\;=\,\Bigl[\>\Rph{(p)}{(p)\9}\vpb{1-n^2}\>\prmn\>
\Bigl(\>{\tht(\om^m)\over\om^m-1}\)\Bigr)^{\!\!\raise 3pt\mbox{\ssize n-m}}\>
\Bigr]^{\tsize{n+\ell-1\choose n}}.
\Tag{Xi}
\endgather
$$
\endpro
Let ${\Sq=\Sq[\al;x;y;\eta\)]}$ be the elliptic \Spair/.
\vv.15>
Let ${G_\lg=G_\lg(t;\al;x;\eta)}$ and
\nline
$G'_\lg=G_\lg(t;\al\1;y;\eta\1)$.
\Prop{detS}
\ifMag
$$
\alignat2
\\
\nn-22>
\det &\Rlap{\bigl[\Sq(G'_\lg,G_\mg)\bigr]_{\lg,\mg\in\Zln}\,={}} &&
\\
\nn6>
&& \}{}=\,\Xi^{-2}\>
(-1)\vpb{n(n-1)/2\,\cdot\!\tsize{n+\ell-1\choose n}}\>
\eta\vpb{\,n(3-n)/2\,\cdot\!\tsize{n+\ell-1\choose n+1}}
\tpron x_m^{(n-1)\tsize{n+\ell-1\choose n}}\,\x{}\! &
\\
\nn3>
&& {}\>\x\,\pros\,\biggl[\,{\tht(\eta)\vpb{n}\tht(\eta^s\al\1)\>
\tht\bigl(\eta^{s+2-2\ell}\al\prod x_m/y_m\bigr)\over
\tho\vpb{n}\tht(\eta^{s+1})\vpb{n}\prod\!\prod\tht(\eta^{-s}x_l/y_m)}
\,\biggr]^{\tsize{n+\ell-s-2\choose n-1}} &.
\endalignat
$$
\else
$$
\align
\\
\nn-18>
\det\bigl[\Sq(G'_\lg,G_\mg)\bigr]_{\lg,\mg\in\Zln}\,=\,\Xi^{-2}\>
(-1)\vpb{n(n-1)/2\,\cdot\!\tsize{n+\ell-1\choose n}}\>
\eta\vpb{\,n(3-n)/2\,\cdot\!\tsize{n+\ell-1\choose n+1}}
\tpron x_m^{(n-1)\tsize{n+\ell-1\choose n}}\,\x{}\! &
\\
\nn3>
{}\>\x\,\pros\,\biggl[\,{\tht(\eta)\vpb{n}\)\tht(\eta^s\al\1)\>
\tht\bigl(\eta^{s+2-2\ell}\al\prod x_m/y_m\bigr)\over
\tho\vpb{n}\)\tht(\eta^{s+1})\vpb{n}\)\prod\!\prod\tht(\eta^{-s}x_l/y_m)}
\,\biggr]^{\tsize{n+\ell-s-2\choose n-1}} &.
\endalign
$$
\fi
Here $\Xi$ is given by \(Xi), $\,\prod$ stands for $\!\pron$ and
$\,\prod\!\prod$ stands for $\prod_{l=1}^n\pron\!$.
\endpro
\goodbm
\Pf.
By Lemma \[ResW] we have that
$$
\align
\det\bigl[\Sq(G'_\lg,G_\mg)\bigr]_{\lg,\mg\in\Zln}\, &{}=\,
\bigl(\det Q(\al;x;y;\eta)\)\det Q(\al\1;y;x;\eta\1)\bigr)\1\,\x{}
\\
\nn4>
&\>{}\x\,\prod_{\lg\in\Zln}\pron\prod_{s=0}^{\lg_m-1}\,{\eta^s\)\tht(\eta)\>
\tht(\eta^s\al_{\lg,m}\1)\>\tht(\eta^{1-s-\lg_m}\al_{\lg,m}\)x_m/y_m)
\over\tho\>\tht(\eta^{s+1})\>\tht(\eta^{-s}x_m/y_m)}\;.
\endalign
$$
To get the final answer we use Proposition \[detQ] and simplify the triple
product changing the order of the products and applying Lemma \[combi] several
times.
\epf
\Lm{combi}
The following identity holds:
$$
\sum_{a=0}^l\>{j+a\choose j}{j+k+a\choose k}{l+m-a\choose m}\,=\,
{j+k\choose k}{j+k+l+m+1\choose j+k+m+1}\,.
$$
\endpro
\nt
The statement can be proved by induction \wrt/ $l$ and $m$.

\Sect[D]{Three determinant formulae}
For any $\lg\in\Zln$ let $g_\lg(t;x;y;\eta)$ be the following \fn/:
\ifMag
$$
\align
\quad & \;g_\lg(t;x;y;\eta)\,={}
\Tagg{glg}
\\
\nn3>
& {}=\;{1\over\lg_1!\ldots\lg_n!}\,
\pron\,\pral\,{1\over t_a-x_m}\;\prab\ {t_a-t_b\over\eta\)t_a-t_b}\;
\susi\,\pron\,\prod_{a\in\Gm_{\Rph lm}}t_{\si_a}^m\,.
\endalign
$$
\else
$$
g_\lg(t;x;y;\eta)\,=\;{1\over\lg_1!\ldots\lg_n!}\,
\pron\,\pral\,{1\over t_a-x_m}\;\prab\ {t_a-t_b\over\eta\)t_a-t_b}\;
\susi\,\pron\,\prod_{a\in\Gm_{\Rph lm}}t_{\si_a}^m\,.
\Tag{glg}
$$
\fi
Here $\Gm_m=\lb 1+\lg^{m-1}\,\lc\lg^m\rb$, $\mn$. The \fn/s $g_\lg(t;x;\eta)$,
${\lg\in\Zln}$, form a basis in the \thgf/ ${\Fc[x;\eta\);\ell\>]}$.
\Par
Let $w_\lg$, ${\lg\in\Zln}$, be the \twf/s. Define a matrix $X(x;y;\eta)$ by
the rule:
$$
\smash{w_\lg(t;x;y;\eta)\,=\smiZ X_{\lg\mg}(x;y;\eta)\>g_\mg(t;x;\eta)\,,
\qqq\lg\in\Zln\,.}\vp\sum\ifMag\else\vrd2ex>\fi
$$
\Prop{detX}
\back\Cite{T},\;\Cite{TV1}
$\ \;\dsize\det X(x;y;\eta)\,=\,\pros\,\plmn\>(\eta^s y_l-x_m)
\vpb{\tsize{n+\ell-s-2\choose n-1}}$.
\endpro
\nt
Let $W_\lg$, ${\lg\in\Zln}$, be the \ewf/s and let ${I=I\)[\al;x;y;\eta\)]}$
be the \hpair/.
\Prop{detIW}
\back\Cite{TV1}
\ifMag
$$
\gather
\\
\nn-24>
\Lline{\det\bigl[I(W_\lg,w_\mg)\bigr]_{\lg,\mg\in\Zln}\,=\,
\eta\vpb{-n\tsize{n+\ell-1\choose n+1}}\prsl\,\prmn\>\tht\bigl(
{\tsize\eta^{s+\ell-1}\al\1\!\!\!\prllm\!\!y_l/x_l}\bigr)^{d(n,m,\ell,s)}
\;\x{}}
\\
\nn4>
\Rline{{}\x\,\pros\,\biggl[\,{(\eta\1)\9^n\>(\eta^s\al\1)\9\>
(p\)\eta^{s+2-2\ell}\al\prod x_m/y_m)\9\over
(\eta^{-s-1})\9^n\>(p)\9^{2n-1}\)\prod(\eta^{-s}x_m/y_m)\9}\!
\plmn{(\eta^sy_l/x_m)\9\over(\eta^{-s}x_l/y_m)\9\!}\;\biggr]
^{\tsize{n+\ell-s-2\choose n-1}}.}
\endgather
$$
\else
$$
\align
\\
\nn-24>
\det & \bigl[I(W_\lg,w_\mg)\bigr]_{\lg,\mg\in\Zln}\,=\,
\eta\vpb{-n\tsize{n+\ell-1\choose n+1}}\prsl\,\prmn\>\tht\bigl(
{\tsize\eta^{s+\ell-1}\al\1\!\!\!\prllm\!\!y_l/x_l}\bigr)^{d(n,m,\ell,s)}
\;\x{}
\\
\nn4>
& {}\!\x\,\pros\,\biggl[\,{(\eta\1)\9^n\>(\eta^s\al\1)\9\>
(p\)\eta^{s+2-2\ell}\al\prod x_m/y_m)\9\over
(\eta^{-s-1})\9^n\>(p)\9^{2n-1}\)\prod(\eta^{-s}x_m/y_m)\9}\!
\plmn{(\eta^sy_l/x_m)\9\over(\eta^{-s}x_l/y_m)\9\!}\;\biggr]
^{\tsize{n+\ell-s-2\choose n-1}}.
\endalign
$$
\fi
Here $\prod$ stands for $\!\pron\!$ \;and the exponents $d(n,m,\ell,s)$
are given by \(dnm).
\endpro
\goodbm
\nt
Let \fn/s $G_\lg$, ${\lg\in\Zln}$, be given by \(Glg).
\Cr{detIG}
\back\Cite{TV1}
\ifMag
$$
\gather
\\
\nn-24>
\Lline{\det\bigl[I(G_\lg,g_\mg)\bigr]_{\lg,\mg\in\Zln}\,=\,\Xi\1\>
\eta\vpb{-n(n+1)/2\,\cdot\!\tsize{n+\ell-1\choose n+1}}\,\x{}}
\\
\nn4>
\Rline{{}\x\,\pros\,\biggl[\,{(\eta\1)\9^n\>(\eta^s\al\1)\9\>
(p\)\eta^{s+2-2\ell}\al\prod x_m/y_m)\9\over(\eta^{-s-1})\9^n\>
(p)\9^{2n-1+n(n-1)/2}\)\prod\!\prod(\eta^{-s}x_l/y_m)\9}\,\biggr]
^{\tsize{n+\ell-s-2\choose n-1}}.}
\endgather
$$
\else
$$
\align
\\
\nn-24>
\det\bigl[I(G_\lg,g_\mg)\bigr]_{\lg,\mg\in\Zln}\, &{}=\,\Xi\1\>
\eta\vpb{-n(n+1)/2\,\cdot\!\tsize{n+\ell-1\choose n+1}}\,\x{}
\\
\nn4>
& {}\>\x\,\pros\,\biggl[\,{(\eta\1)\9^n\>(\eta^s\al\1)\9\>
(p\)\eta^{s+2-2\ell}\al\prod x_m/y_m)\9\over(\eta^{-s-1})\9^n\>
(p)\9^{2n-1+n(n-1)/2}\)\prod\!\prod(\eta^{-s}x_l/y_m)\9}\,\biggr]
^{\tsize{n+\ell-s-2\choose n-1}}.
\endalign
$$
\fi
Here $\,\Xi\,$ is given by \(Xi), $\prod$ stands for $\!\pron\!$ \;and
$\,\prod\!\prod$ stands for $\prod_{l=1}^n\pron\!$.
\endpro
\Pf.
The statement follows from Propositions \[detQ], \[detX] and \[detIW].
\epf

\Sect[J]{The Jackson integrals via the \hint/s}
Consider the \hint/ ${\IH[x;y;\eta\)](f\)\Pht)}$, \cf. \(hint1), for a \fn/
\mmline
$f(\tell)$ of the form
$$
f(\tell)\,=\,P(\tell)\>\Tht(\tell)\prab(t_a/t_b)\9
\Tag{fPT}
$$
where $P(\tell)$ is a \sym/ \pol/ of degree at most $M$ in each of the \var/s
$\tell$ and $\Tht(\tell)$ is a \sym/ \hof/ on $\Cxl\!$ \st/
$$
\Tht(\tpell)\,=\,A\>(-t_a)^{-n}\>\Tht(\tell)
\Tag{TA}
$$
for some constant $A$. The \hint/s which appear in the definition of
the \hpair/s, see \(hpair), fit this case for $M=n-1$ and $A$ determined by
$\al,\eta$, $\xn$, $\yn$.
\goodbreak
\Par
For any ${x\in\Cn\!}$, ${\lg\in\Zln}\!$, ${\sb\in\Z^\ell}\!$, \,let
${\xt(\lg,\sb)\>[\eta\)]\in\Cxl}\!$ be the following point:
\ifMag
$$
\gather
\Lline{
\xt(\lg,\sb)\>[\eta\)]\,=\,(\)p^{s_1\lsym+s_{\lg_1}}\eta^{1-\lg_1}x_1,
\)p^{s_2\lsym+s_{\lg_1}}\eta^{\)2-\lg_1}x_1\lc\)p^{s_{\lg_1}}x_1,}
\\
\nn4>
\Rline{p^{s_{\lg_1\]+1}\lsym+s_{\lg_1\]+\lg_2}}\eta^{1-\lg_2}x_2\lc\)
p^{s_{\lg_1\]+\lg_2}}x_2,\,\ldots\,,\)
p^{s_{\ell-\lg_n\]+1}\lsym+s_\ell}\eta^{1-\lg_n}x_n\lc\)p^{s_\ell}x_n)\,.}
\endgather
$$
\else
$$
\align
\xt(\lg,\sb)\>[\eta\)]\,=\,({}\) & p^{s_1\lsym+s_{\lg_1}}\eta^{1-\lg_1}x_1,
\)p^{s_2\lsym+s_{\lg_1}}\eta^{\)2-\lg_1}x_1\lc\)p^{s_{\lg_1}}x_1,
\\
\nn4>
& p^{s_{\lg_1\]+1}\lsym+s_{\lg_1\]+\lg_2}}\eta^{1-\lg_2}x_2\lc\)
p^{s_{\lg_1\]+\lg_2}}x_2,\,\ldots\,,\)
p^{s_{\ell-\lg_n\]+1}\lsym+s_\ell}\eta^{1-\lg_n}x_n\lc\)p^{s_\ell}x_n)\,.
\endalign
$$
\fi
For instance, if $\sb=(0\lc 0)$, then
$\xt(\lg,\sb)\>[\eta\)]=\xt\lg\>[\eta\)]$, \cf. \(xtry).
\Prop{AJ}
Let the parameters $\eta$, $\xn$, $\yn$ be generic. Let a \fn/ ${f(\tell)}$
have the form \(fPT), \(TA). Assume that
\ifMag
$$
|\)p^n A\}\pron\!x_m\1|\,<\>\min\>(1,|\eta\)|^{1-\ell})\,.
$$
\goodbm\nt
\else
$\,{|\)p^n A\}\pron\!x_m\1|<\min\>(1,|\eta\)|^{1-\ell})}$.
\fi
Then the sum below is convergent and
$$
{1\over\ell\)!}\IH[x;y;\eta\)](f\)\Pht)\,=
\smiZ\]\sum_{\,\,\sb\in\Zp^\ell\!\!\!}
\Res\bigl(t_1\1\!\ldots t_\ell\1f(t)\>\Pht(t;x;y;\eta)\bigr)
\vst{t=\xt(\mg,\sb)[\eta]}\,.
$$
Similarly, if $\,{|\)p^M\]A\}\pron\!y_m\1|>\max\>(1,|\eta\)|^{\ell-1})}$,
\,then the sum below is convergent and
$$
{1\over\ell\)!}\IH[x;y;\eta\)](f\)\Pht)\,=\,(-1)^\ell\smiZ\]
\sum_{\,\,\sb\in\Zn^\ell\!\!\!}\Res
\bigl(t_1\1\!\ldots t_\ell\1f(t)\>\Pht(t;x;y;\eta)\bigr)
\vst{t=\yt(\mg,\sb)[\eta\1]}\,.
$$
\endpro
\nt
The proof is similar to the proof of Theorem F.1 in \Cite{TV1}.
The sums in Proposition \[AJ] coincide with the \sym/ \Atype/ Jackson
integrals, see for example \Cite{AK}.

\ifosaka
\vsk1.5>
\cline{\hbox to.3\hsize{\hrulefill}}
\vsk>
\fi

\myRefs
\ifosaka\osakaRefs\else\widest{WW}\fi

\ref\Key A
\by \Aomoto/
\paper $q$-analogue of de~Rham cohomology associated with Jackson integrals, I
\jour Proceedings of Japan Acad.{} \vol 66 {\rm Ser\&A} \yr 1990
\pages 161--164
\moreref \paper II \jour Proceedings of Japan Acad.{} \vol 66 {\rm Ser\&A}
\yr 1990 \pages 240--244
\endref

\ref\Key AK
\by \Aomoto/ and Y\]\&Kato
\paper Gauss decomposition of connection matrices for \sym/ \Atype/ Jackson
integrals \jour \SMNS/ \vol 1 \yr 1995 \issue 4 \pages 623--666
\endref

\ref\Key CM
\by K\&Cho and K\&Matsumoto
\paper Intersection theory for twisted cohomologies and twisted Riemann's
period relations I
\jour Nagoya Math.\ J. \vol 139 \yr 1995 \pages 67--86
\endref

\ref\Key FR
\by \Fre/ and \Reshy/
\paper Quantum affine algebras and holonomic \dif/ \eq/s
\jour \CMP/ \vol 146 \yr 1992 \pages 1--60
\endref

\ref\Key GR
\by G\&Gasper and M\&Rahman
\book Basic \hgeom/ series \bookinfo Encycl.\ Math.\ Appl.{}
\yr 1990 \publ \CUP/
\endref

\ref\Key KBI
\by \Kor/, N\&M\&Bogolyubov and A\&G\&Izergin
\book Quantum inverse scattering method and correlation \fn/s
\yr 1993 \publ \CUP/ 
\endref

\ref\Key M
\by K\&Matsumoto
\paper Intersection numbers of logarithmic \,$k$ forms
\jour Preprint \yr 1996 \pages 1--22
\endref

\ref\Key MV
\by E\&Mukhin and \Varch/
\paper The quantized \KZv/ \eq/ in tensor products of \irr/ \^{$sl_2$-}modules
\jour Preprint \yr 1997 \pages 1--32
\endref

\ref\Key S
\by \Smirnov/
\paper On the deformation of Abelian integrals
\jour \LMP/ \vol 36 \yr 1996 \pages 267--275    
\endref

\ref\Key T
\by \VoT/
\paper Irreducible monodromy matrices for the \Rm/ of the {\sl XXZ}-model
and lattice local quantum Hamiltonians
\jour \TMP/ \vol 63 \yr 1985 \pages 440--454
\endref

\ref\Key TV1
\by \VT/ and \Varch/
\paper Geometry of $q$-\hgeom/ \fn/s, \qaff/s and \eqg/s
\jour Ast\'erisque \vol 246 \yr 1997 \pages 1--135
\endref

\ref\Key TV2
\by \VT/ and \Varch/
\paper Geometry of \qhgeom/ \fn/s as a bridge between Yangians and \qaff/s
\jour \Inv/ \yr 1997 \vol 128 \issue 3 \pages 501--588
\endref

\ref\Key TV3
\by \VT/ and \Varch/
\paper Asymptotic \sol/ to the quantized \KZv/ \eq/ and \Bv/s
\jour \AMS/ Transl.,\ Ser\&\)2 \vol 174 \yr 1996 \pages 235--273
\endref

\ref\Key V
\by \Varch/
\paper Quantized \KZv/ \eq/s, quantum \YB/, and \deq/s for $q$-\hgeom/ \fn/s
\jour \CMP/ \vol 162 \yr 1994 \pages 499--528
\endref

\endRefs

\ifosaka
\vsk1.5>
\rline{\vtop{\hbox{\SPb/ Branch of}
\hbox{Steklov Mathematical Institute}
\hbox{\SPb/ \,191011, \,Russia}
\hbox{\sl E-mail\/{\rm: }\homemail/}
\vsk.5>
\hbox{\it and}
\vsk.5>
\hbox{Department of Mathematics}
\hbox{Faculty of Science, Osaka University}
\hbox{Toyonaka, Osaka 560, \,Japan}}}
\fi

\bye